\documentclass[12pt,a4paper]{article}
\usepackage{apacite}
\usepackage[leqno]{amsmath}
\usepackage{amssymb,amsthm,natbib,graphicx,mathrsfs,tikz,amsfonts,mathdots}
\usepackage{multirow,enumitem}
\usepackage[colorlinks,citecolor=blue,urlcolor=blue]{hyperref}

\usetikzlibrary[decorations.pathreplacing]
\usetikzlibrary{patterns}

\makeatletter
\newcommand{\leqnomode}{\tagsleft@true}
\newcommand{\reqnomode}{\tagsleft@false}
\makeatother

\newcommand{\R}{\mathbb{R}}

\newcommand{\vertiii}[1]{{\left\vert\kern-0.25ex\left\vert\kern-0.25ex\left\vert #1 
    \right\vert\kern-0.25ex\right\vert\kern-0.25ex\right\vert}}

\begin{document}

\title{The Identification Problem for Linear Rational Expectations Models\footnote{Thanks are due to three anonymous referees, Alexander Meyer-Gohde, Manfred Deistler, Benedikt P\"{o}tscher, Hashem Pesaran, Fabio Canova, Barbara Rossi, Christian Brownlees, Geert Mesters, Davide Debortoli, G\'{a}bor Lugosi, Omiros Papaspiliopoulos, Albert Satorra, Bernd Funovits, Peter C.\ B.\ Phillips, Michele Piffer, Ivana Komunjer, and to seminar participants at Universitat Pompeu Fabra, University of Vienna, UC San Diego, University of Southampton, Universitat de les Illes Balears, NYU Abu Dhabi, London School of Economics, Queen Mary University of London, University of Helsinki, University of Kent, University of Cambridge, University of York, and Georgetown University.}}
\author{Majid M.\ Al-Sadoon\footnote{MMA acknowledges support by Spanish Ministry of Economy and Competitiveness projects ECO2012-33247 and ECO2015-68136-P (MINECO/FEDER, UE) and Fundaci\'{o}n BBVA Scientific Research Grant PR16-DAT-0043.}\\Durham University Business School \and Piotr Zwiernik\\University of Toronto, Universitat Pompeu Fabra, \& BSE}

\theoremstyle{plain}
\newtheorem{thm}{Theorem}
\newtheorem{defn}{Definition}
\newtheorem{lem}{Lemma}
\newtheorem{cor}{Corollary}
\newtheorem{prop}{Proposition}

\theoremstyle{definition}
\newtheorem{exmp}{Example}

\renewcommand{\labelenumi}{(\roman{enumi})}

\maketitle

\thispagestyle{empty}

\newpage

\topskip0pt
\vspace*{\fill}
\abstract{
We consider the problem of the identification of stationary solutions to linear rational expectations models from second moments of observable data. Observational equivalence and structural identification are characterized both algebraically and geometrically. Necessary and sufficient criteria are provided for identification and local identification under general non-linear restrictions. The results directly generalize classical theory for vector autoregressive moving average models.

\bigskip


\bigskip

}
\vspace*{\fill}

\thispagestyle{empty}

\newpage

\section{Introduction}
The linear rational expectations model (LREM) is distinguished among dynamic linear systems in that the present state depends not only on events leading up to the present but also on endogenously formulated expectations of the future. Such models are used today by researchers, practitioners, and policy-makers for causal and counter-factual analysis as well as forecasting. Yet the statistical properties of LREMs are poorly understood and identification in particular has remained an open problem throughout what is now known as the rational expectations revolution. 

Establishing identifiability of this model class is important for several reasons. From the parametrization point of view, the parameters of LREMs codify the decision making of economic actors (e.g.\ households or firms) and it is important to know whether or not this behaviour can be learned from the available data. From the estimation point of view, lack of identifiability leads to ill-conditioned optimization procedures when employing extremum estimators; Bayesian methods are not immune to identification failure either, as the posterior retains the shape of the prior along observationally equivalent directions in the parameter space. Finally, inference is substantially more difficult in the absence of identification in both the frequentist and Bayesian perspectives.

Important early contributions addressings the identification problem of LREMs were \cite{muthestimation}, \cite{wallis80}, \cite{rayner}, and \citep{pesaran81,pesaran1987} (see also, \cite{shiller} and \cite{wickens}). However, these results apply to very restricted LREMs and cannot be employed in the study of modern LREMs. Subsequent econometric work in the area completely ignored identification until \cite{canovasala} called attention to serious identification problems plaguing many LREMs used in practice, a point echoed by \cite{pesaransmith}, \cite{romer}, \cite{blanchard}, and \cite{mcdonaldshalizi}. This spurred a number of researchers to provide computational diagnostics for identification \citep{iskrev,kn,qt,kk,kk23}. Unfortunately, this work has not attempted an analytical examination of the mapping from parameters to observables of LREMs. Instead, it treats LREMs as non-linearly parametrized state space models. Consequently, this literature has failed to uncover the underlying reasons for identification failure. This literature has also relied on rather strong assumptions, like regularity (i.e.\ equal number of shocks and observables) and minimality (i.e. non-redundant dynamics). Finally, this literature has not drawn strong connections to classical identification results for linear dynamic models such as \cite{hannan69,hannan71}, \cite{deistler75,deistler76,deistler78,deistler83}, \cite{ds79}, \cite{sargan}, and \cite{deistlerwang} (see \cite{hd} and \cite{hsiao} for complete reviews of the classical literature).

The present work builds on recent results by \cite{onatski}, \cite{andersonetal2012}, \cite{andersonetal}, and \cite{linsys,spectral} to explain why identification failure occurs in LREMs under very weak assumptions (existence, uniqueness, and invertibility), provide a characterization of observational equivalence, and provide analytical diagnostic tests of identification that extend classical results for vector autoregressive moving average (VARMA) models. The key idea is that the mapping from parameters to spectral densities of LREMs involves an initial Wiener-Hopf factorization but is otherwise identical to the mapping for VARMA. Once this is recognized, the theory proceeds almost exactly analogously to the classical theory.

LREMs are subject to identification failure for some of the same reasons that simultaneous equations models and VARMA models are. This is to be expected since the class of LREMs nests the aforementioned classes of models. However, there is a new source of identification failure that afflicts only LREMs. Expectations in LREMs are endogenous. This necessitates more restrictions than might otherwise be called for in classical models. The core of the identification problem of LREMs is therefore that forward dependence is not identified without imposing additional restrictions. This result has been known to many authors in the literature, although it receives its most general treatment in this paper.

This paper provides algebraic and geometric characterizations of classes of observationally equivalent parameters. The dimensions of these classes, which give the number of restrictions necessary for identifiability, is determined for both specific and generic points in the parameter space. The paper then considers identifiability and local identifiability by general non-linear restrictions. Geometrically, every point in the parameter space can be thought of as falling inside the intersection of the set of observationally equivalent parameters and the set of parameters satisfying the given restrictions. When these two sets intersect at a single point, that point is identified. In certain parameter spaces that we consider and for certain restrictions, the intersection may include a subspace so that identification and local identification fail simultaneously.

We should note, finally, that ours is not the only paper to buck the recent trend of deriving computational diagnostics for the identification of LREMs in favor of analytical understanding. The recent paper of \cite{zadrozny} considers the same problem albeit using matrix polynomial methods rather than the Wiener-Hopf method. The two approaches should be viewed as complementary, as many open problems remain (see the conclusion to this paper for a list) and the solutions to these problems are likely to require diverse approaches.

This paper is organized as follows. Section \ref{sec:prelim} sets the notation. Section \ref{sec:model} introduces the parameter spaces we consider. Section \ref{sec:core} discusses the core of the identification problem of LREMs. Section \ref{sec:oe} simplifies the problem of observational equivalence. Section \ref{sec:structural} provides conditions for structural identifiability under general non-linear restrictions on the parameter space. Section \ref{sec:conclusion} concludes. Appendix \ref{sec:topology} considers topological aspects of our parameter spaces. Appendix \ref{sec:proofs} contains the proofs of our results.

\section{Notation}\label{sec:prelim}
Denote by $\mathbb{Z}\subset\mathbb{R}\subset\mathbb{C}$ the sets of integers, real numbers, and complex numbers respectively. We will need $\mathbb{T}=\{z\in\mathbb{C}: |z|=1\}$, $\mathbb{D}=\{z\in\mathbb{C}: |z|<1\}$, and $\overline{\mathbb{D}}=\{z\in\mathbb{C}: |z|\leq1\}$, the unit circle, the open unit disk, and the closed unit disk respectively. Complements of sets will be denoted by a superscript $c$. When forming arrays populated by elements of a given set, we will use the superscript $n\times m$ (e.g.\ $\mathbb{R}^{n\times m}$ is the set of $n\times m$ real matrices); when $m=1$, we will simply use the superscript $n$ (e.g.\ $\mathbb{R}^n$ is the set of $n$-dimensional real column vectors). We denote the orthogonal group of dimension $m$ by $O_m=\{D\in\mathbb{R}^{m\times m}: D'D=I_m\}$. For $D\in\mathbb{R}^{n\times m}$, define $\mathrm{vec}_{\mathbb{R}^{n\times m}}(D)=(D_{11},\ldots,D_{n1},D_{12},\ldots,D_{n2},\ldots,D_{1m},\ldots,D_{nm})'$ and $\|D\|=\left(\sum_{i=1}^n\sum_{j=1}^m |D_{ij}|^2\right)^{1/2}$. For $d\in\mathbb{R}^{nm}$ define $\mathrm{mat}_{\mathbb{R}^{n\times m}}(d)=\left[\begin{smallmatrix} d_1 & d_{n+1} & \cdots & d_{(m-1)n+1}\\ \vdots & \vdots & & \vdots \\ d_n & d_{2n} & \cdots & d_{nm}\end{smallmatrix}\right]$. Note that $\mathrm{vec}_{\mathbb{R}^{n\times m}}(\mathrm{mat}_{\mathbb{R}^{n\times m}}(d))=d$ and $\mathrm{mat}_{\mathbb{R}^{n\times m}}(\mathrm{vec}_{\mathbb{R}^{n\times m}}(D)=D)$ for all $d\in\mathbb{R}^{nm}$ and $D\in\mathbb{R}^{n\times m}$. We will denote by $I_n$ and $0_{n\times m}$ the $n\times n$ identity matrix and the $n\times m$ matrix of zeros respectively.

Denote by $\mathbb{R}[z]\subset\mathbb{R}[z,z^{-1}]\subset\mathbb{R}(z)$ the sets of real polynomials in $z$, real Laurent polynomials in $z$, and real rational functions in $z$ respectively. Denote by $\mathbb{R}[z^{-1}]$ the set of real polynomials in $z^{-1}$. For a non-zero $D\in\mathbb{R}[z,z^{-1}]^{n\times m}$, we denote by $\max\deg(D)$ the highest power of $z$ that appears in $D$, while $\min\deg(D)$ is the lowest power of $z$ that appears in $D$. For $D\in\mathbb{R}(z)^{n\times m}$, we define $D^\ast(z)=D(z^{-1})'\in\mathbb{R}(z)^{m\times n}$. A pair $(D_1,D_2)\in\mathbb{R}[z]^{n\times n}\times\mathbb{R}[z]^{n\times m}$ is said to be left prime if $\mathrm{rank}[\,D_1(z)\; D_2(z)\,]=n$ for all $z\in\mathbb{C}$. Every $D\in\mathbb{R}(z)^{n\times m}$ can be represented as a ratio $D_1^{-1}D_2$, where $(D_1,D_2)\in\mathbb{R}[z]^{n\times n}\times\mathbb{R}[z]^{n\times m}$ are left prime; the number $\delta(D)=\max\deg(\det(D_1))$ is known as the McMillan degree of $D$; the zeros of $\det(D_1)$ are called the poles of $D$; the points at which $\mathrm{rank}(D_2(z))\leq \max_{z\in\mathbb{C}}\mathrm{rank}(D_2(z))$ are called the zeros of $D$; and none of these definitions depend on the particular choice of left prime ratio representation of $D$ \cite[Section 6.5]{kailath}.

Finally, for each $D\in\mathbb{R}(z)^{n\times m}$ with no poles on $\mathbb{T}$, there is an absolutely convergent Laurent series expansion $D(z)=\sum_{i=-\infty}^\infty D_iz^i$ in a neighbourhood of $\mathbb{T}$. The operator whose action on such functions is $D\mapsto[D]_+=\sum_{i=0}^\infty D_iz^i$ is a well defined linear mapping that can be computed by subtracting terms with poles in $\mathbb{D}$ from the partial fractions representation of $D$ \cite[Appendix A]{hansensargent}. Similar remarks apply to the operator $[D]_-=D-[D]_+$.

\section{The Parameter Space}\label{sec:model}
In order to motivate the parameter spaces that we will study, we will need to start with a more general one.

\begin{defn}\label{defn:lrem}
Given positive integers $n$ and $m$ and non-negative integers $\kappa$ and $\lambda$, the set of formal LREMs of input dimension $m$, output dimension $n$, maximum lag $\kappa$, and maximum horizon $\lambda$ is
\begin{multline*}
\Omega_{n,m,\kappa,\lambda}=\Bigg\{(B,A)\in\mathbb{R}[z,z^{-1}]^{n\times n}\times\mathbb{R}[z]^{n\times m}:\\ \max\deg(B)\leq\kappa,\quad\max\deg(A)\leq\kappa,\quad\min\deg(B)\geq-\lambda\Bigg\}.
\end{multline*}
An element $(B,A)=\left(\sum_{i=-\lambda}^\kappa B_iz^i,\sum_{j=0}^\kappa A_jz^j\right)\in\Omega_{n,m,\kappa,\lambda}$ is expressed formally as
\begin{multline}
B_{-\lambda}E_tY_{t+\lambda}+\cdots+B_{-1}E_tY_{t+1}+B_0Y_t+B_1Y_{t-1}+\cdots+B_\kappa Y_{t-\kappa}\\
=A_0\varepsilon_t+A_1\varepsilon_{t-1}+\cdots+A_\kappa\varepsilon_{t-\kappa}, \qquad t\in\mathbb{Z},\label{eq:lremformal}
\end{multline}
where $\varepsilon$ is the input and $Y$ is the output. We will think of $\Omega_{n,m,\kappa,\lambda}$ as a real vector space of dimension $n^2(\kappa+\lambda+1)+nm(\kappa+1)$ and endow it with the Euclidean metric,
\begin{align*}
d((B,A),(\tilde{B},\tilde{A}))=\left(\sum_{i=-\lambda}^\kappa\|B_i-\tilde{B}_i\|^2+\sum_{i=0}^\kappa\|A_i-\tilde{A}_i\|^2\right)^{1/2}.
\end{align*}
\end{defn}
Expression \eqref{eq:lremformal} is understood as a formal set of structural equations relating current, expected, and lagged values of the output $Y$ to current and lagged values of the input $\varepsilon$. It is formal because we have not yet defined what the input process is, what expectations are, or how a unique output is meant to satisfy \eqref{eq:lremformal}. The class of models considered in Definition \ref{defn:lrem} includes formal VARMA models with a given maximum lag length $\kappa$, $\Omega_{n,m,\kappa,0}$ and formal simultaneous equations models $\Omega_{n,m,0,0}$. While the space $\Omega_{n,m,\kappa,\lambda}$ is clearly larger than necessary (e.g.\ the useless parameter $(0_{n\times m},0_{n\times m})$ belongs to $\Omega_{n,m,\kappa,\lambda}$), it serves as an ambient space in which we study particular LREMs and its topology is quite natural.

\begin{exmp}\label{exmp:1}
\cite{hansensargent} study an LREM of the form 
\begin{align*}
\theta_1 E_tY_{t+1}+\theta_2 Y_t+Y_{t-1}=\theta_3\varepsilon_t,\qquad t\in\mathbb{Z}.
\end{align*}
Here, $Y_t$ is the level of employment of a factor of production at time $t$, while $\varepsilon_t$ is real factor rental. The parameters $(\theta_1,\theta_2,\theta_3)$ are required to fall inside a subset of $\mathbb{R}^3$ that ensures economic interpretability as well as existence and uniqueness of a solution. The purpose of this and similar LREMs is to explain the behaviour of $Y_t$ in terms of $\varepsilon_t$, $\varepsilon_{t-1}$, \ldots. The present model can be expressed succinctly using Definition \ref{defn:lrem} as the pair
\begin{align*}
B(\theta)&=\theta_1z^{-1}+\theta_2+z,& A(\theta)&=\theta_3,
\end{align*}
which is clearly an element of $\Omega_{1,1,1,1}$.
\end{exmp}

The empirical macroeconomics literature (e.g.\ \cite{canova}, \cite{dd}, and \cite{hs}) takes the input process in Definition \ref{defn:lrem} to be i.i.d.\ with mean zero and variance $I_m$. This is stronger than necessary for our formulation of LREM solutions.

\begin{defn}\label{defn:sol}
We say that a solution exists for LREM $(B,A)\in\Omega_{n,m,\kappa,\lambda}$ if there is a sequence $(C_i)_{i\in\mathbb{Z}, i\geq0}$ in $\mathbb{R}^{n\times m}$ such that $\sum_{i=0}^\infty \|C_i\|^2<\infty$ and for any $m$-dimensional stochastic process $(\varepsilon_t)_{t\in\mathbb{Z}}$ satisfying $\mathbb{E}(\varepsilon_{t+1}|\varepsilon_t,\varepsilon_{t-1},\ldots)=0_{m\times 1}$ and $\mathbb{E}(\varepsilon_t\varepsilon_t')=I_m$ for all $t\in\mathbb{Z}$, the linear process $(Y_t)_{t\in\mathbb{Z}}$ defined by
\begin{align}
Y_t=\sum_{i=0}^\infty C_i\varepsilon_{t-i},\qquad t\in\mathbb{Z},\label{eq:ma}
\end{align}
satisfies
\begin{align}
\sum_{i=-\lambda}^\kappa B_i\mathbb{E}(Y_{t-i}|\varepsilon_t,\varepsilon_{t-1},\ldots)=\sum_{i=0}^\kappa A_i\varepsilon_{t-i},\qquad t\in\mathbb{Z}.\label{eq:lrem}
\end{align}
The solution to $(B,A)$ is unique if whenever $(Y_t)_{t\in\mathbb{Z}}$ and $(\tilde{Y}_t)_{t\in\mathbb{Z}}$ are solutions generated by the same input process $(\varepsilon_t)_{t\in\mathbb{Z}}$, then $P(Y_t=\tilde{Y}_t)=1$ for all $t\in\mathbb{Z}$.
\end{defn}

The formulation of existence and uniqueness in Definition \ref{defn:sol} can be viewed as a special case of the definitions in \cite{onatski} and \cite{spectral} even though these works use best linear predictors in place of conditional expectations in \eqref{eq:lrem}. That is because, with the inputs restricted to white noise martingale difference processes (we do not need strict stationarity or even identical distributedness), the conditional expectations in \eqref{eq:lrem} are also best linear predictors \cite[Theorem 1.4.2]{hd}.

When an LREM has a unique solution, its transfer function has a simple algebraic structure.

\begin{defn}\label{defn:S}
Define the set of real stable rational transfer functions as
\begin{align*}
\Sigma=\left\{D\in\mathbb{R}(z): D\text{ has no poles in }\overline{\mathbb{D}}\right\}.
\end{align*}
\end{defn}

Every element of $D\in\Sigma^{n\times m}$ has a representation $D(z)=\sum_{i=0}^\infty D_iz^i$, converging absolutely in a neighbourhood of $\overline{\mathbb{D}}$.

\begin{thm}[Onatski's Criterion]\label{thm:onatski}
Let $(B,A)\in\Omega_{n,m,\kappa,\lambda}$. Then $(B,A)$ has a unique solution if and only if there exists a unique $C\in \Sigma^{n\times m}$ such that
\begin{align}\label{eq:rationaltoeplitz}
[BC]_+=A.
\end{align}
The impulse responses of the linear process \eqref{eq:ma} are the coefficients of the expansion $C(z)=\sum_{i=0}^\infty C_iz^i$ in a neighbourhood of $\overline{\mathbb{D}}$.
\end{thm}

\cite{onatski} derived \eqref{eq:rationaltoeplitz} in a more general setting, where $B$ is an element of a Wiener algebra, the space of solutions is Hardy space, and $[\;\cdot\;]_+$ is substituted by an orthogonal projection operator (see his Equation (39)). Note that in the VARMA case, where $\lambda=0$, equation \eqref{eq:rationaltoeplitz} reduces to the familiar equation $BC=A$ \cite[Theorem 1.2.1]{hd}. This is because $B\in\mathbb{R}[z]^{n\times n}\subseteq\Sigma^{n\times n}$ and $C\in\Sigma^{n\times m}$, so $BC\in\Sigma^{n\times m}$ and therefore $[BC]_+=BC$. Equation \eqref{eq:rationaltoeplitz} will play a crucial role in our development below.

In order to ensure existence and uniqueness of a solution to $(B,A)\in\Omega_{n,m,0,0}$ it is clearly necessary and sufficient for $B$ to be non-singular \cite[p.\ 12]{dhrymesbook}. The VARMA literature, on the other hand, imposes the sufficient but not necessary stability condition: for $(B,A)\in\Omega_{n,m,\kappa,0}$, $\det(B(z))\neq0$ for $|z|\leq1$ \cite[Equation (1.2)]{andersonetal}. For example, $(z,z)\in\Omega_{1,1,1,0}$ has a unique solution but does not satisfy the aforementioned condition. In this work, we will employ a direct generalization of the VARMA condition.

\begin{defn}\label{defn:wh}
For a given $D\in\mathbb{R}[z,z^{-1}]^{n\times n}$, we say that $D$ has a canonical Wiener-Hopf factorization if
\begin{align}
\begin{aligned}
&D=D_-D_+;\\
&D_-\in\mathbb{R}[z^{-1}]^{n\times n},\; \mathrm{rank}(D_-(z))=n \text{ for all }z\in\mathbb{D}^c,\text{ and }\lim_{z\rightarrow\infty}D_-(z)=I_n;\\
&\text{and }D_+\in\mathbb{R}[z]^{n\times n} \text{ and } \mathrm{rank}(D_+(z))=n \text{ for all } z\in\overline{\mathbb{D}}.
\end{aligned}\label{eq:wh}
\end{align}
Define $\Omega_{n,m,\kappa,\lambda}^{EU}$ to be the set of pairs $(B,A)\in\Omega_{n,m,\kappa,\lambda}$ such that $B$ satisfies conditions \eqref{eq:wh}.
\end{defn}

The canonical Wiener-Hopf factorization defined in Definition \ref{defn:wh} is unique whenever it exists \cite[Theorems I.1.1 and I.1.2]{cg}. Like the stability condition for VARMA, the restriction imposed in Definition \ref{defn:wh} concerns only the $B$ part of the LREM parameter pair. While it is rather technical, it has a very simple meaning: $B$ satisfies the conditions in Definition \ref{defn:wh} if and only if $(B,\tilde{A})\in\Omega_{n,m,\kappa,\lambda}$ has a unique solution for every possible choice of $\tilde{A}$ \cite[Theorem 3]{spectral}; that is, if and only if \eqref{eq:lremformal} has a unique solution regardless of the specification of the right hand side. If $(B,A)\in\Omega_{n,m,\kappa,\lambda}^{EU}$, then $B_-$ can be thought of as determining the forward dependence of $Y$ on expectations of $\varepsilon$, while $B_+$ can be thought of as determining the backwards dependence of $Y$ on $\varepsilon$. The unique solution to $(B,A)$ has a transfer function $C$ given by Onatski's formula \cite[p.\ 342]{onatski},
\begin{align}
C=B_+^{-1}[B_-^{-1}A]_+.\label{eq:onatski}
\end{align}
It is an easy exercise to check that $(B_+,[B_-^{-1}A]_+)\in\Omega_{n,m,\kappa,0}^{EU}$ whenever $(B,A)\in\Omega_{n,m,\kappa,\lambda}^{EU}$. It is also easy to check that $\Omega_{n,m,0,0}^{EU}$ is precisely the set of $(B,A)\in\Omega_{n,m,0,0}$ such that $\det(B)\neq0$ and $\Omega_{n,m,\kappa,0}^{EU}$ is precisely the set of $(B,A)\in\Omega_{n,m,\kappa,0}$ such that $\det(B(z))\neq0$ for $|z|\leq1$ (one takes $B_-=I_n$ and $B_+=B$). Thus, the restriction in Definition \ref{defn:wh} is a direct generalization of conditions in the literature. Finally, we mention that $\Omega_{n,m,\kappa,\lambda}^{EU}$ is an open subset of $\Omega_{n,m,\kappa,\lambda}$ (see Theorem \ref{thm:topology} of Appendix \ref{sec:topology}), a fact that will be important for later results.

\begin{exmp}\label{exmp:2}
Consider $(B,A)\in\Omega_{1,1,1,1}^{EU}$. Then $B=B_-B_+$ with $B_-(z)=1-b_-z^{-1}$, $B_+(z)=b_0(1-b_+ z)$, $b_0\neq0$, and $|b_\pm|<1$. Note that in this parameterization $B_{-1}=-b_-b_0$, $B_0=b_0(1+b_-b_+)$, and $B_1=-b_0b_+$. Thus, the parameter displays forward dependence of $Y$ on expectations of $\varepsilon$ if and only if $b_-\neq0$ and backward dependence of $Y$ on $\varepsilon$ if and only if $b_+\neq0$. Writing $A(z)=a_0+a_1z$, we have
\begin{align*}
B_-(z)^{-1}A(z)=\left(\sum_{i=0}^\infty b_-^iz^{-i}\right)(a_0+a_1 z)
\end{align*}
for $z\in\mathbb{D}^c$. This implies that
\begin{align*}
\left[B_-^{-1}A\right]_+(z)=a_1 b_-+a_0+a_1 z.
\end{align*}
In consequence, the solution to $(B,A)$ has transfer function
\begin{align*}
C(z)&=\frac{a_1 b_-+a_0+a_1z}{b_0(1-b_+z)},
\end{align*}
which is the same as the solution of the ARMA model with autoregressive part $B_+(z)=b_0(1-b_+z)$ and the moving average part $[B_-^{-1}A]_+(z)=a_1b_-+a_0+a_1z$.
\end{exmp}

The rest of the restrictions we need to impose on $\Omega_{n,m,\kappa,\lambda}$ are standard in the literature and address the identifiability of the transfer function of the linear process \eqref{eq:ma} from its spectral density. Supposing, as the literature does, that the input $(\varepsilon_t)_{t\in\mathbb{Z}}$ is unobserved and that inference is conducted using only the second moments of $(Y_t)_{t\in\mathbb{Z}}$ (see e.g.\ \cite{iskrev}, \cite{qt} \cite{kn}, \cite{kk,kk23}), then the identification problem is concerned with the injectivity of the mapping that assigns to each element of $\Omega_{n,m,\kappa,\lambda}^{EU}$ the spectral density function of its solution. Since solutions to elements of $\Omega_{n,m,\kappa,\lambda}^{EU}$ always have real stable rational transfer functions by Theorem \ref{thm:onatski}, we can just as well consider autocovariance generating functions.

\begin{defn}\label{defn:oe}
Denote the space of autocovariance generating functions of solutions to elements of $\Omega_{n,m,\kappa,\lambda}^{EU}$ by 
\begin{align*}
\Delta_{n,m}=\{DD^\ast: D\in\Sigma^{n\times m}\}.
\end{align*}
Let $\varphi:\Omega_{n,m,\kappa,\lambda}^{EU}\rightarrow\Delta_{n,m}$ be the mapping
\begin{align*}
\varphi:(B,A)\mapsto CC^\ast,
\end{align*}
where $C$ is the transfer function of the solution to $(B,A)$. We say that the $(B,A)$ and $(\tilde{B},\tilde{A})$ in $\Omega_{n,m,\kappa,\lambda}^{EU}$ are observationally equivalent if $\varphi(B,A)=\varphi(\tilde{B},\tilde{A})$ and we express this using the equivalence relation $(B,A)\sim(\tilde{B},\tilde{A})$. We denote by $(B,A)/\sim$ the set of all $(\tilde{B}, \tilde{A})\in\Omega_{n,m,\kappa,\lambda}^{EU}$ such that $(\tilde{B},\tilde{A})\sim (B,A)$.
\end{defn}

The linear systems literature has traditionally addressed the identification problem by factorizing $\varphi$ as $\varphi=\varphi_1\circ\varphi_2$, where
\begin{align*}
\Omega_{n,m,\kappa,\lambda}^{EU}\overset{\varphi_2}{\longrightarrow} \Sigma^{n\times m}\overset{\varphi_1}{\longrightarrow}\Delta_{n,m}.
\end{align*}
The identification problem can then be broken down into two parts: one part concerning the injectivity of $\varphi_1$ and another concerning the injectivity of $\varphi_2$. For example, \cite{ds} consider $\varphi$ restricted to $\Omega_{n,n,\kappa,0}^{EU}$ in Section 6.1, $\varphi_1^{-1}$ (the spectral factorization problem over $\Delta_{n,n}$) in Section 6.2, and $\varphi_2^{-1}$ (the identifiability of elements of $\Omega_{n,n,\kappa,0}^{EU}$ from transfer functions) in Section 6.3. We will follow the same approach, although we will require a more delicate factorization of $\varphi$ because we allow for singular spectral densities and because solving LREMs is more complicated than solving VARMA models.

Consider the spectral factorization problem first. The problem consists of describing the set of all $\tilde{C}\in \Sigma^{n\times m}$ such that $\tilde{C}\tilde{C}^\ast=CC^\ast$ for a given $C\in \Sigma^{n\times m}$. If $\mathrm{rank}(C)<m$ uniformly over $\mathbb{C}$ (including if $n<m$), then one of the shocks (more precisely, a linear combination of $\varepsilon$) plays no role in determining $Y$; such redundancy is likely to be eliminated in the modelling process and we will rule it out here. It is easy to show then (see Theorem 2 (4) of \cite{baggioferrante}) that the set of all factors of $CC^\ast$ is $C\cdot\{V\in \Sigma^{m\times m}: V^\ast V=I_m \text{ on }\mathbb{T}\}$. \cite{ds} refer to the latter set as the set of Blaschke matrix functions (they are also known as all-pass, inner, and para-unitary matrix functions). In order to make progress on the identification problem, we need to introduce further restrictions that reduce the set of possible factors to a more manageable set. In the time series literature \citep{rosenblatt,hd,bd,ds}, this is accomplished by imposing the following restriction.

\begin{defn}\label{defn:invertible}
If $C\in \Sigma^{n\times m}$ and
\begin{align*}
\mathrm{rank}(C(z))=m\text{ for all }z\in\mathbb{D} \text{ (resp.\ $z\in\overline{\mathbb{D}}$)},
\end{align*}
we say that $C$ is invertible (strictly invertible). Let $\Omega_{n,m,\kappa,\lambda}^{EUI}$ be the subset of $\Omega_{n,m,\kappa,\lambda}^{EU}$ whose solutions have transfer functions that are invertible.
\end{defn}

Theorem 2 (3) of \cite{baggioferrante} proves that the set of invertible factors of $CC^\ast$ for $C\in \Sigma^{n\times m}$ is $C\cdot O_m$. We will be able to work with this set but we can go further (analytically at least) if we reduce this set even further. If $C$ is invertible, then the first impulse response $C(0)$ must be of rank $m$. By applying the Gram-Schmidt algorithm to the rows of $C(0)$ from top to bottom, we can reduce it to the following form.

\begin{defn}\label{defn:cqlt}
A matrix $D\in\mathbb{R}^{n\times m}$ of rank $m$ is said to be in canonical quasi-lower triangular (CQLT) form if the first non-zero element of column $j$ is positive and occurs in row $i_j$ with $1\leq i_1<\cdots<i_m\leq n$. 
\end{defn}

For example, in the case of $n=4$ and $m=3$, the following matrices are CQLT
\begin{align*}
\left[\begin{array}{cccc}
+ & 0 & 0 \\
\ast & + & 0 \\
\ast & \ast & + \\
\ast & \ast & \ast
\end{array}\right]
,\qquad\left[\begin{array}{cccc}
0 & 0 & 0 \\
+ & 0 & 0 \\
\ast & + & 0 \\
\ast & \ast & + 
\end{array}\right]
,\qquad\left[\begin{array}{cccc}
+ & 0 & 0 \\
\ast & 0 & 0 \\
\ast & + & 0 \\
\ast & \ast & +
\end{array}\right],\qquad
\left[\begin{array}{cccc}
+ & 0 & 0 \\
\ast & + & 0 \\
\ast & \ast & 0 \\
\ast & \ast & +
\end{array}\right],
\end{align*}
where $+$ represents a positive number and $\ast$ represents an arbitrary number. It is easy to show that for any $X\in\mathbb{R}^{n\times m}$ of rank $m$ and any $V\in O_m$, $X$ and $XV$ have the same CQLT form. Thus, if $C\in\Sigma^{n\times m}$ is invertible and $C(0)$ is CQLT, then $C$ can be retrieved from $CC^\ast$. Note that in the special case $\Omega_{n,n,\kappa,0}^{EUI}$, the CQLT restriction on $C(0)$ is simply the recursive identification restriction common in the structrual VAR literature (we reconsider structural VAR models in Section \ref{sec:structural}). A similar restriction is used in \cite{andersonetal2012} although they do not impose the positivity restriction on the first non-zero element of each column, which means they identify each column only up to a sign.

\begin{defn}\label{defn:0}
Let $\Omega_{n,m,\kappa,\lambda}^{EUI0}$ be the subset of $\Omega_{n,m,\kappa,\lambda}^{EUI}$ whose solutions have transfer functions with first impulse response in CQLT form.
\end{defn}

The additional restrictions imposed on $\Omega_{n,m,\kappa,\lambda}^{EUI}$ and $\Omega_{n,m,\kappa,\lambda}^{EUI0}$ over and above $\Omega_{n,m,\kappa,\lambda}^{EU}$ may seem peculiar in that they are not imposed directly on the parameter $(B,A)$ but rather indirectly on the transfer function of the solution $C$. Indeed, Section 3 of \cite{andersonetal} provides a parametrization of $\Omega_{n,m,\kappa,0}^{EUI}$   and Section 12.6 of \cite{lutkepohl} (implicitly) provides a parametrization of $\Omega_{n,n,\kappa,0}^{EUI0}$ by imposing restrictions directly on $(B,A)$ in both cases. It is easy to check that 
\begin{align*}
\Omega_{n,m,\kappa,\lambda}^{EUI}&=\left\{(B,A)\in\Omega_{n,m,\kappa,\lambda}^{EU}: [B_-^{-1}A]_+\text{ is invertible}\right\},\\
\Omega_{n,m,\kappa,\lambda}^{EUI0}&=\left\{(B,A)\in\Omega_{n,m,\kappa,\lambda}^{EUI}: B_+(0)^{-1}[B_-^{-1}A]_+(0)\text{ is CQLT}\right\}.
\end{align*}
In particular, $(B_+,[B_-^{-1}A]_+)\in\Omega_{n,m,\kappa,0}^{EUI}$ (resp.\ $\Omega_{n,m,\kappa,0}^{EUI0}$) whenever $(B,A)\in\Omega_{n,m,\kappa,\lambda}^{EUI}$ (resp.\ $\Omega_{n,m,\kappa,\lambda}^{EUI0}$). However, these equivalent formulations have neither simplicity nor convenience to recommend them.

It is important to highlight that the CQLT assumption is interesting primarily from the theoretical standpoint because it allows the transfer function of the solution to be identifiable from its autocovariance generating function (Zadrozny's Assumption 5 serves the same purpose \citep{zadrozny}). LREMs used in practice may not satisfy the CQLT assumption. However, (i) the analytical results we are able to obtain under this assumption are noteworthy and (ii) the CQLT restriction in univariate LREMs will usually be innocuous as the researcher is likely to know a priori the sign of the first impulse response; if positive, this is already guaranteed by the CQLT restriction and, if negative, one can work with $(-\varepsilon_t)_{t\in\mathbb{Z}}$ instead of $(\varepsilon_t)_{t\in\mathbb{Z}}$. In Example \ref{exmp:1}, \cite{hansensargent} presume $\theta$ to belong to a subset of $\mathbb{R}^3$ such that the first impulse response is negative (an unexpected increase in real factor rental will reduce the employment of that factor); in that case, we can simply modify the model to $(\theta_1 z^{-1}+\theta_2+z,-\theta_3)$.

Combining Theorem \ref{thm:onatski} with the classical results on observational equivalence reviewed above, we arrive at the following.

\begin{thm}\label{thm:oe}
\begin{enumerate}
\item If $(B,A),(\tilde B,\tilde A)\in\Omega_{n,m,\kappa,\lambda}^{EUI}$ and the solution to $(B,A)$ has transfer function $C$, then $(B,A)\sim(\tilde B,\tilde A)$ if and only if
\begin{align*}
[\tilde{B}C]_+V=\tilde{A}
\end{align*}
for some $V\in O_m$.
\item If $(B,A),(\tilde B,\tilde A)\in\Omega_{n,m,\kappa,\lambda}^{EUI0}$ and the solution to $(B,A)$ has transfer function $C$, then $(B,A)\sim(\tilde B,\tilde A)$ if and only if
\begin{align*}
[\tilde{B}C]_+=\tilde{A}.
\end{align*}
\end{enumerate}
\end{thm}

We mention in closing that Appendix \ref{sec:topology} provides topological characterizations of $\Omega_{n,m,\kappa,\lambda}^{EU}$, $\Omega_{n,m,\kappa,\lambda}^{EUI}$, and $\Omega_{n,m,\kappa,\lambda}^{EUI0}$.

\section{Core of the Identification Problem for LREMs}\label{sec:core}

In order to understand how identification failure happens in LREMs, it will be helpful to develop the content of Theorem \ref{thm:oe} step by step. This is because identification problems accumulate with increased complexity of the parameter space under consideration and it helps to see where each type of identification problem arises.

Consider the case of simultaneous equations models $(B,A),(\tilde{B},\tilde{A})\in\Omega_{n,m,0,0}^{EUI}$. By Theorem \ref{thm:oe} (i), $(\tilde{B},\tilde{A})\sim(B,A)$ if and only if $\tilde A=\tilde{B}(B^{-1}A)V$ for some $V\in O_m$. This is equivalent to the existence of an invertible $G\in\mathbb{R}^{n\times n}$ and a $V\in O_m$ such that $(\tilde{B},\tilde{A})=G(B,AV)$. Clearly, $G=\tilde{B}B^{-1}$.

Now consider VARMA models $(B,A),(\tilde{B},\tilde{A})\in\Omega_{n,m,\kappa,0}^{EUI}$. Theorem \ref{thm:oe} (i) implies that $(\tilde{B},\tilde{A})\sim(B,A)$ if and only if $\tilde A=\tilde{B}(B^{-1}A)V$ for some $V\in O_m$. This then is equivalent to the existence of an invertible $G\in\Sigma^{n\times n}$ and a $V\in O_m$ such that $(\tilde{B},\tilde{A})=G(B,AV)$. Again $G=\tilde{B}B^{-1}$.

Finally, consider the case of LREMs $(B,A),(\tilde{B},\tilde{A})\in\Omega_{n,m,\kappa,\lambda}^{EUI}$. Theorem \ref{thm:oe} (i) together with Onatski's formula \eqref{eq:onatski} implies that $(\tilde{B},\tilde{A})\sim(B,A)$ if and only if for some $V\in O_m$,
\begin{align*}
\tilde A=[\tilde{B}B_+^{-1}[B_-^{-1}AV]_+]_+.
\end{align*}
To make sense of this equation, let us borrow some notation from \cite{spectral}.

\begin{defn}\label{defn:toeplitz}
For $M\in\mathbb{R}(z)$ with no poles on $\mathbb{T}$, define the associated operator
\begin{align*}
\boldsymbol{M}:\Sigma\rightarrow\Sigma,\qquad\boldsymbol{M}:D\mapsto[MD]_+.
\end{align*}
For $M\in\mathbb{R}(z)^{k\times l}$ with no poles on $\mathbb{T}$, define
\begin{align*}
\boldsymbol{M}:\Sigma^{k\times l}\rightarrow\Sigma^{k\times l},\qquad\boldsymbol{M}:D\mapsto[MD]_+,
\qquad\boldsymbol{M}=\left[\begin{array}{ccc}
\boldsymbol{M}_{11} & \cdots & \boldsymbol{M}_{1l}\\
\vdots & & \vdots\\
\boldsymbol{M}_{k1} & \cdots & \boldsymbol{M}_{kl}
\end{array}\right].
\end{align*}
\end{defn}

The operators introduced in Definition \ref{defn:toeplitz} are known as Toeplitz operators \citep{nikolski} and are very well-understood in functional analysis \citep{gf,cg}.

Then if $(B,A)\in\Omega_{n,m,\kappa,\lambda}^{EUI}$, the operator associated with $B$, $\boldsymbol{B}$, is invertible and its inverse is given by
\begin{align*}
\boldsymbol{B}^{-1}:D\mapsto B_+^{-1}[B_-^{-1}D]_+.
\end{align*}
We can now rewrite the criterion of Theorem \ref{thm:oe} (i) as 
\begin{align*}
\boldsymbol{\tilde{A}}=\boldsymbol{\tilde{B}}\boldsymbol{B}^{-1}\boldsymbol{A}\boldsymbol{V}.
\end{align*}
Note that $(B,A)$ is retrievable from $(\boldsymbol{B},\boldsymbol{A})$ by the fact that $B=z^{-\lambda}\boldsymbol{B}(z^\lambda I_n)$ and $A=\boldsymbol{A}(I_m)$. We have therefore proven the following fact about observational equivalence in $\Omega_{n,m,\kappa,\lambda}^{EUI}$ (and by extension $\Omega_{n,m,\kappa,\lambda}^{EUI0}$).

\begin{thm}\label{thm:G}
\begin{enumerate}
\item If $(B,A),(\tilde B,\tilde A)\in\Omega_{n,m,\kappa,\lambda}^{EUI}$, then $(B,A)\sim(\tilde B,\tilde A)$ if and only if there exists a linear mapping $\boldsymbol{G}:\Sigma^{n\times (n+m)}\rightarrow\Sigma^{n\times (n+m)}$ and a $V\in O_m$ such that
\begin{align*}
(\boldsymbol{\tilde{B}},\boldsymbol{\tilde{A}})=\boldsymbol{G}(\boldsymbol{B},\boldsymbol{A}\boldsymbol{V}).
\end{align*}
Clearly, $\boldsymbol{G}=\boldsymbol{\tilde{B}}\boldsymbol{B}^{-1}$.
\item If $(B,A),(\tilde B,\tilde A)\in\Omega_{n,m,\kappa,\lambda}^{EUI0}$, then $(B,A)\sim(\tilde B,\tilde A)$ if and only if there exists a linear mapping $\boldsymbol{G}:\Sigma^{n\times (n+m)}\rightarrow\Sigma^{n\times (n+m)}$ such that 
\begin{align*}
(\boldsymbol{\tilde{B}},\boldsymbol{\tilde{A}})=\boldsymbol{G}(\boldsymbol{B},\boldsymbol{A}).
\end{align*}
Clearly, $\boldsymbol{G}=\boldsymbol{\tilde{B}}\boldsymbol{B}^{-1}$.
\end{enumerate}
\end{thm}

Theorem \ref{thm:G} makes clear that identification failure in $\Omega_{n,m,\kappa,\lambda}^{EUI}$ occurs because of the modulo-$O_m$ uniqueness of spectral factorization and because $\boldsymbol{B}^{-1}\boldsymbol{A}=(\boldsymbol{G}\boldsymbol{B})^{-1}(\boldsymbol{G}\boldsymbol{A})$ for certain operators $\boldsymbol{G}$ such that $\boldsymbol{G}(\boldsymbol{B},\boldsymbol{A})$ is associated with an element of $\Omega_{n,m,\kappa,\lambda}^{EUI}$. Thus, Theorem \ref{thm:G} shows that the identification problem for LREMs is just an extension of that for VARMA and simultaneous equations models, albeit with slightly more advanced mathematics. Indeed, by specializing Theorem \ref{thm:G} to the cases $\lambda=0$ and $\lambda=\kappa=0$ we obtain the results reviewed above for VARMA and simultaneous equations models respectively.

While Theorem \ref{thm:G} makes clear why identification failure happens in LREMs, it does not clarify how. In order to see how identification failure occurs in LREMs, it is helpful to trace step-by-step the process of obtaining $CC^\ast$ from $(B,A)\in\Omega_{n,m,\kappa,\lambda}^{EUI}$, where $C$ is the transfer function of the solution to $(B,A)$. We have already obtained one factorization of this mapping, $\varphi=\varphi_1\circ\varphi_2$, in Section \ref{sec:model}, where we saw that $\varphi_1$ introduced the problem of identifiability up to right multiplication by $O_m$. We now explore the identifiablility problems that arise in $\varphi_2$ by breaking up the process of obtaining $C$ from $(B,A)$ into three intermediate steps: the first step obtains the canonical Wiener-Hopf factorization $B=B_-B_+$ and the numerator matrix polynomial in Onatski's formula $[B_-^{-1}A]_+$; the second step discards $B_-$ as it is no longer needed; and the third step divides $[B_-^{-1}A]_+$ by $B_+$. Thus, $C=\varphi_2(B,A)=\varphi_{2a}\circ\varphi_{2b}\circ\varphi_{2c}(B,A)$, where
\begin{align}
\begin{aligned}
\varphi_{2c}&:\Omega_{n,m,\kappa,\lambda}^{EUI}\rightarrow\Omega_{n,m,\kappa,\lambda}^{EUIWH}, &(\mathsf{X},\mathsf{Y})&\overset{\varphi_{2c}}{\longmapsto}(\mathsf{X}_-,\mathsf{X}_+,[\mathsf{X}_-^{-1}\mathsf{Y}]_+),\\
\varphi_{2b}&:\Omega_{n,m,\kappa,\lambda}^{EUIWH}\rightarrow \Omega_{n,m,\kappa,0}^{EUI}, &(\mathsf{Z},\mathsf{X},\mathsf{Y})&\overset{\varphi_{2b}}{\longmapsto}(\mathsf{X},\mathsf{Y}),\\
\varphi_{2a}&:\Omega_{n,m,\kappa,0}^{EUI}\rightarrow\Sigma^{n\times m}, &(\mathsf{X},\mathsf{Y})&\overset{\varphi_{2a}}{\longmapsto}\mathsf{X}^{-1}\mathsf{Y},
\end{aligned}
\label{eq:sequence}
\end{align}
and
\begin{align*}
\Omega_{n,m,\kappa,\lambda}^{EUIWH}=\Big\{&(\mathsf{Z},\mathsf{X},\mathsf{Y})\in\mathbb{R}[z^{-1}]^{n\times n}\times\in\mathbb{R}[z]^{n\times n}\times\in\mathbb{R}[z]^{n\times m}:\\
&\min\deg(\mathsf{Z})\geq-\lambda,\; \mathrm{rank}(\mathsf{Z}(z))=n \text{ for all }z\in\mathbb{D}^c,\text{ and }\lim_{z\rightarrow\infty}\mathsf{Z}(z)=I_n;\\
&\max\deg(\mathsf{X})\leq\kappa \text{ and } \mathrm{rank}(\mathsf{X}(z))=n \text{ for all } z\in\overline{\mathbb{D}};\\
&\text{ and }\mathsf{Y}\text{ is invertible }\Big\}.
\end{align*}
It is an easy exercise to verify that all the mappings above are onto. Note that $\varphi_{2b}\circ\varphi_{2c}$ restricted to $\Omega_{n,m,\kappa,0}^{EUI}$ is just the identity map and $\varphi_{2a}$ is the only common solution step for both VARMA and LREMs. The $\varphi_{2a}$ mapping is precisely where we encounter the identification failures of VARMA because a given transfer function may be associated with many different VARMA parameters. Thus, whatever identification problems are particular to LREMs must arise in $\varphi_{2b}\circ\varphi_{2c}$. Let us consider each step separately. It is easy to check that $\varphi_{2c}^{L}:(\mathsf{Z},\mathsf{X},\mathsf{Y})\mapsto(\mathsf{Z}\mathsf{X},[\mathsf{Z}\mathsf{Y}]_+)$ is a left inverse of $\varphi_{2c}$ and therefore $\varphi_{2c}$ is one-to-one and cannot introduce identification problems. This leaves $\varphi_{2b}$ as the only remaining suspect. Since $\varphi_{2b}$ drops the first argument of $\varphi_{2c}(\mathsf{X},\mathsf{Y})=(\mathsf{X}_-,\mathsf{X}_+,[\mathsf{X}_-^{-1}\mathsf{Y}]_+)$ and $\mathsf{X}_-$ determines the forward dependence of the solution to $(\mathsf{X},\mathsf{Y})\in\Omega_{n,m,\kappa,\lambda}^{EUI}$, we arrive at the core of the identification problem for LREMs, that forward dependence is not identified without imposing additional restrictions. 

In the next section, we provide a complete characterization of the fibers of $\varphi$ (i.e.\ the sets of observational equivalence) in $\Omega_{n,m,\kappa,\lambda}^{EUI}$ and $\Omega_{n,m,\kappa,\lambda}^{EUI0}$.

\section{Observational Equivalence}\label{sec:oe}
Fix $(B,A)\in\Omega_{n,m,\kappa,\lambda}^{EUI}$ and let its solution have transfer function $C$. Then Theorem \ref{thm:oe} (i) implies that $\big((B,A)/\sim\big)\cap\;\Omega_{n,m,\kappa,\lambda}^{EUI}$ is the set of all $(\tilde B,\tilde A)\in\Omega_{n,m,\kappa,\lambda}^{EUI}$ satisfying
\begin{align*}
[\tilde{B}C]_+=\tilde{A}V'
\end{align*}
for some $V\in O_m$. Obtaining the Taylor series expansion about $z=0$ of both sides and matching coefficients we see that this condition can be written equivalently as
\begin{align*}
\left[
\begin{array}{cccccccc}
\tilde B_{-\lambda} & \cdots & \tilde B_\kappa & 0_{n\times n} & \cdots
\end{array}
\right]
\left[
\begin{array}{cccccccc}
C_\lambda & C_{\lambda+1} & C_{\lambda+2} & \cdots \\
C_{\lambda-1} & C_\lambda & C_{\lambda+1} & \cdots \\
\vdots & \ddots & \ddots & \ddots \\
\vdots & \ddots & \ddots & \ddots \\
C_1 & C_2 & C_3 & \ddots \\
C_0 & C_1 & C_2 & \ddots \\
0_{n\times m} & C_0 & C_1  & \ddots\\
0_{n\times m} & 0_{n\times m} & C_0 & \ddots \\
\vdots & \ddots & \ddots & \ddots
\end{array}
\right]=\left[\begin{array}{ccccc}
\tilde A_0 V' & \cdots & \tilde A_\kappa V' & 0_{n\times m} & \cdots 
\end{array}
\right].
\end{align*}
It appears then, that $\big((B,A)/\sim\big)\cap\;\Omega_{n,m,\kappa,\lambda}^{EUI}$ is determined by an infinite number of equations. However, exploiting the rationality of $\tilde{B}$, $\tilde{A}$, and $C$, we can truncate the infinite system of equations to a finite one. In order to do that, we will need some additional notation.

\begin{defn}\label{defn:vecmat}
The vectorization operator on $\Omega_{n,m,\kappa,\lambda}$ is given by
\begin{multline*}
\mathrm{vec}_{\Omega_{n,m,\kappa,\lambda}}\left(\sum_{i=-\lambda}^\kappa \mathsf{X}_iz^i,\sum_{i=0}^\kappa \mathsf{Y}_iz^i\right)\\
=(\mathrm{vec}_{\mathbb{R}^{n\times n}}(\mathsf{X}_{-\lambda})',\ldots,\mathrm{vec}_{\mathbb{R}^{n\times n}}(\mathsf{X}_{\kappa})',\mathrm{vec}_{\mathbb{R}^{n\times m}}(\mathsf{Y}_{0})',\ldots,\mathrm{vec}_{\mathbb{R}^{n\times m}}(\mathsf{Y}_{\kappa})')',
\end{multline*}
defined for all $\mathsf{X}_{-\lambda},\ldots, \mathsf{X}_\kappa\in\mathbb{R}^{n\times n}$ and all $\mathsf{Y}_0,\ldots, \mathsf{Y}_\kappa\in\mathbb{R}^{n\times m}$. It obtains the vectorization of all of the coefficient matrices of a parameter. The operator $\mathrm{mat}_{\Omega_{n,m,\kappa,\lambda}}$ reconstructs parameters from their vectorizations,
\begin{align*}
\mathrm{mat}_{\Omega_{n,m,\kappa,\lambda}}\left(\left(\mathsf{x}_{-\lambda}',\ldots,\mathsf{x}_{\kappa}',\mathsf{y}_{0}',\ldots,\mathsf{y}_{\kappa}\right)'\right)=\left(\sum_{i=-\lambda}^\kappa \mathrm{mat}_{\mathbb{R}^{n\times n}}(\mathsf{x}_i)z^i,\sum_{i=0}^\kappa \mathrm{mat}_{\mathbb{R}^{n\times m}}(\mathsf{y}_i)z^i\right),
\end{align*}
defined for all $\mathsf{x}_{-\lambda},\ldots, \mathsf{x}_\kappa\in\mathbb{R}^{n^2}$ and all $\mathsf{y}_0,\ldots, \mathsf{y}_\kappa\in\mathbb{R}^{nm}$.
\end{defn}

\begin{thm}\label{thm:oeEUI}
Let $(B,A)\in\Omega_{n,m,\kappa,\lambda}^{EUI}$, let $C(z)=\sum_{i=0}^\infty C_iz^i$ be the transfer function of the solution to $(B,A)$, let
\begin{gather*}
T=\left[\begin{array}{ccccccc}
C_\lambda & \cdots & C_{\kappa+\lambda} \\
\vdots &  & \vdots \\
C_0 & \cdots & C_\kappa \\
 & \ddots & \vdots \\
0 & & C_0 
\end{array}\right], \qquad
H=\left[\begin{array}{ccccc}
C_{\kappa+\lambda+1} & \cdots & C_{(n+1)\kappa+\lambda} \\
\vdots & & \vdots \\
C_{\kappa+1} & \cdots & C_{(n+1)\kappa} \\
\vdots & & \vdots \\
C_1 & \cdots & C_{n\kappa}
\end{array}\right],\\ 
\intertext{and, for $V\in O_m$, let}
P(V)=\left[\begin{array}{cc}
-T & -H\\	
I_{\kappa+1}\otimes V' & 0_{m(\kappa+1)\times nm\kappa}\\	
\end{array}\right].
\end{gather*}
Then:
\begin{enumerate}
\item If $(\tilde B,\tilde A)\in\Omega_{n,m,\kappa,\lambda}^{EUI}$, then $(\tilde B,\tilde A)\sim(B,A)$ if and only if there exists a $V\in O_m$ such that
\begin{align*}
\mathrm{vec}_{\Omega_{n,m,\kappa,\lambda}}(\tilde{B},\tilde{A})\in\ker\left(P(V)'\otimes I_n\right).
\end{align*}
\item $\big((B,A)/\sim\big)\cap\;\Omega_{n,m,\kappa,\lambda}^{EUI}$ is given by
\begin{align*}
\left(\bigcup_{V\in O_m}\mathrm{mat}_{\Omega_{n,m,\kappa,\lambda}}\left(\ker\left(P(V)'\otimes I_n\right)\right)\right)\cap\Omega_{n,m,\kappa,\lambda}^{EU}.
\end{align*}
\item $\big((B,A)/\sim\big)\cap\;\Omega_{n,m,\kappa,\lambda}^{EUI}$ is a differentiable manifold of dimension
\begin{align*}
\frac{1}{2}m(m-1)+n^2(\kappa+\lambda+1)-n\delta\left(C(z^{-1})-C(0)\right)\geq \frac{1}{2}m(m-1)+n^2(\lambda+1).
\end{align*}
Equality holds generically in $\Omega_{n,m,\kappa,\lambda}^{EUI}$.
\end{enumerate}
\end{thm}

Theorem \ref{thm:oeEUI} (i) reduces the infinite dimensional system \eqref{eq:rationaltoeplitz} to a finite dimensional criterion. In particular, observational equivalence is determined by a non-linear criterion in $(\tilde{B},\tilde{A})$ and a $V\in O_m$ (note that $V'$ multiplies $\tilde{A}$). Theorem \ref{thm:oeEUI} (ii) characterizes a set of observationally equivalent parameters in $\Omega_{n,m,\kappa,\lambda}^{EUI}$ as the intersection of two cones in $\Omega_{n,m,\kappa,\lambda}$, $\left(\bigcup_{V\in O_m}\mathrm{mat}_{\Omega_{n,m,\kappa,\lambda}}\left(\ker\left(P(V)'\otimes I_n\right)\right)\right)$ and $\Omega_{n,m,\kappa,\lambda}^{EU}$ (a cone is a set invariant to multiplication by positive scalars). Finally, Theorem \ref{thm:oeEUI} (iii) gives the dimension of the set of observationally equivalent parameters. Note that cones are textbook examples of geometric objects that are not manifolds because ``no neighborhood of the vertex point looks like a simple piece of the plane'' \cite[p.\ 1]{gp}; however, the vertex point is the useless parameter $(0_{n\times n},0_{n\times m})\not\in\Omega_{n,m,\kappa,\lambda}^{EUI}$ and therefore excluded. The dimensions reported in Theorem \ref{thm:oeEUI} (iii) can be understood as the sum of the dimensions of the fibres of $\varphi_1$ (i.e.\ $\frac{1}{2}m(m-1)$), $\varphi_{2c}$ (i.e.\ $n^2(\kappa+1)-n\delta\left(C(z^{-1})-C(0)\right)$), $\varphi_{2b}$ (i.e.\ $n^2\lambda$), and $\varphi_{2a}$ (i.e.\ $0$). It can also be understood as the number of restrictions that must be imposed on the parameter space in order to identify a parameter. It is given by the effective number of free parameters of the model, $n^2(\kappa+\lambda+1)+nm(\kappa+1)$, plus the dimension of $O_m$ (recall that $O_m$ is a differentiable manifold of dimension $\frac{1}{2}m(m-1)$ \cite[p.\ 22]{gp}), minus the complexity of the transfer function $C$ of the solution to $(B,A)$, as measured by $n\delta(C(z^{-1})-C(0))$. The generic subset of $\Omega_{n,m,\kappa,\lambda}^{EUI}$ identified in Theorem \ref{thm:oeEUI} (iii) is the set of all pairs $(B,A)$ whose solution is strictly invertible, its VARMA representation $(B_+,[B_-^{-1}A]_+)$ is left prime, $\det(B_+)$ has distinct zeros, and $\det(B_\kappa)\neq0$.

\begin{exmp}\label{exmp:3}
If $(B,A)\in\Omega_{n,n,0,0}^{EUI}$ has the solution $C\in\mathbb{R}^{n\times n}$, then $Y_t=C\varepsilon_t$ for $t\in\mathbb{Z}$ and the only information available is the covariance matrix of $Y_t$, $CC'$. Theorem \ref{thm:oe} (i) then implies that if $(\tilde{B},\tilde{A})\in\Omega_{n,n,0,0}^{EUI}$, then $(B,A)\sim(\tilde{B},\tilde{A})$ if and only if $\tilde{A}=\tilde{B}CV$ for some $V\in O_m$ and Theorem \ref{thm:oeEUI} (i) is just a vectorization of that criterion. More interestingly, Theorem \ref{thm:oeEUI} (ii) and (iii) provide a geometric description of $\big((B,A)/\sim\big)\cap\;\Omega_{n,m,\kappa,\lambda}^{EUI}$, whose dimension, $2n^2-\frac{1}{2}n(n+1)$, is exactly the number of necessary restrictions to obtain local identification computed in Section 9.1.3 of \cite{lutkepohl}.
\end{exmp}

\begin{exmp}\label{exmp:4}
Consider $\Omega_{1,1,1,1}^{EUI}$ again. If $(\tilde B,\tilde A)\sim(B,A)$ then, by Theorem \ref{thm:oe} (i),
\begin{align*}
\left[(\tilde B_{-1}z^{-1}+\tilde B_0+\tilde B_1z)C(z)\right]_+V=\tilde A_0+\tilde A_1z,\qquad V^2=1.
\end{align*}
Noting that $[z^{-1}C(z)]_+=z^{-1}(C(z)-C(0))$ and $[z^iC(z)]_+=z^iC(z)$ for $i\geq0$, we may rewrite the equation above as
\begin{align*}
\tilde{B}_{-1}z^{-1}(C(z)-C(0))+\tilde{B}_0C(z)+\tilde{B}_1zC(z)=\tilde{A}_0V+\tilde{A}_1Vz,\qquad V^2=1.
\end{align*}
Equating the first three Taylor series coefficients of both sides we obtain
\begin{align*}
\left[\begin{array}{cccccc}
-C_1 & -C_0 & 0   & V & 0 \\
-C_2 & -C_1 & -C_0 & 0 & V \\
-C_3 & -C_2 & -C_1 & 0 & 0 
\end{array}\right]\left[\begin{array}{c}
\tilde{B}_{-1} \\ \tilde{B}_0 \\ \tilde{B}_1 \\ \tilde{A}_0 \\ \tilde{A}_1
\end{array}\right]&=0_{3\times 1}, & &V^2=1.
\end{align*}
This is the content of Theorem \ref{thm:oeEUI} (i). To see that, indeed, these are the only equations that determine observational equivalence, we can make use of the fact that $C=\frac{a_1 b_-+a_0+a_1z}{b_0(1-b_+z)}$, obtained in Example \ref{exmp:2}. If we multiply both sides of $[\tilde{B}C]_+V=\tilde{A}$ by $b_0(1-b_+z)$, we obtain an equality of two polynomials of degree 2; equating the coefficients and rearranging, we obtain exactly the equations above (this computation is clearly tedious so we omit it).

When $(B,A)=(1,1)$ so that $C_0=1$ and $C_1=C_2=C_3=0$. Then $(1,1)/\sim$ is a subset of
\begin{align*}
\mathrm{mat}_{\Omega_{1,1,1,1}}\left(\{(\mathsf{x}_1,\mathsf{x}_2,\mathsf{x}_3,\mathsf{x}_1,\mathsf{x}_3):\mathsf{x}\in\mathbb{R}^3\}\cup\{(\mathsf{x}_1,\mathsf{x}_2,\mathsf{x}_3,-\mathsf{x}_1,-\mathsf{x}_3):\mathsf{x}\in\mathbb{R}^3\}\right).
\end{align*}
Not every element of this set is in $\Omega_{1,1,1,1}^{EUI}$ (e.g.\ the useless parameter $(0,0)\not\in\Omega_{1,1,1,1}^{EUI}$); however, its intersection with $\Omega_{1,1,1,1}^{EU}$ is precisely $\big((B,A)/\sim\big)\cap\;\Omega_{1,1,1,1}^{EUI}$ according to Theorem \ref{thm:oeEUI} (ii). For example, $(1,1)\sim\left(1+\frac{1}{2}z,-1-\frac{1}{2}z\right)\in\Omega_{1,1,1,1}^{EUI}$ owing to the fact that $\left(0,1,\frac{1}{2},-1,-\frac{1}{2}\right)'\in \ker(P(-1)'\otimes I_n)$; this illustrates the identification failure familiar for ARMA models due to the two polynomials not being left prime \citep[p.\ 123]{ds}. A more interesting example in our context is $\left(\frac{1}{3}z^{-1}+1+\frac{1}{2}z,1+\frac{1}{2}z\right)\in\Omega_{n,m,\kappa,\lambda}^{EUI}$, which is observationally equivalent to $(1,1)$ because $\left(\frac{1}{3},1,\frac{1}{2},1,\frac{1}{2}\right)'\in \ker(P(+1)'\otimes I_n)$.

We may also consider more complicated parameters. Take Example \ref{exmp:1} with $\theta=\left(\frac{2}{3}, -\frac{7}{3},-1\right)$. Then $(B(\theta),A(\theta))=\left(\frac{2}{3}z^{-1}-\frac{7}{3}+z,-1\right)$ and following the computations in Example \ref{exmp:2}, we have that $C(z)=\frac{1}{2-z}$ so that $C_i=2^{-1-i}$ for $i\geq0$. Then $\left(\frac{2}{3}z^{-1}-\frac{7}{3}+z,-1\right)/\sim$ is a subset of
\begin{multline*}
\mathrm{mat}_{\Omega_{1,1,1,1}}\bigg(\bigg\{\bigg(\mathsf{x}_1,\mathsf{x}_2,-\frac{1}{4}\mathsf{x}_1-\frac{1}{2}\mathsf{x}_2,\frac{1}{4}\mathsf{x}_1+\frac{1}{2}\mathsf{x}_2,0\bigg):\mathsf{x}\in\mathbb{R}^2\bigg\}\cup\\
\bigg\{\bigg(\mathsf{x}_1,\mathsf{x}_2,-\frac{1}{4}\mathsf{x}_1-\frac{1}{2}\mathsf{x}_2,-\frac{1}{4}\mathsf{x}_1-\frac{1}{2}\mathsf{x}_2,0\bigg):\mathsf{x}\in\mathbb{R}^2\bigg\}
\bigg).
\end{multline*}
Note the lower dimension of the two subspaces. The parameter $\left(\frac{2}{3}z^{-1}-\frac{7}{3}+z,-1\right)$ is more complex as an element of $\Omega_{1,1,1,1}^{EUI}$ than $(1,1)$ and therefore requires fewer restrictions in order to identify it. In contrast, $(1,1)$ is a simple model in a large space with significantly more observationally equivalent parameters.
\end{exmp}

In the case of $(B,A)\in\Omega_{n,m,\kappa,\lambda}^{EUI0}$, observational equivalence is much simpler because the transfer function of the solution is identified from the associated autocovariance generating function.

\begin{thm}\label{thm:oeEUI0}
Let $(B,A)\in\Omega_{n,m,\kappa,\lambda}^{EUI0}$, let $C(z)=\sum_{i=0}^\infty C_iz^i$ be the transfer function of the solution to $(B,A)$ and let $P$ be as in Theorem \ref{thm:oeEUI}. Then:
\begin{enumerate}
\item If $(\tilde B,\tilde A)\in\Omega_{n,m,\kappa,\lambda}^{EUI0}$, then $(\tilde B,\tilde A)\sim(B,A)$ if and only if
\begin{align*}
\mathrm{vec}_{\Omega_{n,m,\kappa,\lambda}}(\tilde{B},\tilde{A})\in\ker\left(P(I_m)'\otimes I_n\right).
\end{align*}
\item $\big((B,A)/\sim\big)\cap\;\Omega_{n,m,\kappa,\lambda}^{EUI0}$ is given by
\begin{align*}
\mathrm{mat}_{\Omega_{n,m,\kappa,\lambda}}\left(\ker\left(P(I_m)'\otimes I_n\right)\right)\cap\Omega_{n,m,\kappa,\lambda}^{EU}.
\end{align*}
\item $\big((B,A)/\sim\big)\cap\;\Omega_{n,m,\kappa,\lambda}^{EUI0}$ is a differentiable manifold of dimension
\begin{align*}
n^2(\kappa+\lambda+1)-n\delta\left(C(z^{-1})-C(0)\right)\geq n^2(\lambda+1).
\end{align*}
Equality holds generically in $\Omega_{n,m,\kappa,\lambda}^{EUI0}$.
\end{enumerate}
\end{thm}

Observational equivalence in $\Omega_{n,m,\kappa,\lambda}^{EUI0}$ is much more straightforward. In particular, if $(B,A)\in\Omega_{n,m,\kappa,\lambda}^{EUI0}$, every $(\tilde{B},\tilde{A})\in\Omega_{n,m,\kappa,\lambda}^{EUI0}$ observationally equivalent to $(B,A)$ has the same transfer function. This makes the observational equivalence criterion $[\tilde{B}C]_+=\tilde{A}$ linear in the coefficients of $(\tilde{B},\tilde{A})$ and its finite dimensional counterpart, developed in Theorem \ref{thm:oeEUI0} (i), is also linear. Theorem \ref{thm:oeEUI0} (ii) characterizes a subspace that contains $\big((B,A)/\sim\big)\cap\;\Omega_{n,m,\kappa,\lambda}^{EUI0}$ and this subspace is a subset of the cone characterized in Theorem \ref{thm:oeEUI} (ii). As before, $\big((B,A)/\sim\big)\cap\;\Omega_{n,m,\kappa,\lambda}^{EUI0}$ is relatively open and so small perturbations of $(B,A)$ along the aforementioned subspace do not lead to departure from $\Omega_{n,m,\kappa,\lambda}^{EUI0}$. The reductions in the dimensions given in Theorem \ref{thm:oeEUI0} (iii) relative to Theorem \ref{thm:oeEUI} (iii) are due to the normalization of $C_0$ to CQLT form. The generic subset of $\Omega_{n,m,\kappa,\lambda}^{EUI0}$ identified in Theorem \ref{thm:oeEUI0} (iii) is, again, the set of all pairs $(B,A)$ whose solution is strictly invertible, its VARMA representation $(B_+,[B_-^{-1}A]_+)$ is left prime, $\det(B_+)$ has distinct zeros, and $\det(B_\kappa)\neq0$.

Theorems \ref{thm:oeEUI} and \ref{thm:oeEUI0} make clear how severe the identification problem in LREMs actually is. The generic VARMA model in $\Omega_{n,m,\kappa,0}^{EUI}$ requires $\frac{1}{2}m(m-1)+n^2$ restrictions but the generic LREM in $\Omega_{n,m,\kappa,\lambda}^{EUI}$ requires an additional $n^2\lambda$ restrictions. The additional number of restrictions is exactly the number of additional endogenous variables introduced by the expectation terms $E_tY_{t+i}$ for $i>0$ in \eqref{eq:lremformal} and a facet of the core of the identification problem for LREMs, that forward dependence is not identified without imposing additional restrictions.

Having simplified the characterization of observational equivalence, we are now in a position to study particular models such as the one given in Example \ref{exmp:1} by imposing the necessary additional restrictions. 

\section{Structural Identifiability}\label{sec:structural}
In practice, LREMs are restricted in a variety of ways such as exclusion (setting a parameter to zero), normalization (setting a parameter to 1), and, more generally, restrictions that set functions of the parameters (possibly non-linear and possibly across equations) to fixed values. Here we consider the ability of such restrictions to identify a single parameter.

\begin{defn}\label{defn:ident}
Let $(B,A)\in\Lambda\subset\Omega_{n,m,\kappa,\lambda}^{EU}$. We say that $(B,A)$ is identified in $\Lambda$ if 
\begin{align*}
(B,A)=\big((B,A)/\sim\big)\cap\;\Lambda.
\end{align*}
We say that a parameter $(B,A)$ is locally identified in $\Lambda$ if there is an open set $N\subset\Omega_{n,m,\kappa,\lambda}$ containing $(B,A)$ such that 
\begin{align*}
(B,A)=\big((B,A)/\sim\big)\cap\;\Lambda\cap N.
\end{align*}
\end{defn}

The basic idea in resolving identifiability queries is then to solve for the intersection $\big((B,A)/\sim\!\big)\cap\;\Lambda$ in Definition \ref{defn:ident}. In particular, the identification results below stack equations restricting the parameters to $\Lambda$ along with the equations characterizing observational equivalence from Theorem \ref{thm:oeEUI} (ii) (or Theorem \ref{thm:oeEUI0} (ii)) then verify whether or not $(B,A)$ is the unique solution. Local identification is verified by local uniqueness of the solution $(B,A)$ to these systems of equations. The results below define $\Lambda$ implicitly in terms of the inverse of some given function $R$ but the general principle applies just as well if the restrictions are parametrized as $\Lambda=\{(B(\theta),A(\theta))\in\Omega_{n,m,\kappa,\lambda}: \theta\in\Theta\}$ for some $\Theta\subset\mathbb{R}^d$. In that case, we may substitute $(B(\tilde\theta),A(\tilde\theta))$ into the equations characterizing observational equivalence from Theorem \ref{thm:oeEUI} (ii) (or Theorem \ref{thm:oeEUI0} (ii)) and solve for $\tilde\theta$; if $\tilde\theta=\theta$ is the only solution in $\Theta$ such that $(B(\tilde \theta),A(\tilde \theta))\in\Omega_{n,m,\kappa,\lambda}^{EUI}$ (or $\Omega_{n,m,\kappa,\lambda}^{EUI0}$), then $(B(\theta),A(\theta))$ is identified.

\begin{thm}\label{thm:identEUI}
Let $R:\mathbb{R}^{n^2(\kappa+\lambda+1)+nm(\kappa+1)}\rightarrow \mathbb{R}^r$ and let
\begin{align}
\Omega^R=\left\{(\mathsf{X},\mathsf{Y})\in\Omega_{n,m,\kappa,\lambda}: R\left(\mathrm{vec}_{\Omega_{n,m,\kappa,\lambda}}(\mathsf{X},\mathsf{Y})\right)=0_{r\times 1}\right\}.\label{eq:nonlinear}
\end{align}
If $(B,A)\in\Omega^R\cap\Omega_{n,m,\kappa,\lambda}^{EUI}$ and $P$ is defined as in Theorem \ref{thm:oeEUI}, then 
$(B,A)$ is identified in $\Omega^R\cap\Omega_{n,m,\kappa,\lambda}^{EUI}$ if and only if
\begin{align}
\begin{aligned}
\left(P(\mathsf{V})'\otimes I_n\right)
\mathrm{vec}_{\Omega_{n,m,\kappa,\lambda}}(\mathsf{X},\mathsf{Y})&=0_{nm((n+1)\kappa+1)\times 1},\\
R\left(\mathrm{vec}_{\Omega_{n,m,\kappa,\lambda}}(\mathsf{X},\mathsf{Y})\right)&=0_{r\times 1},\\
\mathsf{V}'\mathsf{V}&=I_m,
\end{aligned}
\label{eq:identcriterionEUI}
\end{align}
has a unique solution in $\Omega_{n,m,\kappa,\lambda}^{EU}\times\mathbb{R}^{m\times m}$, $\left((\mathsf{X},\mathsf{Y}),\mathsf{V}\right)=((B,A),I_m)$.
\end{thm}

From an algebraic perspective, Theorem \ref{thm:identEUI} (i) reduces the problem of identifiability of a point in $\Omega_{n,m,\kappa,\lambda}^{EUI}$ to the problem of solving the non-linear system of  equations obtained from stacking the observational equivalence equations in Theorem \ref{thm:oeEUI} (i) together with the equations that determine the imposed restrictions in \eqref{eq:nonlinear}. Geometrically, the problem reduces to that of finding the points of intersection of the space $\Omega^R$ with the cone of observationally equivalent models in $\Omega_{n,m,\kappa,\lambda}^{EUI}$. In the special case when $R$ is an affine function, $R\left(\mathrm{vec}_{\Omega_{n,m,\kappa,\lambda}}(\mathsf{X},\mathsf{Y})\right)=R\;\mathrm{vec}_{\Omega_{n,m,\kappa,\lambda}}(\mathsf{X},\mathsf{Y})-u$, for some $R\in \mathbb{R}^{r\times n^2(\kappa+\lambda+1)+nm(\kappa+1)}$ and some $u\in \mathbb{R}^r$, a necessary condition for the point $(B,A)$ in Theorem \ref{thm:identEUI} to be identified is that $\left[\begin{smallmatrix} P(V)'\otimes I_n\\
R \end{smallmatrix}\right]$ be of full rank for all $V\in O_m$. Although there is little more that we can say analytically, this system of polynomial equations is easy to solve numerically because: (i) it only requires the first $(n+1)\kappa+\lambda+1$ impulse responses of the solution to $(B,A)$ as inputs, (ii) it has a simple structure that is easy to code, and (iii) solutions are sought in an open set $\Omega_{n,m,\kappa,\lambda}^{EU}\times\mathbb{R}^{m\times m}$. Indeed, when $R$ is a polynomial, \eqref{eq:identcriterionEUI} is a system of polynomial equations and there are efficient algorithms for obtaining the set of all solutions to such systems (\cite{kk23} have proposed using one such algorithm).

\begin{thm}\label{thm:localEUI}
Let $R$ and $\Omega^R$ be as in Theorem \ref{thm:identEUI} and suppose $R$ is continuously differentiable. If $(B,A)\in\Omega^R\cap\Omega_{n,m,\kappa,\lambda}^{EUI}$ and $P$ is defined as in Theorem \ref{thm:oeEUI}, then 
$(B,A)$ is locally identified in $\Omega^R\cap\Omega_{n,m,\kappa,\lambda}^{EUI}$ if 
\begin{align}
\begin{aligned}
\left(P(I_m)'\otimes I_n\right)\mathrm{vec}_{\Omega_{n,m,\kappa,\lambda}}(\mathsf{dX},\mathsf{dY})+K\mathrm{vec}_{\mathbb{R}^{m\times m}}(\mathsf{dV})&=0_{nm((n+1)\kappa+1)\times 1},\\
\nabla R(B,A)\mathrm{vec}_{\Omega_{n,m,\kappa,\lambda}}(\mathsf{dX},\mathsf{dY})&=0_{r\times 1},\\
\mathsf{dV}'+\mathsf{dV}&=0_{m\times m},
\end{aligned}
\label{eq:nonlinearcriterionEUI}
\end{align}
has a unique solution in $\Omega_{n,m,\kappa,\lambda}\times\mathbb{R}^{m\times m}$, $\left((\mathsf{dX},\mathsf{dY}),\mathsf{dV}\right)=((0_{n\times n},0_{n\times m}),0_{m\times m})$, where
\begin{align*}
K=\left[\begin{array}{c}
I_m\otimes A_0\\
\vdots\\
I_m\otimes A_\kappa\\
0_{n^2m\kappa\times m^2}\\
\end{array}
\right].
\end{align*}
Conversely, $(B,A)$ is not locally identified in $\Omega^R\cap\Omega_{n,m,\kappa,\lambda}^{EUI}$ if the dimension of solutions to \eqref{eq:nonlinearcriterionEUI} is non-zero and constant in a neighborhood of $(B,A)$ in  $\Omega^R\cap\Omega_{n,m,\kappa,\lambda}^{EUI}$.
\end{thm}

Local identifiability is addressed in Theorem \ref{thm:localEUI} by the purely linear criterion \eqref{eq:nonlinearcriterionEUI}. The reason why the converse requires a more stringent condition is well understood in the identification literature \cite[Section 5.1]{hsiao}. The condition is known as regularity. Without it, it is not possible to conclude non-identifiabilty because the intersection of $\Omega^R$ with $\big((B,A)/\sim\big)\cap\;\Omega_{n,m,\kappa,\lambda}^{EUI}$ may be non-transversal. Note that when $m=1$, $\mathsf{dV}=0$ and the criterion is simplified because $O_1=\{+1,-1\}$ is discrete and so any solutions to \eqref{eq:identcriterionEUI} local to $((\mathsf{X},\mathsf{Y}),\mathsf{V})=((B,A),+1)$ must keep $\mathsf{V}$ fixed at $+1$. \cite{darkside} consider the problem of locally identified models that are not identified in a susbset of $\Omega_{n,m,\kappa,0}^{EUI}$.

\begin{exmp}\label{exmp:6}
Let us continue where we left off in Example \ref{exmp:4}. By Theorem~\ref{thm:oeEUI} and the discussion following it, to obtain identifiability we need to impose at least $n^2(\lambda+1)=2$ restrictions. Thus, we may consider whether $(B(\theta),A(\theta))$ with $\theta=\left(\frac{2}{3}, -\frac{7}{3},-1\right)$ is identified in $\Omega_{1,1,1,1}^{EUI}$ by the additional restrictions
\begin{align*}
\Omega^R=\left\{(\mathsf{X},\mathsf{Y})\in\Omega_{1,1,1,1}: \mathsf{X}_1=1, \mathsf{Y}_1=0\right\}.
\end{align*}
This implies that \eqref{eq:identcriterionEUI} can be written as
\begin{align*}
\left[\begin{array}{cccccc}
-\frac{1}{4} & -\frac{1}{2} & 0   &  V & 0 \\
-\frac{1}{8} & -\frac{1}{4} & -\frac{1}{2} &  0 & V \\
-\frac{1}{16} & -\frac{1}{8} & -\frac{1}{4} &  0 & 0 \\
0 & 0 & 1 &  0 & 0 \\
0 & 0 & 0 &  0 & 1
\end{array}\right]\left[\begin{array}{c}
\mathsf{X}_{-1}\\ \mathsf{X}_0 \\ \mathsf{X}_1 \\ \mathsf{Y}_0 \\ \mathsf{Y}_1
\end{array}\right]&=
\left[\begin{array}{c}
0 \\ 0 \\ 0 \\ 1 \\ 0
\end{array}\right],\qquad V^2=1.
\end{align*}
The determinant of the matrix on the left hand side is $0$ so $(B(\theta),A(\theta))$ is not identifiable in $\Omega^R\cap\Omega_{1,1,1,1}^{EUI}$ or $\Omega^R\cap\Omega_{1,1,1,1}^{EUI0}$. It was already evident from the parameterization of its observationally equivalent parameters that such restrictions would not be sufficient. Suppose we modify our restrictions as
\begin{align*}
\Omega^{\tilde{R}}=\left\{(\mathsf{X},\mathsf{Y})\in\Omega_{1,1,1,1}: \mathsf{X}_1=1, \mathsf{Y}_1=0, \mathsf{X}_{-1}=\frac{2}{3}\right\}.
\end{align*}
Then it is easy to check that \eqref{eq:identcriterionEUI} has two solutions 
\begin{align*}
((\mathsf{X},\mathsf{Y}),\mathsf{V})\in\left\{\left(\left(\frac{2}{3}z^{-1} -\frac{7}{3}+z,-1\right),+1\right),\left(\left(\frac{2}{3}z^{-1} -\frac{7}{3}+z,+1\right),-1\right)
\right\}.
\end{align*}
This type of identification failure can be overcome by an additional inequality restriction (e.g.\ $\mathsf{Y}_0\leq0$) or by an additional equality restriction (e.g.\ $\mathsf{Y}_0=-1$). Note that $(B(\theta),A(\theta))$ is nevertheless locally identified in $\Omega^{\tilde{R}}\cap\Omega_{n,m,\kappa,\lambda}^{EUI}$.
\end{exmp}

\begin{thm}\label{thm:identEUI0}
Let $R$ and $\Omega^R$ be as in Theorem \ref{thm:identEUI}. If $(B,A)\in\Omega^R\cap\Omega_{n,m,\kappa,\lambda}^{EUI0}$ and $P$ is defined as in Theorem \ref{thm:oeEUI}, then $(B,A)$ is identified in $\Omega^R\cap\Omega_{n,m,\kappa,\lambda}^{EUI0}$ if and only if
\begin{align}
\begin{aligned}
\left(P(I_m)'\otimes I_n\right)\mathrm{vec}_{\Omega_{n,m,\kappa,\lambda}}(\mathsf{X},\mathsf{Y})&=0_{nm((n+1)\kappa+1)\times 1},\\
R\left(\mathrm{vec}_{\Omega_{n,m,\kappa,\lambda}}(\mathsf{X},\mathsf{Y})\right)&=0_{r\times 1},
\end{aligned}
\label{eq:identcriterionEUI0}
\end{align}
has a unique solution in $\Omega_{n,m,\kappa,\lambda}^{EU}$,  $\left(\mathsf{X},\mathsf{Y}\right)=(B,A)$.
\end{thm}

Theorem \ref{thm:identEUI0} has the same geometric interpretation as Theorem \ref{thm:identEUI}. In the special case when $R$ is an affine function, $R\left(\mathrm{vec}_{\Omega_{n,m,\kappa,\lambda}}(\mathsf{X},\mathsf{Y})\right)=R\;\mathrm{vec}_{\Omega_{n,m,\kappa,\lambda}}(\mathsf{X},\mathsf{Y})-u$, for some $R\in \mathbb{R}^{r\times n^2(\kappa+\lambda+1)+nm(\kappa+1)}$ and some $u\in \mathbb{R}^r$, Theorem \ref{thm:identEUI0} shows that identifiability of a point in $\Omega_{n,m,\kappa,\lambda}^{EUI0}$ is assessed by a linear criterion. Indeed, since system \eqref{eq:identcriterionEUI0} always has the solution $(B,A)$, one can simply check whether or not $\left[\begin{smallmatrix} P(I_m)'\otimes I_n\\ R \end{smallmatrix}\right]$ is of full column rank. Additionally, if identification fails, the intersection set must contain a subspace and Theorem \ref{thm:oeEUI} (ii) then implies that every neighbourhood of $(B,A)$ contains infinitely many observationally equivalent parameters that also satisfy the given restrictions. That is, a point in $\Omega_{n,m,\kappa,\lambda}^{EUI0}$ is identifiable by affine restrictions if and only if it is locally identifiable.

\begin{exmp}\label{exmp:svar}
Consider the parameter space $\Omega_{n,n,\kappa,0}^{EUI0}$ with 
\begin{align*}
\Omega^R=\left\{(\mathsf{X},\mathsf{Y})\in\Omega_{n,n,\kappa,0}: \mathsf{Y}_1=\cdots=\mathsf{Y}_\kappa=0_{n\times n} \text{ and }\mathsf{X}_0=I_n\right\}.
\end{align*}
Substituting the latter equations of \eqref{eq:identcriterionEUI0} into the observational equivalence equation, we obtain
\begin{align*}
[\;\mathsf{Y}_0\;\;0_{n\times n}\;\;\cdots\;\;0_{n\times n}\;]=[\;I_n\;\;\mathsf{X}_1\;\;\mathsf{X}_2\;\;\cdots\;\;\mathsf{X}_\kappa\;][\;T\;\;H\;].
\end{align*}
Since $T$ is upper triangular with non-singular diagonal blocks, $[\;T\;\;H\;]$ is of full row rank. Thus, there exists a unique solution to $\mathsf{Y}_0$, $\mathsf{X}_1$,\ldots, $\mathsf{X}_\kappa$. The set $\Omega^R\cap\Omega_{n,n,\kappa,0}^{EUI0}$ is commonly known as the recursively identified structural VAR model popularized by \cite{macroreality}, (see Section 2.3.2 of \cite{lutkepohl} for an overview).
\end{exmp}

\begin{thm}\label{thm:localEUI0}
Let $R$ and $\Omega^R$ be as in Theorem \ref{thm:localEUI}. If $(B,A)\in\Omega^R\cap\Omega_{n,m,\kappa,\lambda}^{EUI0}$ and $P$ is defined as in Theorem \ref{thm:oeEUI}, then $(B,A)$ is locally identified in $\Omega^R\cap\Omega_{n,m,\kappa,\lambda}^{EUI0}$ if
\begin{align}
\left[\begin{array}{c}
P(I_m)'\otimes I_n\\
\nabla R(B,A)
\end{array}\right]
\mathrm{vec}_{\Omega_{n,m,\kappa,\lambda}}(\mathsf{dX},\mathsf{dY})=0_{r+nm((n+1)\kappa+1)\times 1}\label{eq:nonlinearcriterionEUI0}
\end{align}
has a unique solution in $\Omega_{n,m,\kappa,\lambda}^{EU}$,  $\left(\mathsf{dX},\mathsf{dY}\right)=(0_{n\times n},0_{n\times m})$.

Conversely, $(B,A)$ is not locally identified in $\Omega^R\cap\Omega_{n,m,\kappa,\lambda}^{EUI0}$ if the dimension of solutions to \eqref{eq:nonlinearcriterionEUI0} is non-zero and constant in a neighborhood of $(B,A)$ in $\Omega^R\cap\Omega_{n,m,\kappa,\lambda}^{EUI0}$.
\end{thm}

For completeness, we provide results for the identifiability of the $i$-th equation of \eqref{eq:lremformal}.

\begin{defn}\label{defn:identi}
Let $(B,A)\in\Lambda\subset\Omega_{n,m,\kappa,\lambda}^{EU}$, let $i\in\{1,\ldots,n\}$, and let $e_i$ be the $i$-th standard basis vector of $\mathbb{R}^n$. We say that the $i$-th equation of $(B,A)$ is identified in $\Lambda$ if 
\begin{align*}
e_i'(B,A)=e_i'\left(\big((B,A)/\sim\big)\cap\;\Lambda\right).
\end{align*}
We say that the $i$-th equation of $(B,A)$ is locally identified in $\Lambda$ if there is an open set $N\subset\Omega_{n,m,\kappa,\lambda}$ containing $(B,A)$ such that
\begin{align*}
e_i'(B,A)=e_i'\left(\big((B,A)/\sim\big)\cap\;\Lambda\cap N\right).
\end{align*}
\end{defn}

Before we state the next result, we will need the following notation.

\begin{defn}\label{defn:evecemat}
The equation vectorization operator on $\Omega_{n,m,\kappa,\lambda}$ is given by
\begin{align*}
\mathrm{evec}_{\Omega_{n,m,\kappa,\lambda}}\left(\sum_{i=-\lambda}^\kappa \mathsf{X}_iz^i,\sum_{i=0}^\kappa \mathsf{Y}_iz^i\right)=(\mathsf{X}_{-\lambda},\ldots,\mathsf{X}_{\kappa},\mathsf{Y}_{0},\ldots,\mathsf{Y}_{\kappa})',
\end{align*}
defined for all $\mathsf{X}_{-\lambda},\ldots, \mathsf{X}_\kappa\in\mathbb{R}^{1\times n}$ and all $\mathsf{Y}_0,\ldots, \mathsf{Y}_\kappa\in\mathbb{R}^{1\times m}$. The equation reconstruction operator on $\Omega_{n,m,\kappa,\lambda}$ is given by
\begin{align*}
\mathrm{emat}_{\Omega_{n,m,\kappa,\lambda}}\left(\left(\mathsf{x}_{-\lambda}',\ldots,\mathsf{x}_{\kappa}',\mathsf{y}_{0}',\ldots,\mathsf{y}_{\kappa}\right)'\right)=\left(\sum_{i=-\lambda}^\kappa \mathsf{x}_i'z^i,\sum_{i=0}^\kappa \mathsf{y}_i'z^i\right),
\end{align*}
defined for all $\mathsf{x}_{-\lambda},\ldots, \mathsf{x}_\kappa\in\mathbb{R}^n$ and all $\mathsf{y}_0,\ldots, \mathsf{y}_\kappa\in\mathbb{R}^m$.
\end{defn}

\begin{thm}\label{thm:identEUIi}
Let $R_i:\mathbb{R}^{n(\kappa+\lambda+1)+m(\kappa+1)}\rightarrow \mathbb{R}^r$ and let
\begin{align}
\Omega^{R_i}=\left\{(\mathsf{X},\mathsf{Y})\in\Omega_{n,m,\kappa,\lambda}: R_i\left(\mathrm{evec}_{\Omega_{n,m,\kappa,\lambda}}\left(e_i'(\mathsf{X},\mathsf{Y})\right)\right)=0_{r\times 1}\right\}.\label{eq:nonlineari}
\end{align}
If $(B,A)\in\Omega^{R_i}\cap\Omega_{n,m,\kappa,\lambda}^{EUI}$ and $P$ is defined as in Theorem \ref{thm:oeEUI}, then the $i$-th equation of $(B,A)$ is identified in $\Omega^{R_i}\cap\Omega_{n,m,\kappa,\lambda}^{EUI}$ if and only if
\begin{align}
\begin{aligned}
P(V)'\mathrm{evec}_{\Omega_{n,m,\kappa,\lambda}}(\mathsf{x},\mathsf{y})&=0_{m((n+1)\kappa+1)\times 1},\\
R_i\left(\mathrm{evec}_{\Omega_{n,m,\kappa,\lambda}}(\mathsf{x},\mathsf{y})\right)&=0_{r\times1},\\
\mathsf{V}'\mathsf{V}&=I_m,
\end{aligned}
\label{eq:identcriterionEUIi}
\end{align}
has a unique solution in $e_i'\Omega_{n,m,\kappa,\lambda}^{EU}\times\mathbb{R}^{m\times m}$, $\left((\mathsf{x},\mathsf{y}),\mathsf{V}\right)=(e_i'(B,A),I_m)$.
\end{thm}

\begin{thm}\label{thm:localEUIi}
Let $R_i$ and $\Omega^{R_i}$ be as in Theorem \ref{thm:identEUIi} and suppose $R_i$ is continuously differentiable. If $(B,A)\in\Omega^{R_i}\cap\Omega_{n,m,\kappa,\lambda}^{EUI}$ and $P$ is defined as in Theorem \ref{thm:oeEUI}, then the $i$-th equation of $(B,A)$ is locally identified in $\Omega^{R_i}\cap\Omega_{n,m,\kappa,\lambda}^{EUI}$ if
\begin{align}
\begin{aligned}
P(I_m)'\mathrm{evec}_{\Omega_{n,m,\kappa,\lambda}}(\mathsf{dx},\mathsf{dy})+K_i\mathsf{dV}&=0_{m((n+1)\kappa+1)\times1},\\
\nabla R_i(B,A)\mathrm{evec}_{\Omega_{n,m,\kappa,\lambda}}(\mathsf{dx},\mathsf{dy})&=0_{r\times 1},
\\ 
\mathsf{dV}'+\mathsf{dV}&=0_{m\times m},
\end{aligned}
\label{eq:nonlinearcriterionEUIi}
\end{align}
has a unique solution in $e_i'\Omega_{n,m,\kappa,\lambda}\times\mathbb{R}^{m\times m}$, $\left((\mathsf{dx},\mathsf{dy}),\mathsf{dV}\right)=((0_{1\times n},0_{1\times m}),0_{m\times m})$, where
\begin{align*}
K_i=\left[\begin{array}{c}
I_m\otimes e_i'A_0\\
\vdots\\
I_m\otimes e_i'A_\kappa\\
0_{nm\kappa\times m^2}\\
\end{array}
\right].
\end{align*}
Conversely, the $i$-th equation of $(B,A)$ is not locally identified in $\Omega^R\cap\Omega_{n,m,\kappa,\lambda}^{EUI}$ if the dimension of solutions to \eqref{eq:nonlinearcriterionEUIi} is non-zero and constant in a neighborhood of $(B,A)$ in  $\Omega^R\cap\Omega_{n,m,\kappa,\lambda}^{EUI}$.
\end{thm}

\begin{thm}\label{thm:identEUI0i}
Let $R_i$ and $\Omega^{R_i}$ be as in Theorem \ref{thm:identEUIi}. If $(B,A)\in\Omega^{R_i}\cap\Omega_{n,m,\kappa,\lambda}^{EUI0}$ and $P$ is defined as in Theorem \ref{thm:oeEUI}, then the $i$-th equation of $(B,A)$ is identified in $\Omega^{R_i}\cap\Omega_{n,m,\kappa,\lambda}^{EUI0}$ if and only if
\begin{align}
\begin{aligned}
P(I_m)'\mathrm{evec}_{\Omega_{n,m,\kappa,\lambda}}(\mathsf{x},\mathsf{y})&=0_{m((n+1)\kappa+1)\times 1},\\
R_i\left(\mathrm{evec}_{\Omega_{n,m,\kappa,\lambda}}(\mathsf{x},\mathsf{y})\right)&=0_{r\times1},
\end{aligned}
\label{eq:identcriterionEUI0i}
\end{align}
has a unique solution in $e_i'\Omega_{n,m,\kappa,\lambda}^{EU}$, $(\mathsf{x},\mathsf{y})=e_i'(B,A)$.
\end{thm}

\begin{thm}\label{thm:localEUI0i}
Let $R_i$ and $\Omega^{R_i}$ be as in Theorem \ref{thm:localEUIi}. If $(B,A)\in\Omega^{R_i}\cap\Omega_{n,m,\kappa,\lambda}^{EUI0}$ and $P$ is defined as in Theorem \ref{thm:oeEUI}, then the $i$-th equation of $(B,A)$ is locally identified in $\Omega^{R_i}\cap\Omega_{n,m,\kappa,\lambda}^{EUI0}$ if
\begin{align}
\begin{aligned}
P(I_m)'\mathrm{evec}_{\Omega_{n,m,\kappa,\lambda}}(\mathsf{dx},\mathsf{dy})&=0_{m((n+1)\kappa+1)\times 1},\\
\nabla R_i(B,A)\mathrm{evec}_{\Omega_{n,m,\kappa,\lambda}}(\mathsf{dx},\mathsf{dy})&=0_{r\times 1},
\end{aligned}
\label{eq:nonlinearcriterionEUI0i}
\end{align}
has a unique solution in $e_i'\Omega_{n,m,\kappa,\lambda}^{EU}$, $(\mathsf{dx},\mathsf{dy})=(0_{n\times n},0_{n\times m})$.

Conversely, the $i$-th equation of $(B,A)$ is not locally identified in $\Omega^R\cap\Omega_{n,m,\kappa,\lambda}^{EUI0}$ if the dimension of solutions to \eqref{eq:nonlinearcriterionEUI0i} is non-zero and constant in a neighborhood of $(B,A)$ in $\Omega^R\cap\Omega_{n,m,\kappa,\lambda}^{EUI0}$.
\end{thm}

Theorems \ref{thm:identEUIi}-\ref{thm:localEUI0i} have analogous interpretations and similar geometry to Theorems \ref{thm:identEUI}-\ref{thm:localEUI0}.

\section{Conclusion}\label{sec:conclusion}
This paper's title is an hommage to the seminal paper of the VARMA identification literature \citep{hannan71}. Like Hannan's paper, the present work characterizes observational equivalence and provides conditions for identification in a variety of settings. More importantly, and again much like Hannan's paper, the present work has not succeeded in answering all of the questions surrounding the identification of the model under study.  We now turn to some of the pending issues.

Recognizing the difficulty of the identification problem for VARMA models, the literature proposed a variety of canonical parametrizations \cite[p.\ 67]{hd}. These parametrizations allow the researcher to specify a model without having to worry about identification. It would be quite useful for empirical work to find similar parametrizations for LREMs.

Our framework has excluded measurement errors, which are commonly used in the LREM literature. This is not an insurmountable challenge as the literature on latent variables and measurement error is very well developed \citep{bollen,fuller}. Treating this material in the present work would have made it prohibitively long and complicated, not to mention distracting from the primary theoretical difficulties of the identification problem for LREMs. Thus, this is also left for further work.

Having derived conditions for identifiability, the natural next question is to consider near non-identifiability. It is an empirical fact that the parameters of LREMs are frequently weakly identified \citep{mavroeidis04,mavroeidis05}. The results of this paper should be able to shed light on the structural reasons for weak identifiability.

\cite{cochrane} has pointed out that certain parameters of an LREM are important not so much because they directly affect the solution but because they rule out certain classes of solutions. Thus, a given parameter can be important economically at the same time as it is statistically not identified. In this case, it would be useful to identify these parameters analytically rather than following the common practice of assessing identifiability computationally.

Finally, this paper makes clear, as pointed out to us by Peter C.\ B.\ Phillips, that what is needed above and beyond an understanding of the identification problem for LREMs is an encompassing mathematical approach to modeling in which one can address all, point-, set-, weak-, and near-identification issues.

\newpage

\appendix

\section*{Appendix}

\section{Topological Aspects of the Parameter Spaces}\label{sec:topology}
This section develops the topological attributes of our parameter spaces. In particular, it derives their boundaries, interiors, and the number of connected components of the interiors. Such results are important because they can have significant implications for optimization and sampling algorithms. For example, it will always be computationally more convenient to compute likelihood functions and objective functions at interior points, owing to the imprecision of floating point arithmetic. For another example, it is important to be mindful of the fact that over some parameter spaces, algorithms can fail if initialized in the wrong connected component, and caution would then suggest to initialize these algorithms in each connected component of the parameter space.

To get some intuition about the topology of the parameter spaces in this paper, it helps to consider the simplest possible examples.

\begin{exmp}
\begin{enumerate}
\item $\Omega_{1,1,0,0}^{EU}$ is the set of all pairs $(B_0,A_0)\in\mathbb{R}^2$ such that $B_0\neq0$ and consists of two connected components according to whether $B_0\gtrless0$. 

$\Omega_{1,1,0,0}^{EUI}$ is the set of pairs $(B_0,A_0)\in\Omega_{1,1,0,0}^{EU}$ such that, additionally, $A_0\neq0$; this set is open and consists of four connected components according to whether $B_0\gtrless0$ and $A_0\gtrless0$ (the interiors of the standard quadrants of $\mathbb{R}^2$).

$\Omega_{1,1,0,0}^{EUI0}$ is the set of pairs $(B_0,A_0)\in\Omega_{1,1,0,0}^{EUI}$ such that $A_0/B_0>0$; this set is open and consists of two connected components according to whether $B_0\gtrless0$ with $A_0$ of the same sign as $B_0$ (the interiors of the north east and south west standard quadrants of $\mathbb{R}^2$). Notice that the two components of $\Omega_{1,1,0,0}^{EU}$ and the four components of $\Omega_{1,1,0,0}^{EUI}$ are separated by nowhere dense sets but the components of $\Omega_{1,1,0,0}^{EUI0}$ are not.

\item $\Omega_{2,1,0,0}^{EU}$ is the set of pairs $(B_0,A_0)\in\mathbb{R}^{2\times 2}\times\mathbb{R}^{2\times 1}$ such that $\det(B_0)\neq0$; this set is open and consists of two connected components according to whether $\det(B_0)\gtrless0$, separated by a nowhere dense set.

$\Omega_{2,1,0,0}^{EUI}$ is the set of pairs $(B_0,A_0)\in\Omega_{2,1,0,0}^{EUI}$ such that $A_0\neq0_{2\times 1}$; this set is open and consists of two connected components according to whether $\det(B_0)\gtrless0$, separated by a nowhere dense set. 

$\Omega_{2,1,0,0}^{EUI0}$ is the set of pairs $(B_0,A_0)\in\Omega_{2,1,0,0}^{EUI}$ such that $B_0^{-1}A_0$ is CQLT; this set is open and also consists of two connected components separated by a set that is not nowhere dense.

\item $\Omega_{1,1,1,0}^{EU}$ is the set of pairs $(B_0+B_1z,A_0+A_1z)\in\mathbb{R}[z]^2$ such that $|B_1/B_0|<1$. This set is open and its boundary is the set of parameters satisfying $|B_1/B_0|=1$. Every $(B_0+B_1z,A_0+A_1z)\in\Omega_{1,1,1,0}^{EU}$ is path-connected to $(B_0,A_0)$ in $\Omega_{1,1,1,0}^{EU}$ via the path $t\mapsto (B_0+B_1(1-t)z,A_0+A_1(1-t)z)$ for $t\in[0,1]$. It follows that $\Omega_{1,1,1,0}^{EU}$ consists of two connected components like $\Omega_{1,1,0,0}^{EU}$.

$\Omega_{1,1,1,0}^{EUI}$ is the set of pairs $(B_0+B_1z,A_0+A_1z)\in\Omega_{1,1,1,0}^{EU}$ such that $|A_1/A_0|\leq1$. This set is not open (e.g.\ $(1,1-z)\in\Omega_{1,1,1,0}^{EUI}$ is a boundary point) and it is not closed (e.g.\ $(1-z,1)\not\in\Omega_{1,1,1,0}^{EU}$ is a boundary point). Every $(B_0+B_1z,A_0+A_1z)\in\Omega_{1,1,1,0}^{EUI}$ is path-connected to $(B_0,A_0)$ in $\Omega_{1,1,1,0}^{EUI}$ via the same path considered for $\Omega_{1,1,1,0}^{EU}$. It follows that $\Omega_{1,1,1,0}^{EUI}$ consists of four connected components like $\Omega_{1,1,0,0}^{EUI}$.

$\Omega_{1,1,1,0}^{EUI0}$ is the set of pairs $(B_0+B_1z,A_0+A_1z)\in\Omega_{1,1,1,0}^{EUI}$ such that $A_0/B_0>0$. By similar arguments to the above, it is not open or closed and it consists of two connected components. Note that removing the points where $|A_1/A_0|=1$, we obtain the interiors of $\Omega_{1,1,1,0}^{EUI}$ and $\Omega_{1,1,1,0}^{EUI0}$.
\end{enumerate}
\end{exmp}

\begin{thm}\label{thm:topology}
\begin{enumerate}
\item $\Omega_{n,m,\kappa,\lambda}^{EU}$ is an open subset of $\Omega_{n,m,\kappa,\lambda}$ and consists of two connected components.
\item $\Omega_{n,m,\kappa,\lambda}^{EUI}$ is homeomorphic to a subset of $\mathbb{R}^{n^2(\kappa+\lambda+1)+nm(\kappa+1)}$. The interior of $\Omega_{n,m,\kappa,\lambda}^{EUI}$ consists of four connected components if $n=m$ and two connected components if $n>m$.
\item $\Omega_{n,m,\kappa,\lambda}^{EUI0}$ is homeomorphic to a subset of $\mathbb{R}^{n^2(\kappa+\lambda+1)+nm(\kappa+1)-\frac{1}{2}m(m-1)}$. The interior of $\Omega_{n,m,\kappa,\lambda}^{EUI0}$ consists of two connected components.
\end{enumerate}
\end{thm}
\begin{proof}
(i) That $\Omega_{n,m,\kappa,\lambda}^{EU}$ is an open subset of $\Omega_{n,m,\kappa,\lambda}$ follows from Proposition X.1.1 of \cite{cg}. Given $(B,A)\in\Omega_{n,m,\kappa,\lambda}^{EU}$, let $B=B_-B_+$ be the canonical Wiener-Hopf factorization of $B$. Then the path 
\begin{align*}
t\mapsto (B_-(z/(1-t))B_+((1-t)z),(1-t)A(z))
\end{align*}
maps $[0,1]$ to $\Omega_{n,m,\kappa,\lambda}^{EU}$. To see this, note that for every $t\in[0,1]$, $B_-(z/(1-t))B_+((1-t)z)$ is a canonical Wiener-Hopf factorization. The advantage of this path is that it connects $(B,A)$ to $(B_+(0),0_{n\times m})$. Since $B_+(0)$ is invertible, it falls into exactly one of two connected components \citep[Theorem 4]{components}.

(ii) The key to parametrizing $\Omega_{n,m,\kappa,\lambda}^{EUI}$ is to make use of the canonical Wiener-Hopf factorization again and take advantage of the straightforward parametrization of VARMA models $\Omega_{n,m,\kappa,0}^{EUI}$. Let
\begin{align*}
\Theta_{n,m,\kappa,\lambda}^{EUI}=\Bigg\{(\mathsf{F}_\lambda,\ldots,\mathsf{F}_1,\;&\mathsf{B}_0,\ldots,\mathsf{B}_\kappa,\mathsf{A}_0,\mathsf{A}_1,\ldots,\mathsf{A}_\kappa):\\
&\mathsf{F}_\lambda,\ldots,\mathsf{F}_1,\mathsf{B}_0,\ldots,\mathsf{B}_\kappa\in\mathbb{R}^{n\times n},\mathsf{A}_0,\mathsf{A}_1,\ldots,\mathsf{A}_\kappa\in\mathbb{R}^{n\times m},\\
&\mathrm{rank}(\mathsf{F}(z))=n\text{ for all }z\in \mathbb{D}^c,\text{ where }\mathsf{F}(z)=I_n+\sum_{i=1}^\lambda \mathsf{F}_iz^{-i},\\
&\mathrm{rank}(\mathsf{B}(z))=n\text{ for all }z\in\overline{\mathbb{D}},\text{ where }\mathsf{B}(z)=\sum_{i=0}^\kappa \mathsf{B}_iz^i,\\
&\mathrm{rank}(\mathsf{A}(z))=n\text{ for all }z\in\mathbb{D},\text{ where }\mathsf{A}(z)=\sum_{i=0}^\kappa \mathsf{A}_iz^i\Bigg\}.
\end{align*}
Then $\Theta_{n,m,\kappa,\lambda}^{EUI}$ can be viewed as a subset of $\mathbb{R}^{n^2(\kappa+\lambda+1)+nm(\kappa+1)}$. We claim that the mapping $\phi:\Theta_{n,m,\kappa,\lambda}^{EUI}\rightarrow\Omega_{n,m,\kappa,\lambda}^{EUI}$, defined by
\begin{multline*}
(\mathsf{F}_\lambda,\ldots,\mathsf{F}_1,\mathsf{B}_0,\ldots,\mathsf{B}_\kappa,\mathsf{A}_0,\mathsf{A}_1,\ldots,\mathsf{A}_\kappa)\mapsto (\mathsf{F}\mathsf{B},[\mathsf{F}\mathsf{A}]_+)\\
=\left(\left(I_n+\sum_{i=1}^\lambda \mathsf{F}_iz^{-i}\right)\left(\sum_{i=0}^\kappa \mathsf{B}_iz^i\right),\left[\left(I_n+\sum_{i=1}^\lambda \mathsf{F}_iz^{-i}\right)\left(\sum_{i=0}^\kappa \mathsf{A}_iz^i\right)\right]_+\right),
\end{multline*}
is a homeomorphism. Note that $[\mathsf{F}\mathsf{A}]_+\in\mathbb{R}[z]^{n\times m}$ and $\max\deg\left([\mathsf{F}\mathsf{A}]_+\right)\leq\kappa$ and the solution to $(\mathsf{F}\mathsf{B},[\mathsf{F}\mathsf{A}]_+)$ is unique because the factorization $\mathsf{F}\mathsf{B}$ is a canonical Wiener-Hopf factorization. Thus, the transfer function of the solution is
\begin{align*}
\mathsf{B}^{-1}[\mathsf{F}^{-1}[\mathsf{F}\mathsf{A}]_+]_+&=\mathsf{B}^{-1}[\mathsf{F}^{-1}\mathsf{F}\mathsf{A}]_+-\mathsf{B}^{-1}[\mathsf{F}^{-1}[\mathsf{F}\mathsf{A}]_-]_+\\
&=\mathsf{B}^{-1}[\mathsf{A}]_+-0_{n\times m}\\
&=\mathsf{B}^{-1}\mathsf{A},
\end{align*}
where we have used the fact that $\mathsf{F}^{-1}[\mathsf{F}\mathsf{A}]_-$ is strictly proper and has no poles outside $\mathbb{D}$. Clearly, $\phi$ is continuous. To check injectivity, let $(\mathsf{F}\mathsf{B},[\mathsf{F}\mathsf{A}]_+)=(\tilde{\mathsf{F}}\tilde{\mathsf{B}},[\tilde{\mathsf{F}}\tilde{\mathsf{A}}]_+)$. Then $\mathsf{F}=\tilde{\mathsf{F}}$ and $\mathsf{B}=\tilde{\mathsf{B}}$ by the uniqueness of Wiener-Hopf factorization \cite[Theorems I.1.1 and I.1.2]{cg}. On the other hand, $[\mathsf{F}\mathsf{A}]_+=[\tilde{\mathsf{F}}\tilde{\mathsf{A}}]_+$, implies that $[\mathsf{F}(\mathsf{A}-\tilde{\mathsf{A}})]_+=0_{n\times m}$. Multiplying by $\mathsf{F}^{-1}$ and applying the $[\;\cdot\;]_+$ operator, we obtain 
\begin{align*}
0_{n\times m}&=[\mathsf{F}^{-1}[\mathsf{F}(\mathsf{A}-\tilde{\mathsf{A}})]_+]_+\\
&=[\mathsf{F}^{-1}\mathsf{F}(\mathsf{A}-\tilde{\mathsf{A}})]_+-[\mathsf{F}^{-1}[\mathsf{F}(\mathsf{A}-\tilde{\mathsf{A}})]_-]_+\\
&=[(\mathsf{A}-\tilde{\mathsf{A}})]_+-0_{n\times m}\\
&=(\mathsf{A}-\tilde{\mathsf{A}}).
\end{align*}
To check surjectivity, let $(B,A)\in\Omega_{n,m,\kappa,\lambda}^{EUI}$, then we can simply set $\mathsf{F}=B_-$, $\mathsf{B}=B_+$, and $\mathsf{A}=[\mathsf{F}^{-1}A]_+$. Thus, $\phi^{-1}$ is a well-defined function. 
Proposition X.1.1 of \cite{cg} shows that the mapping $B\mapsto(\mathsf{F},\mathsf{B})$ is continuous in the sense that small changes in $B$ in the $L_\infty$ norm on $\mathbb{T}$ lead to small changes in $\mathsf{F}$ and $\mathsf{B}$ in the $L_2$ norm on $\mathbb{T}$. Since $B$, $\mathsf{F}$, and $\mathsf{B}$ are Laurent polynomials of bounded degrees, this implies that the mapping $(B_{-\lambda},\ldots,B_\kappa)\mapsto(\mathsf{F}_\lambda,\ldots,\mathsf{F}_1,\mathsf{B}_0,\ldots,\mathsf{B}_\kappa)$ is continuous. On the other hand, if we expand $\mathsf{F}(z)^{-1}=\sum_{i=0}^\infty \mathsf{F}^iz^{-i}$ for $z\in\mathbb{D}^c$ and consider the coefficients of $I_n=\mathsf{F}(z)\mathsf{F}(z)^{-1}$, we have that $\mathsf{F}^0=I_n$ and $\mathsf{F}^i=-\sum_{j=1}^{\min\{\lambda,i\}}\mathsf{F}_j\mathsf{F}^{i-j}$ for $i\geq1$ and so the elements of $\mathsf{F}^1,\ldots,\mathsf{F}^\kappa$ are polynomials in the elements of $\mathsf{F}_1,\ldots,\mathsf{F}_\lambda$ and are therefore continuous. Since $[\mathsf{F}^{-1}A]_+=\sum_{i=0}^\kappa \left(\sum_{j=i}^\kappa \mathsf{F}^{j-i}A_j\right)z^i$, it follows that $\mathsf{A}$ is also continuous in $(B,A)$. Thus, $\phi^{-1}$ is continuous.

Next, we claim that
\begin{multline*}
\Theta_{n,m,\kappa,\lambda}^{EUI\circ}=\Bigg\{(\mathsf{F}_\lambda,\ldots,\mathsf{F}_1,\mathsf{B}_0,\ldots,\mathsf{B}_\kappa,\mathsf{A}_0,\mathsf{A}_1,\ldots,\mathsf{A}_\kappa)\in\Theta_{n,m,\kappa,\lambda}^{EUI}:\\
\mathrm{rank}(\mathsf{A}(z))=m\text{ for all }z\in\mathbb{T}\Bigg\}
\end{multline*}
is the interior of $\Theta_{n,m,\kappa,\lambda}^{EUI}$. By the continuity of 
zeros of a polynomial with respect to its coefficients \cite[Appendix D]{hj1}, $\Theta_{n,m,\kappa,\lambda}^{EUI\circ}$ is open. Now pick any point $(\mathsf{F}_\lambda,\ldots, \mathsf{F}_1,\mathsf{B}_0,\ldots,\mathsf{B}_\kappa,\mathsf{A}_0,\mathsf{A}_1,\ldots,\mathsf{A}_\kappa)\in\Theta_{n,m,\kappa,\lambda}^{EUI}\backslash\Theta_{n,m,\kappa,\lambda}^{EUI\circ}$, then $\mathrm{rank}(\mathsf{A}(z_0))<m$ for some $z_0\in\mathbb{T}$. Now consider the point $(\mathsf{F}_\lambda,\ldots, \mathsf{F}_1,\mathsf{B}_0,\ldots,\mathsf{B}_\kappa,\mathsf{A}_0,\rho\mathsf{A}_1,\ldots,\rho^\kappa\mathsf{A}_\kappa)$
for some $\rho>1$. This point falls outside $\Theta_{n,m,\kappa,\lambda}^{EUI}$ because $\sum_{i=0}^\kappa \rho^i\mathsf{A}_i z^i=\mathsf{A}(\rho z)$ has a zero in $\mathbb{D}$. It follows that $\Theta_{n,m,\kappa,\lambda}^{EUI}\backslash\Theta_{n,m,\kappa,\lambda}^{EUI\circ}$ are boundary points of $\Theta_{n,m,\kappa,\lambda}^{EUI}$.

Finally, let $(\mathsf{F}_\lambda,\ldots,\mathsf{F}_1,\mathsf{B}_0,\ldots,\mathsf{B}_\kappa,\mathsf{A}_0,\mathsf{A}_1,\ldots,\mathsf{A}_\kappa)\in\Theta_{n,m,\kappa,\lambda}^{EUI\circ}$. Then, we may follow a similar reasoning to the above to show that for $t\in[0,1]$,
\begin{align*}
t\mapsto((1-t)^\kappa\mathsf{F}_\lambda,\ldots, (1-t)\mathsf{F}_1,\mathsf{B}_0,(1-t)\mathsf{B}_1,\ldots,(1-t)^\kappa\mathsf{B}_\kappa,\mathsf{A}_0,(1-t)\mathsf{A}_1,\ldots,(1-t)^\kappa \mathsf{A}_\kappa)
\end{align*}
defines a path in $\Theta_{n,m,\kappa,\lambda}^{EUI\circ}$. Thus, $(\mathsf{F}_\lambda,\ldots,\mathsf{F}_1,\mathsf{B}_0,\ldots,\mathsf{B}_\kappa,\mathsf{A}_0,\mathsf{A}_1,\ldots,\mathsf{A}_\kappa)\in\Theta_{n,m,\kappa,\lambda}^{EUI\circ}$ is in the same connected component as $(0,\ldots,0,\mathsf{B}_0,0,\ldots,0,\mathsf{A}_0,0,\ldots,0)$. It remains to consider the connectedness of the space 
\begin{align*}
\{(\mathsf{B}_0,\mathsf{A}_0)\in\mathbb{R}^{n\times n}\times \mathbb{R}^{n\times m}: \mathrm{rank}(\mathsf{B}_0)=n\text{ and }\mathrm{rank}(\mathsf{A}_0)=m\}.
\end{align*}
The set of real $n\times m$ matrices of rank $m$ consists of two components if $n=m$ and one component if $n>m$ \cite[Theorem 4]{components}. Thus, the set above consists of four components if $n=m$ and two components if $n>m$. The same is true of the interior of $\Omega_{n,m,\kappa,\lambda}^{EUI}$.

(iii) is proven by the same arguments with a few minor changes. We parametrize $\Omega_{n,m,\kappa,\lambda}^{EUI0}$ and its interior by
\begin{align*}
\Theta_{n,m,\kappa,\lambda}^{EUI0}=\Bigg\{(\mathsf{F}_\lambda,\ldots,\mathsf{F}_1,\;&\mathsf{B}_0,\ldots,\mathsf{B}_\kappa,\mathsf{C}_0,\mathsf{A}_1,\ldots,\mathsf{A}_\kappa):\mathsf{C}_0\text{ is CQLT and}\\
&(\mathsf{F}_\lambda,\ldots,\mathsf{F}_1,\mathsf{B}_0,\ldots,\mathsf{B}_\kappa,\mathsf{B}_0\mathsf{C}_0,\mathsf{A}_1,\ldots,\mathsf{A}_\kappa)\in \Theta_{n,m,\kappa,\lambda}^{EUI}\Bigg\},\\
\Theta_{n,m,\kappa,\lambda}^{EUI0\circ}=\Bigg\{(\mathsf{F}_\lambda,\ldots,\mathsf{F}_1,\;&\mathsf{B}_0,\ldots,\mathsf{B}_\kappa,\mathsf{C}_0,\mathsf{A}_1,\ldots,\mathsf{A}_\kappa):\mathsf{C}_0\text{ is CQLT and}\\
&(\mathsf{F}_\lambda,\ldots,\mathsf{F}_1,\mathsf{B}_0,\ldots,\mathsf{B}_\kappa,\mathsf{B}_0\mathsf{C}_0,\mathsf{A}_1,\ldots,\mathsf{A}_\kappa)\in \Theta_{n,m,\kappa,\lambda}^{EUI\circ}\Bigg\}.
\end{align*}
There are now $\frac{1}{2}m(m-1)$ fewer free parameter than before because $\mathsf{C}_0$ is restricted to be CQLT. Then tracing the same steps as in (ii), we find that the mapping
\begin{multline*}
(\mathsf{F}_\lambda,\ldots,\mathsf{F}_1,\mathsf{B}_0,\ldots,\mathsf{B}_\kappa,\mathsf{C}_0,\mathsf{A}_1,\ldots,\mathsf{A}_\kappa)\\
\mapsto\left(\left(I_n+\sum_{i=1}^\lambda \mathsf{F}_iz^{-i}\right)\left(\sum_{i=0}^\kappa \mathsf{B}_iz^i\right),\left[\left(I_n+\sum_{i=1}^\lambda \mathsf{F}_iz^{-i}\right)\left(\mathsf{B}_0\mathsf{C}_0+\sum_{i=1}^\kappa \mathsf{A}_iz^i\right)\right]_+\right)
\end{multline*}
is a homeomorphism from $\Theta_{n,m,\kappa,\lambda}^{EUI0}$ to $\Omega_{n,m,\kappa,\lambda}^{EUI0}$ (resp.\ from $\Theta_{n,m,\kappa,\lambda}^{EUI0\circ}$ to $\Omega_{n,m,\kappa,\lambda}^{EUI0\circ}$). In the final step, we have the problem of finding the number of connected components of the set
\begin{align*}
\{(\mathsf{B}_0,\mathsf{C}_0)\in\mathbb{R}^{n\times n}\times \mathbb{R}^{n\times m}: \mathrm{rank}(\mathsf{B}_0)=n\text{ and }\mathsf{C}_0\text{ is CQLT}\}.
\end{align*}
The set of $n\times n$ matrices of rank $n$ consists of two connected components \cite[Theorem 4]{components}, while the set of $n\times m$ CQLT matrices consists of one connected component because every such matrix is connected by a straight line to $\left[\begin{smallmatrix} I_n \\ 0_{(n-m)\times m}\end{smallmatrix}\right]$. Thus $\Omega_{n,m,\kappa,\lambda}^{EUI0}$ consists of two connected components.
\end{proof}

\section{Proofs}\label{sec:proofs}

\begin{proof}[Proof of Theorem \ref{thm:onatski}]
In order to prove this result, we will need to review some concepts from Fourier analysis. Let $\mu$ be normalized Lebesgue measure on $\mathbb{T}$,
\begin{align*}
\mu\left(\{e^{\mathrm{i}s}: a\leq s\leq b\}\right)=\frac{1}{2\pi}(b-a),\qquad 0\leq b-a\leq 2\pi.
\end{align*}
For positive integers $n$ and $m$, the frequency domain of $n$-dimensional random vectors generated by a white noise process of mean zero and variance $I_m$ is the space $L_2^{n\times m}$ of equivalence classes of measurable mappings $D:\mathbb{T}\rightarrow\mathbb{C}^{n\times m}$ such that
\begin{align*}
\int \|D\|^2d\mu<\infty,
\end{align*}
where any two such mappings $D$ and $\tilde{D}$ are equivalent if $\int\|D-\tilde{D}\|^2d\mu=0$. With addition and scalar multiplication defined in the usual way and the inner product
\begin{align*}
\langle D,\tilde{D}\rangle=\int \mathrm{tr}(D^\ast\tilde{D}) d\mu,\qquad D,\tilde{D}\in L_2^{n\times m},
\end{align*}
$L_2^{n\times m}$ is a complex Hilbert space.

Define $\chi^k:z\mapsto z^k$ for $z\in\mathbb{T}$ and $k\in\mathbb{Z}$ and $E_{ij}\in\mathbb{C}^{n\times m}$ to be the matrix with all zero entries except for the $(i,j)$-th which is equal to 1. Then
\begin{align*}
\left\{E_{ij}\chi^k: k\in\mathbb{Z}, i=1,\ldots,n, j=1,\ldots, m\right\}
\end{align*}
is an orthonormal basis for $L_2^{n\times m}$. Thus, every $D\in L_2^{n\times m}$ has a well-defined Fourier series expansion $\sum_{i=-\infty}^\infty D_i\chi^i$, which converges in $L_2^{n\times m}$, $\int\|D\|^2d\mu=\sum_{i=-\infty}^\infty \|D_i\|^2<\infty$. We denote by $H_2^{n\times m}$ the closure in $L_2^{n\times m}$ of the linear span of the orthonormal set
\begin{align*}
\left\{E_{ij}\chi^k: k\geq0, i=1,\ldots,n, j=1,\ldots, m\right\}.
\end{align*}
The orthogonal projection of $D\in L_2^{n\times m}$ onto $H_2^{n\times m}$ is denoted by $P(D|H_2^{n\times m})$ and we set $Q(D|H_2^{n\times m})=D-P(D|H_2^{n\times m})$. If $D$ has Fourier series expansion $\sum_{i=-\infty}^\infty D_i\chi^i$, then $P(D|H_2^{n\times m})=\sum_{i=0}^\infty D_i\chi^i$.

Finally, if $D\in\mathbb{R}(z)^{n\times m}$ has no poles on $\mathbb{T}$, then it has a Laurent series expansion in a neighbourhood of $\mathbb{T}$, $\sum_{i=-\infty}^\infty D_iz^i$, and therefore its Fourier series expansion is $\sum_{i=-\infty}^\infty D_i\chi^i$. This then implies that $P(D|H_2^{n\times m})=[D]_+$ on $\mathbb{T}$.

We can now proceed to the proof, which is in three steps.

\bigskip

\noindent\emph{STEP 1: $(B,A)\in\Omega_{n,m,\kappa,\lambda}$ has a unique solution if and only if there exists a unique $C\in H_2^{n\times m}$ satisfying.}
\begin{align}
P(BC|H_2^{n\times m})=A.\label{eq:toeplitz}
\end{align}

\bigskip

This fact was proved by \cite{onatski} in the course of proving his Proposition 1. We present his argument here for the sake of completeness. Recall that the frequency domain representation of $\mathbb{E}(Y_{t-i}|\varepsilon_t,\varepsilon_{t-1},\ldots)$ for $(Y_t)_{t\in\mathbb{Z}}$ defined by \eqref{eq:ma} is $P(\chi^iC|H_2^{n\times m})$ \cite[eq.\ (4.22)]{rozanov}. Therefore, the frequency domain representation of the left hand side of \eqref{eq:lrem} is $\sum_{i=-\lambda}^\kappa B_iP(\chi^iC|H_2^{n\times m})=P(BC|H_2^{n\times m})$, while that of the right hand side is $A$. In other words, \eqref{eq:toeplitz} is the frequency domain representation of equations \eqref{eq:lrem}. The result then follows from the spectral representation property \cite[p.\ 33]{rozanov}.

\bigskip

\noindent\emph{STEP 2: The only if part.}

\bigskip

Suppose $(B,A)\in\Omega_{n,m,\kappa,\lambda}$ has a unique solution with transfer function $C\in H_2^{n\times m}$. From step 1, $C$ is the unique element of $H_2^{n\times m}$ satisfying \eqref{eq:toeplitz}. We may assume that $\det(B)$ is not identically zero, for otherwise $C+U$ would also be the transfer function of a solution to $(B,A)$ for any $0_{n\times m}\neq U\in\mathbb{R}^{n\times m}$ whose columns lay in the right null space of $B$. That is, $(B,A)$ cannot have a unique solution if $\det(B)$ is identically zero. It follows from Theorem 3.2 of \cite{ls} that $B$ has an inner-limit Wiener-Hopf factorization with respect to $\mathbb{T}$,
\begin{align*}
B=B_-B_0B_+,
\end{align*}
where $B_-\in\mathbb{R}[z^{-1}]^{n\times n}$ is of full rank in $\mathbb{D}^c$ and $\lim_{z\rightarrow\infty}B_-(z)$ is non-singular, $B_+\in\mathbb{R}[z]^{n\times n}$ is of full rank in $\mathbb{D}$, and $B_0$ is a unique diagonal matrix with diagonal elements $z^{\tilde{\kappa}_1},\ldots,z^{\tilde{\kappa}_n}$ such that $\tilde{\kappa}_1\geq\cdots\geq\tilde{\kappa}_n$ are integers known as partial indices. Now multiplying by $B_-^{-1}$ and projecting onto $H_2^{n\times m}$ both sides of \eqref{eq:toeplitz}, we obtain
\begin{align*}
P(B_0B_+C|H_2^{n\times m})&=P(B_-^{-1}A|H_2^{n\times m}).
\end{align*}
Indeed, the Fourier expansion of $BC$ has only finitely many terms with negative powers of $\chi$ so that the Fourier series expansion of $Q(BC|H_2^{n\times m})$ consists of finitely many terms with negative powers of $\chi$, while the Fourier series expansion of $B_-^{-1}$ has no terms with positive powers of $\chi$. Thus, the Fourier series expansion of $B_-^{-1}Q(BC|H_2^{n\times m})$ only has terms with negative powers of $\chi$ and it must be orthogonal to $H_2^{n\times m}$. We have
\begin{align*}
P(B_-^{-1}P(BC|H_2^{n\times m})|H_2^{n\times m})&=P(B_-^{-1}BC|H_2^{n\times m})-P(B_-^{-1}Q(BC|H_2^{n\times m})|H_2^{n\times m})\\
&=P(B_-^{-1}B_-B_0B_+C|H_2^{n\times m})-0_{n\times m}\\
&=P(B_0B_+C|H_2^{n\times m}).
\end{align*}
If $\tilde{\kappa}_n<0$, there cannot be a unique solution because then $C+B_+^{-1}U\chi$ would also be the transfer function of a solution to $(B,A)$, where $U\in\mathbb{R}^{n\times m}$ is any matrix with non-zero $n$-th row and zero everywhere else. Therefore, $\tilde{\kappa}_n\geq0$ and so $P(B_0B_+C|H_2^{n\times m})=B_0B_+C$. It follows that
\begin{align*}
C=B_+^{-1}B_0^{-1}P(B_-^{-1}A|H_2^{n\times m}).
\end{align*}
By the same argument used to show that $Q(BC|H_2^{n\times m})\in\mathbb{R}[\chi^{-1}]^{n\times m}$, $P(B_-^{-1}A|H_2^{n\times m})\in\mathbb{R}[\chi]^{n\times m}$ and so $C\in\mathbb{R}(\chi)^{n\times m}$. Since a rational function is uniquely determined by its values on an infinite set, such as $\mathbb{T}$, there exists a unique element of $\mathbb{R}(z)^{n\times m}$ that coincides with $C$ on $\mathbb{T}$. If the reader will forgive the abuse of notation in referring to this extension as $C$ as well, then $C$ must be an element of $\Sigma^{n\times m}$ and
\begin{align*}
[BC]_+=P(BC|H_2^{n\times m})=A
\end{align*}
at every points of $\mathbb{T}$. Since the left-most and right-most expressions above are rational and agree on $\mathbb{T}$, they must agree everywhere (by the aforementioned uniqueness property again) and \eqref{eq:rationaltoeplitz} follows. There cannot be another solution in $\Sigma^{n\times m}$ because $\{D(\chi): D\in\Sigma^{n\times m}\}\subset H_2^{n\times m}$ and there are not other solutions in $H_2^{n\times m}$.

\bigskip

\noindent\emph{STEP 3: The if part.}

\bigskip

Suppose $C$ is the unique element of $\Sigma^{n\times m}$ that solves \eqref{eq:rationaltoeplitz}. Then 
\begin{align*}
P(BC|H_2^{n\times m})=A.
\end{align*}
If $\tilde{C}\in H_2^{n\times m}$ is a solution to \eqref{eq:toeplitz}, then
\begin{align*}
P(B\tilde{C}|H_2^{n\times m})=A.
\end{align*}
Combining the last two equations, $P(B(C-\tilde{C})|H_2^{n\times m})=0_{n\times m}$. The same arguments used in step 2, imply that $\det(B)$ is not identically zero and $B$ has an inner-limit Wiener-Hopf factorization $B_-B_0B_+$ with respect to $\mathbb{T}$ with non-negative partial indices. We therefore arrive at $B_0B_+(C-\tilde{C})=0_{n\times m}$, which implies that $\tilde{C}=C$.
\end{proof}

\begin{proof}[Proof of Theorem \ref{thm:oe}]
(i) Let $(B,A),(\tilde B,\tilde A)\in\Omega_{n,m,\kappa,\lambda}^{EUI}$ have solutions with transfer functions $C$ and $\tilde{C}$ respectively. If $[\tilde{B}CV]_+=\tilde{A}$ for some $V\in O_m$, then $\tilde{C}=CV$ by Theorem \ref{thm:onatski}. It follows that $\tilde{C}\tilde{C}^\ast=CC^\ast$ and so $(B,A)\sim(\tilde B,\tilde A)$.

Conversely, if $(B,A)\sim(\tilde B,\tilde A)$, then $CC^\ast=\tilde{C}\tilde{C}^\ast$ and so
\begin{align*}
CV=\tilde{C}
\end{align*}
for some $V\in O_m$ \citep[Theorem 2 (3)]{baggioferrante}. If we multiply both sides by $\tilde{B}$ and apply the $[\;\cdot\;]_+$ operator, we obtain
\begin{align*}
[\tilde{B}CV]_+=[\tilde{B}\tilde{C}]_+=\tilde{A}.
\end{align*}

(ii) is proven by exactly the same arguments.
\end{proof}

For the next result, we will need the following lemma.

\begin{lem}\label{lem:trunc}
Let $f\in\mathbb{R}[z,z^{-1}]$ and $g\in\mathbb{R}[z]$ with $g^{-1}\in\Sigma$. Then, $[f/g]_+=0$ if and only if the first $\max(\max\deg(f)+1,\max\deg(g))$ coefficients of the Taylor series expansion of $[f/g]_+$ in a neighborhood of $z=0$ are equal to zero.
\end{lem}
\begin{proof}
Since $g$ and $[f/g]_+$ are elements of $\Sigma$, $h=g[f/g]_+\in\Sigma$. If we can show that $[z^{-\delta}h]_+=0$ for some non-negative integer $\delta$, then $h\in\mathbb{R}[z]$ and $\max\deg(h)<\delta$. \begin{align*}
\left[z^{-\delta}h\right]_+&=\left[z^{-\delta}g\left[f/g\right]_+\right]_+\\
&=\left[z^{-\delta}g(f/g)\right]_+-\left[z^{-\delta}g\left[f/g\right]_-\right]_+& &(f/g=[f/g]_-+[f/g]_+)\\
&=\left[z^{-\delta}f\right]_+& &\text{(if $\delta\geq\max\deg(g)$)}\\
&=0& &\text{(if $\delta>\max\deg(f)$).}
\end{align*}
Therefore, $\delta=\max(\max\deg(f)+1,\max\deg(g))$ provides the desired integer. It follows that $[f/g]_+=h/g$, where $\max\deg(h)\leq\max(\max\deg(f),\max\deg(g)-1)$. The claim then follows from the well-known fact that a polynomial ratio $h/g$ with $g(0)\neq0$ is identically zero if the first $\max\deg(h)+1$ Taylor series coefficients of $h/g$ are equal to zero \citep[Lemma A.4]{dr}.
\end{proof}

We will also need to derive some properties of a submatrix of the infinite system of equations identified in Section \ref{sec:oe}.

\begin{lem}\label{lem:hankellrem}
Let $(B,A)\in\Omega_{n,m,\kappa,\lambda}^{EU}$, let $C(z)=\sum_{i=0}^\infty C_iz^i$ be the transfer function of its solution, and let
\begin{align*}
H=\left[\begin{array}{cccc}
	C_{\kappa+\lambda+1} & C_{\kappa+\lambda+2} & \cdots & C_{(n+1)\kappa+\lambda}\\
	\ddots & \ddots & \ddots & \vdots \\
	C_2 & \ddots & \ddots & C_{n\kappa+1}\\
	C_1 & C_2 & \ddots & C_{n\kappa}
\end{array}\right].
\end{align*}
Then:
\begin{enumerate}
\item $\mathrm{rank}(H)=\delta(C(z^{-1})-C(0))\leq n\kappa$.
\item The set of parameters $(B,A)\in\Omega_{n,m,\kappa,\lambda}^{EUI}$ satisfying $\mathrm{rank}(H)=n\kappa$ contains an open and dense subset of $\Omega_{n,m,\kappa,\lambda}^{EUI}$.
\item The set of parameters $(B,A)\in\Omega_{n,m,\kappa,\lambda}^{EUI0}$ satisfying $\mathrm{rank}(H)=n\kappa$ contains an open and dense subset of $\Omega_{n,m,\kappa,\lambda}^{EUI0}$.
\end{enumerate}
\end{lem}

\begin{proof}[Proof of Lemma \ref{lem:hankellrem}]
(i) The rank of the infinite Hankel matrix
\begin{align*}
\left[\begin{array}{cccc}
	C_1 & C_2 & C_3 & \iddots \\
	C_2 & C_3 & C_4 & \iddots \\
	C_3 & C_4 & C_5 & \iddots \\
	\iddots & \iddots & \iddots & \iddots
\end{array}\right]
\end{align*}
is equal to $\delta(C(z^{-1})-C(0))$ \citep[Theorem~2.4.1 (iii)]{hd}. Now write
\begin{align*}
C(z^{-1})-C(0)&=B_+^{-1}(z^{-1})\left([B_-^{-1}A]_+(z^{-1})-B_+(z^{-1})C(0)\right)\\
&=\left(z^\kappa B_+(z^{-1})\right)^{-1}\left(z^\kappa [B_-^{-1}A]_+(z^{-1})-z^\kappa B_+(z^{-1})C(0)\right).
\end{align*}
It is easy to check that $z^\kappa B_+(z^{-1})\in\mathbb{R}[z]^{n\times n}$ and $z^\kappa [B_-^{-1}A]_+(z^{-1})\in\mathbb{R}[z]^{n\times m}$. Since $B_+(0)$ is non-singular, $\max\deg\left(\det\left(z^\kappa B_+(z^{-1})\right)\right)=n\kappa$ \citep[p.\ 42]{hd}. It follows that $\delta(C(z^{-1})-C(0))\leq n\kappa$ \cite[Lemma~2.2.1 (e)]{hd}. By Theorem~2.4.1 (iii) of \cite{hd} again and reordering the blocks, $\delta(C(z^{-1})-C(0))$ is the rank of the matrix
\begin{align*}
Q=\left[\begin{array}{cccc}
	C_{n\kappa} & C_{n\kappa+1} & \cdots & C_{2n\kappa-1}\\
	\ddots & \ddots & \ddots & \vdots \\
	C_2 & \ddots & \ddots & C_{n\kappa+1}\\
	C_1 & C_2 & \ddots & C_{n\kappa}
\end{array}\right].
\end{align*}
Since $[B_-^{-1}A]_+=B_+C$ is a polynomial matrix of degree at most $\kappa$, the $\kappa+1,\kappa+2,\ldots, 2n\kappa-1$ coefficient matrices of $B_+C$ are all zero. That is,
\begin{align*}
\sum_{i=0}^\kappa B_{+,i}C_{j-i}=0,\qquad\kappa+1\leq j\leq 2n\kappa-1.
\end{align*}
Since $B_{+,0}=B_+(0)$ is non-singular,
\begin{align*}
C_j\;=-B_{+,0}^{-1}\sum_{i=1}^\kappa B_{+,i}C_{j-i},\qquad\kappa+1\leq j\leq 2n\kappa-1.
\end{align*}
This implies that all of the top blocks of $Q$ are linear dependent on the bottom $\kappa$ blocks of $Q$. This implies that 
\begin{align*}
\delta(C(z^{-1})-C(0))=\mathrm{rank}(Q)=\mathrm{rank}(H)=\mathrm{rank}(\check H),
\end{align*}
where
\begin{align*}
\check H=\left[\begin{array}{cccc}
	C_{\kappa} & C_{\kappa+1} & \cdots & C_{(n+1)\kappa-1}\\
	\ddots & \ddots & \ddots & \vdots \\
	C_2 & \ddots & \ddots & C_{n\kappa+1}\\
	C_1 & C_2 & \ddots & C_{n\kappa}
\end{array}\right]
\end{align*}
will be needed later.

(ii) The proof is in two steps.

\bigskip

\noindent \emph{STEP 1: The set of $(B,A)\in\Omega_{n,m,\kappa,\lambda}^{EUI}$ such that $\det(B_+)$ has distinct zeros, $\det(B_\kappa)\neq0$, the transfer function of the solution is strictly invertible, and $(B_+,[B_-^{-1}A]_+)$ is left prime is open and dense in $\Omega_{n,m,\kappa,\lambda}^{EUI}$.}

\bigskip

Define the set
\begin{align*}
\Omega_{n,m,\kappa,\lambda}^{EUI\circ}=\Big\{(B,A)&\in\Omega_{n,m,\kappa,\lambda}^{EUI}:\det(B_+)\text{ has distinct zeros}, \det(B_\kappa)\neq0,\\
&\mathrm{rank}(B_+^{-1}[B_-^{-1}A](z))=m\text{ for all }z\in\mathbb{T},\text{ and }(B_+,[B_-^{-1}A]_+)\text{ is left prime}\Big\}.
\end{align*}
To see that $\Omega_{n,m,\kappa,\lambda}^{EUI\circ}$ is open, let $(B,A)\in\Omega_{n,m,\kappa,\lambda}^{EUI\circ}$. We will construct a neighbourhood of $(B,A)$ in $\Omega_{n,m,\kappa,\lambda}^{EUI\circ}$ as $N=N_1\cap N_2\cap N_3\cap N_4$, where $N_i$ is a neighbourhood of $(B,A)$ satisfying the $i$-th additional condition in $\Omega_{n,m,\kappa,\lambda}^{EUI\circ}$.

By the continuity of Wiener-Hopf factorization \cite[Proposition X 1.1]{cg} and the continuity of zeros of polynomials \citep[Appendix D]{hj1}, there is a neighbourhood $N_1$ of $(B,A)$ in $\Omega_{n,m,\kappa,\lambda}^{EUI\circ}$ such that for every $(\check B,\check A)\in N_1$, $\det(\check B_+)$ has distinct zeros.

Since $\det(B_\kappa)\neq0$, the continuity of the determinant implies that there is a neighbourhood $N_2$ of $(B,A)$ in $\Omega_{n,m,\kappa,\lambda}^{EUI\circ}$ such that for every $(\check B,\check A)\in N_2$, $\det(\check B_\kappa)\neq0$.

Since $\mathrm{rank}([B_-^{-1}A]_+)=m$ for all $z\in\overline{\mathbb{D}}$, it has a minor of order $m$ that has no zeros in $\overline{\mathbb{D}}$. Using the continuity of Wiener-Hopf factorization again and the continuity of zeros of a polynomial, there is a neighbourhood $N_3$ of $(B,A)$ in $\Omega_{n,m,\kappa,\lambda}^{EUI\circ}$ such that for every $(\check B,\check A)\in N_3$, $\mathrm{rank}([\check B_-^{-1}\check A]_+)=m$ for all $z\in\overline{\mathbb{D}}$.

Finally, since $\det(B_+)$ has distinct zeros and $B_{+,\kappa}=B_\kappa$ is non-singular, $\det(B_+)$ has $n\kappa$ distinct zeros $z_1,\ldots,z_{n\kappa}\in\mathbb{C}$. The fact that $(B_+,[B_-^{-1}A]_+)$ is left prime then implies that
\begin{align*}
\Pi_i\;[B_-^{-1}A]_+(z_i)\neq0, \qquad i=1,\ldots, n\kappa,
\end{align*}
where $\Pi_i$ is the orthogonal projection matrix onto the left null space of $B_+(z_i)$. Since
\begin{align*}
\mathrm{rank}(B_+(z_i))=n-1,\qquad i=1,\ldots,n\kappa,
\end{align*}
small perturbations to $B_+(z_i)$ that leave the rank fixed at $n-1$ lead to small perturbations to $\Pi_i$ \citep[Theorem~13.5.1]{glr}. Thus, there exists a neighbourhood $N_4$ of $(B,A)$ in $\Omega_{n,m,\kappa,\lambda}^{EUI\circ}$ such that for every $(\check B,\check A)\in N_4$, 
\begin{align*}
\check\Pi_i[\check B_-^{-1}\check A]_+(\check z_i)\neq0, \qquad i=1,\ldots, n\kappa,
\end{align*}
where $\check z_1,\ldots,\check z_{n\kappa}$ are the zeros of $\det(\check B_+)$ and $\check\Pi_i$ is the orthogonal projection matrix onto the left null space of $\check B_+(\check z_i)$. In other words, every $(\check B,\check A)\in N_4$ is left prime.

To see that $\Omega_{n,m,\kappa,\lambda}^{EUI\circ}$ is dense, let $(B,A)\in\Omega_{n,m,\kappa,\lambda}^{EUI}\backslash\Omega_{n,m,\kappa,\lambda}^{EUI\circ}$. We propose an infinitesimal perturbation of $B$, followed by an infinitesimal perturbation $A$ that leads to a point in $\Omega_{n,m,\kappa,0}^{EUI\circ}$. Thus, any neighbourhood of $(B,A)$ in $\Omega_{n,m,\kappa,\lambda}^{EUI}$ will contain an element of $\Omega_{n,m,\kappa,\lambda}^{EUI\circ}$.

If $\det(B_+)$ has a zero of multiplicity greater than 1, \cite{andersonetal} show that there is an infinitesimal perturbation of $B_+$ that maintains maximum degree of at most $\kappa$, and that leads to a matrix polynomial whose determinants has distinct zeros. Such a perturbation can be made small enough that no zero moves in inside $\overline{\mathbb{D}}$.

Next, if $B_\kappa$ is singular, we can infinitesimally perturb its singular values at zero to obtain a non-singular matrix. By the continuity of zeros of polynomials, this infinitesimal perturbation does not interfere with the $\det(B_+)$ having distinct zeros outside $\overline{\mathbb{D}}$.

Next, if $\mathrm{rank}([B_-^{-1}A]_+(z_0))<m$ for some $z_0\in\mathbb{T}$ then for $\rho<1$ arbitrarily near 1, we can perturb $B_-(z)$ to $B_-^\rho(z)=B_-(\rho z)$ and $A$ to $A^\rho(z)=A(\rho z)$ such that $[B_-^{\rho,-1}A^\rho]_+(z)=[B_-^{-1}A]_+(\rho z)$ is of full rank in $\overline{\mathbb{D}}$.

Finally, suppose that after the sequence of infinitesimal perturbations above we arrive at a $(B,A)\in\Omega_{n,m,\kappa,\lambda}^{EUI}$ such that $(B_+,[B_-^{-1}A]_+)$ is not left prime. Enumerate the zeros of $\det(B_+)$ as $z_1,\ldots,z_{n\kappa}$ and choose non-zero vectors $v_1,\ldots,v_{n\kappa}\in\mathbb{C}^n$ spanning the left null spaces of $B_+(z_1),\ldots, B_+(z_{n\kappa})$ respectively. Since $(B_+,[B_-^{-1}A]_+)$ is not left prime, there is an index $i$ such that $v_i'[B_-^{-1}A]_+(z_i)=0$. Choose $U\in\mathbb{R}^{n\times m}$ satisfying $v_i'U\neq0$ for all $i$ such that $v_i'[B_-^{-1}A]_+(z_i)=0$. Then an infinitesimal perturbation of $A_\kappa$ in the direction of $U$ is sufficient to produce a left prime element pair. This infinitesimal perturbation should maintain strict invertibility by the continuity of zeros of polynomials with respect to their coefficients.

\bigskip

\noindent STEP 2: \emph{For all parameters in $\Omega_{n,m,\kappa,\lambda}^{EUI\circ}$, $\mathrm{rank}(H)=n\kappa$.}

\bigskip

Let $(B,A)\in\Omega_{n,m,\kappa,\lambda}^{EUI\circ}$ have a solution with transfer function $C$. We have already established that the matrix $\check H$ encountered in (i) is of rank $\delta(C(z^{-1})-C(0))\leq n\kappa$. If $\mathrm{rank}(\check H)<n\kappa$, then there exist vectors $x_i\in\R^n$, $i=0,\ldots,\kappa-1$, not all zero, such that
\begin{align*}
(x_0',\ldots,x_{\kappa-1}')\check H=0_{1\times nm\kappa}.
\end{align*}
Setting
\begin{align*}
x(z)=\sum_{i=0}^{\kappa-1}x_iz^i,
\end{align*}
this implies that the terms of degree $\kappa,\kappa+1,\ldots,(n+1)\kappa-1$ of $x'C$ vanish. To see that indeed all higher degree terms vanish as well, notice that each element of the numerator of
\begin{align*}
y'=x'C=\frac{x'\mathrm{adj}(B_+)[B_-^{-1}A]_+}{\det(B_+)}
\end{align*}
is expressible as a polynomial of degree bounded above by $\max\deg(x)+\max\deg(\mathrm{adj}(B_+))+\max\deg([B_-^{-1}A]_+)\leq (\kappa-1)+(n-1)\kappa+\kappa=(n+1)\kappa-1$. It follows from Lemma \ref{lem:trunc} that $y\in\mathbb{R}^m[z]$ and $\max\deg(y)\leq\kappa-1$. Now setting
\begin{align*}
U=\left[\begin{smallmatrix} x'B_+^{-1}\\ S\end{smallmatrix}\right]
\end{align*}
with $S\in\R^{(n-1)\times n}$ chosen so that $\det(U)$ is not identically zero (e.g.\ choose $z_0\in\mathbb{D}$ such that $x(z_0)\neq0$ and choose $S$ as an orthogonal complement to $x'(z_0)B_+^{-1}(z_0)$). Then
\begin{gather*}
\dot B=UB_+=\left[\begin{smallmatrix} x'\\ SB_+\end{smallmatrix}\right]\in\R[z]^{n\times n},\qquad 
\dot A=U[B_-^{-1}A]_+=\left[\begin{smallmatrix} y'\\ S[B_-^{-1}A]_+\end{smallmatrix}\right]\in\R[z]^{n\times m},\\
\dot B^{-1}\dot A=B_+^{-1}[B_-^{-1}A]_+=C,
\end{gather*}
and
\begin{align*}
\max\deg(\det(\dot B))\leq(n-1)\kappa<n\kappa.
\end{align*}
On the other hand, $(B,A)\in\Omega_{n,m,\kappa,0}^{EUI\circ}$ implies that $(B_+,[B_-^{-1}A]_+)$ is left prime and also $\max\deg(\det(B_+))=n\kappa$. By Lemma~2.2.1 (e) \cite{hd}, then
\begin{align*}
\max\deg(\det(\dot B))\geq n\kappa.
\end{align*}
We have reached a contradiction. Therefore, $\mathrm{rank}(\check H)=\mathrm{rank}(H)=n\kappa$.

(iii) By the same arguments used in (ii), the set of $(B,A)\in\Omega_{n,m,\kappa,\lambda}^{EUI0}$ such that $\det(B_+)$ has distinct zeros, $\det(B_\kappa)\neq0$, the transfer function of the solution is strictly invertible, and $(B_+,[B_-^{-1}A]_+)$ is left prime is open and dense in $\Omega_{n,m,\kappa,\lambda}^{EUI0}$. The only arguments that require additional care are the ones that concern perturbations. However, none of the perturbations considered affect the CQLT restriction. Step 2 requires no changes at all.
\end{proof}

\begin{proof}[Proof of Theorem \ref{thm:oeEUI}]
(i) By Theorem \ref{thm:oe} (i) and using \eqref{eq:onatski} $(\tilde B,\tilde A)\sim(B,A)$ is equivalent to
\begin{align*}
0&=\left[-\tilde BC+\tilde AV'\right]_+\\
&=\left[-\tilde BB_+^{-1}[B_-^{-1}A]_++\tilde AV'\right]_+\\
&=\left[\frac{-\tilde B \mathrm{adj}(B_+)[B_-^{-1}A]_++\det(B_+)\tilde AV'}{\det(B_+)}\right]_+.
\end{align*}
Thus, each element of $-\tilde BC+\tilde AV'$ can be expressed as a rational function with numerator of maximum degree at most 
\begin{multline*}
\max\Big\{\max\deg(\tilde B)+\max\deg(\mathrm{adj}(B_+))+\max\deg([B_-^{-1}A]_+),\\
\max\deg(\det(B_+))+\max\deg(\tilde A V')\Big\}\leq\max\{\kappa+(n-1)\kappa+\kappa,n\kappa+\kappa\}=(n+1)\kappa
\end{multline*}
and denominator $\det(B_+)$. By Lemma \ref{lem:trunc}, this is equivalent to the first $(n+1)\kappa+1$ Taylor series coefficients of $[-\tilde BC+\tilde AV']_+$ equating to zero. Thus, observational equivalence is equivalent to
\begin{equation*}
-\left[\begin{array}{ccc}
	\tilde B_{-\lambda} & \ldots & \tilde B_\kappa
\end{array}\right] \left[\begin{array}{cc}
 T &  H
\end{array}\right]+\left[\begin{array}{cccccc}
	\tilde A_0V' & \ldots & \tilde A_\kappa V' & 0_{n\times m} & \ldots & 0_{n\times m}
\end{array}\right]=0_{n\times ((n+1)\kappa+1)m}
\end{equation*}
or equivalently 
\begin{align*}
	\left[\begin{array}{cccccc}
	\tilde B_{-\lambda} & \cdots &\tilde  B_{\kappa} &\tilde  A_0 & \cdots & \tilde A_\kappa
\end{array}\right] P(V)= 0_{n\times ((n+1)\kappa+1)m}.
\end{align*}
Vectorizing we obtain
\begin{align*}
\left(P(V)'\otimes I_n\right)
\mathrm{vec}_{\Omega_{n,m,\kappa,\lambda}}(\tilde{B},\tilde{A})
=0_{nm((n+1)\kappa+1)\times 1}.
\end{align*}

(ii) Let $(\tilde B,\tilde A)\in\big((B,A)/\sim\big)\cap\;\Omega_{n,m,\kappa,\lambda}^{EUI}$. Then $(\tilde B,\tilde A)\in\Omega_{n,m,\kappa,\lambda}^{EU}$ by inclusion. In turn, (i) implies that there exists a $V\in O_m$ such that $(\tilde B,\tilde A)\in\mathrm{mat}_{\Omega_{n,m,\kappa,\lambda}}(\ker(P(V)'\otimes I_n))$. Thus,
\begin{align*}
(\tilde B,\tilde A)\in\left(\bigcup_{V\in\mathbb{R}^{m\times m}, V'V=I_m}\mathrm{mat}_{\Omega_{n,m,\kappa,\lambda}}\left(\ker\left(P(V)'\otimes I_n\right)\right)\right)\cap\Omega_{n,m,\kappa,\lambda}^{EU}.
\end{align*}
If, on the other hand, $(\tilde B,\tilde A)\in\mathrm{mat}_{\Omega_{n,m,\kappa,\lambda}}\left(\ker\left(P(V)'\otimes I_n\right)\right)\cap\Omega_{n,m,\kappa,\lambda}^{EU}$ for some $V\in O_m$, then the first $(n+1)\kappa+1$ Taylor series coefficients of $\tilde AV'-[\tilde BC]_+$ equate to zero and so, following the same argument as used in (i), $[\tilde BCV]_+=\tilde A$. By Theorem \ref{thm:onatski} then, $CV$ is the transfer function of the solution to $(\tilde B,\tilde A)$, which is clearly invertible. Thus, $(\tilde B,\tilde A)\in\Omega_{n,m,\kappa,\lambda}^{EUI}$ and by Theorem \ref{thm:oe} (i), $(\tilde B,\tilde A)\sim(B,A)$.

(iii) Lemma~\ref{lem:hankellrem} (i) implies that
\begin{align*}
\dim(\ker(P(V)'))&=\dim\left(\ker\left(\left[\begin{array}{cc}
-T' & I_{\kappa+1}\otimes V' \\
-H' & 0_{nm\kappa\times m(\kappa+1)}
\end{array}\right]\right)\right)\\
&=\dim(\ker(H'))\\
&=n(\kappa+\lambda+1)-\mathrm{rank}(H)\\
&=n(\kappa+\lambda+1)-\delta(C(z^{-1})-C(0)).
\end{align*}
By Theorem 4.2.15 of \cite{hj2} then
\begin{align*}
\dim(\ker(P(V)'\otimes I_n))&=n^2(\kappa+\lambda+1)-n\delta(C(z^{-1})-C(0)).
\end{align*}
To parametrize $\big((B,A)/\sim\big)\cap\;\Omega_{n,m,\kappa,\lambda}^{EUI}$ we may proceed as follows. By (ii) every point $(\tilde B,\tilde A)\in\big((B,A)/\sim\big)\cap\;\Omega_{n,m,\kappa,\lambda}^{EUI}$ belongs to some subspace $\mathrm{mat}_{\Omega_{n,m,\kappa,\lambda}}\left(\ker\left(P(V)'\otimes I_n\right)\right)$ for some $V\in O_m$. This $V$ is unique because if $(\tilde B,\tilde A)$ were to belong to two subspaces indexed by $V_1,V_2\in O_m$, then by (i), Theorem \ref{thm:oe} (i), and Theorem \ref{thm:onatski}, it would have to be the case that $CV_1=CV_2$ and, since $C_0$ is of full column rank, it would have to be the case that $V_1=V_2$. Thus, we can associate with every $(\tilde B,\tilde A)\in\big((B,A)/\sim\big)\cap\;\Omega_{n,m,\kappa,\lambda}^{EUI}$ a unique $V\in O_m$ that produces the first impulse response of its transfer function. Given this $V$, there is now a unique element of $\mathrm{mat}_{\Omega_{n,m,\kappa,\lambda}}\left(\ker\left(P(V)'\otimes I_n\right)\right)$ that produces $(\tilde B,\tilde A)$. Since $O_m$ is a differentiable manifold of dimension $\frac{1}{2}m(m-1)$ \cite[p.\ 22]{gp}, $\dim\left(\mathrm{mat}_{\Omega_{n,m,\kappa,\lambda}}\left(\ker\left(P(V)'\otimes I_n\right)\right)\right)=n^2(\kappa+\lambda+1)-n\delta(C(z^{-1})-C(0))$, and $\Omega_{n,m,\kappa,\lambda}^{EU}$ is an open set, our result follows.

Finally, $\delta(C(z^{-1})-C(0))$ is generically equal to $n\kappa$ by Lemma \ref{lem:hankellrem} (ii) and so the generic dimension of $\big((B,A)/\sim\big)\cap\;\Omega_{n,m,\kappa,\lambda}^{EUI}$ is $\frac{1}{2}m(m-1)+n^2(\lambda+1)$.
\end{proof}

\begin{proof}[Proof of Theorem \ref{thm:oeEUI0}]
Identical to the proof of Theorem \ref{thm:oeEUI}. Simply take $V=I_m$.
\end{proof}

\begin{proof}[Proof of Theorem \ref{thm:identEUI}]
(i) Suppose that $((\mathsf{X},\mathsf{Y}),\mathsf{V})=((B,A),I_m)$ is the unique solution to \eqref{eq:identcriterionEUI} in $\Omega_{n,m,\kappa,\lambda}^{EU}\times\mathbb{R}^{m\times m}$. If $(\tilde{B},\tilde{A})\in\big((B,A)/\sim\big)\cap\;\Omega^R\cap\Omega_{n,m,\kappa,\lambda}^{EUI}$ then Theorem \ref{thm:oeEUI} (i) implies that there exists a $V\in\mathbb{R}^{m\times m}$ such that $((\tilde{B},\tilde{A}),V)$ satisfies \eqref{eq:identcriterionEUI}. Uniqueness then implies that $((\tilde{B},\tilde{A}),V)=((B,A),I_m)$. Thus, $(B,A)$ is identified in $\Omega^R\cap\Omega_{n,m,\kappa,\lambda}^{EUI}$.

If there exists a solution $((\tilde{B},\tilde{A}),V)\neq((B,A),I_m)$ to \eqref{eq:identcriterionEUI}, then $(\tilde{B},\tilde{A})\in\Omega_{n,m,\kappa,\lambda}^{EU}$ and Theorem \ref{thm:oeEUI} (i) implies that $(\tilde{B},\tilde{A})\in\big((B,A)/\sim\big)\cap\;\Omega_{n,m,\kappa,\lambda}^{EUI0}$. The latter restriction equations of system \eqref{eq:identcriterionEUI} imply that $(\tilde{B},\tilde{A})\in\Omega^R$. It follows that $(\tilde{B},\tilde{A})\in\big((B,A)/\sim\big)\cap\;\Omega^R\cap\Omega_{n,m,\kappa,\lambda}^{EU}$ so that $(B,A)$ is not identified in $\Omega^R\cap\Omega_{n,m,\kappa,\lambda}^{EUI}$. 
\end{proof}

\begin{proof}[Proof of Theorem \ref{thm:localEUI}]
The result follows from Theorem 5.A.1 of \cite{fisher}, which characterizes isolated solutions to non-linear equations. Simply note that \eqref{eq:nonlinearcriterionEUI} can be expressed as
\begin{align*}
J(\mathrm{vec}_{\Omega_{n,m,\kappa,\lambda}}(\mathsf{dX},\mathsf{dY}),\mathrm{vec}_{\mathbb{R}^{m\times m}}(\mathsf{dV}))=0,
\end{align*}
where $J$ is the Jacobian of the mapping $\Omega_{n,m,\kappa,\lambda}\times \mathbb{R}^{m\times m}\rightarrow \mathbb{R}^{r+nm((n+1)\kappa+1)}\times \mathbb{R}^{m\times m}$,
\begin{align*}
((\mathsf{X},\mathsf{Y}),\mathsf{V})\mapsto\left(\left[\begin{array}{c}
\left(P(\mathsf{V})'\otimes I_n\right)\mathrm{vec}_{\Omega_{n,m,\kappa,\lambda}}(\mathsf{X},\mathsf{Y})\\
R\left(\mathrm{vec}_{\Omega_{n,m,\kappa,\lambda}}(\mathsf{X},\mathsf{Y})\right)
\end{array}\right]
-0_{r+nm((n+1)\kappa+1)\times 1}, \mathsf{V}'\mathsf{V}-I_m\right),
\end{align*}
evaluated at $((\mathsf{X},\mathsf{Y}),\mathsf{V})=((B,A),I_m)$. Note that the zeros in $\Omega_{n,m,\kappa,\lambda}^{EU}\times \mathbb{R}^{m\times m}$ of the mapping above are precisely $\big((B,A)/\sim\big)\cap\;\Omega^R\cap\Omega_{n,m,\kappa,\lambda}^{EUI}$.
\end{proof}

\begin{proof}[Proof of Theorem \ref{thm:identEUI0}]
Suppose that $(\mathsf{X},\mathsf{Y})=(B,A)$ is the unique solution to \eqref{eq:identcriterionEUI0} in $\Omega_{n,m,\kappa,\lambda}^{EU}$. If $(\tilde{B},\tilde{A})\in\big((B,A)/\sim\big)\cap\;\Omega^R\cap\Omega_{n,m,\kappa,\lambda}^{EUI0}$ then Theorem \ref{thm:oeEUI0} (i) implies that $(\tilde{B},\tilde{A})$ satisfies \eqref{eq:identcriterionEUI0}. Uniqueness then implies that $(\tilde{B},\tilde{A})=(B,A)$. Thus, $(B,A)$ is identified in $\Omega^R\cap\Omega_{n,m,\kappa,\lambda}^{EUI0}$.

If there exists a solution $(\tilde{B},\tilde{A})\neq(B,A)$ to \eqref{eq:identcriterionEUI0}, then $(\tilde{B},\tilde{A})\in\Omega_{n,m,\kappa,\lambda}^{EU}$ and Theorem \ref{thm:oeEUI0} (i) implies that $(\tilde{B},\tilde{A})\in\big((B,A)/\sim\big)\cap\;\Omega_{n,m,\kappa,\lambda}^{EUI0}$. The restriction equations of system \eqref{eq:identcriterionEUI0} imply that $(\tilde{B},\tilde{A})\in\Omega^R$. It follows that $(\tilde{B},\tilde{A})\in\big((B,A)/\sim\big)\cap\;\Omega^R\cap\Omega_{n,m,\kappa,\lambda}^{EU}$ so that $(B,A)$ is not identified in $\Omega^R\cap\Omega_{n,m,\kappa,\lambda}^{EUI0}$. 
\end{proof}

\begin{proof}[Proof of Theorem \ref{thm:localEUI0}]
The result follows again from Theorem 5.A.1 of \cite{fisher}. Equations \eqref{eq:nonlinearcriterionEUI0} can be expressed as
\begin{align*}
J(\mathrm{vec}_{\Omega_{n,m,\kappa,\lambda}}(\mathsf{dX},\mathsf{dY}))=0,
\end{align*}
where $J$ is the Jacobian of the mapping $\Omega_{n,m,\kappa,\lambda}\rightarrow \mathbb{R}^{r+nm((n+1)\kappa+1)}$,
\begin{align*}
(\mathsf{X},\mathsf{Y})\mapsto\left(\left[\begin{array}{c}
\left(P(I_m)'\otimes I_n\right)\mathrm{vec}_{\Omega_{n,m,\kappa,\lambda}}(\mathsf{X},\mathsf{Y})\\
R\left(\mathrm{vec}_{\Omega_{n,m,\kappa,\lambda}}(\mathsf{X},\mathsf{Y})\right)
\end{array}\right]
-0_{r+nm((n+1)\kappa+1)\times 1}\right),
\end{align*}
evaluated at $(\mathsf{X},\mathsf{Y})=(B,A)$. Clearly, the zeros in $\Omega_{n,m,\kappa,\lambda}^{EU}\times \mathbb{R}^{m\times m}$ of the mapping above are precisely $\big((B,A)/\sim\big)\cap\;\Omega^R\cap\Omega_{n,m,\kappa,\lambda}^{EUI0}$.
\end{proof}

\begin{proof}[Proof of Theorem \ref{thm:identEUIi}]
Suppose that $((\mathsf{x},\mathsf{y}),\mathsf{V})=(e_i'(B,A),I_m)$ is the unique solution to \eqref{eq:identcriterionEUIi} in $e_i'\Omega_{n,m,\kappa,\lambda}^{EU}\times\mathbb{R}^{m\times m}$. If $(\tilde{B},\tilde{A})\in\big((B,A)/\sim\big)\cap\;\Omega^{R_i}\cap\Omega_{n,m,\kappa,\lambda}^{EUI}$ then Theorem \ref{thm:oeEUI} (i) implies that there exists a $V\in\mathbb{R}^{m\times m}$ such that $((\tilde{B},\tilde{A}),V)$ satisfies \eqref{eq:identcriterionEUIi}. In particular,
\begin{align*}
P(V)'\mathrm{evec}_{\Omega_{n,m,\kappa,\lambda}}(e_i'(\tilde{B},\tilde{A}))&=P(V)'(I_{n(\kappa+\lambda+1)+m(\kappa+1)}\otimes e_i')\mathrm{vec}_{\Omega_{n,m,\kappa,\lambda}}(\tilde{B},\tilde{A})\\
&=(P(V)'\otimes e_i')\mathrm{vec}_{\Omega_{n,m,\kappa,\lambda}}(\tilde{B},\tilde{A})\\
&=(I_{m((n+1)\kappa+1)}\otimes e_i')(P(V)'\otimes I_n)\mathrm{vec}_{\Omega_{n,m,\kappa,\lambda}}(\tilde{B},\tilde{A})\\
&=0_{m((n+1)\kappa+1)\times 1}.
\end{align*}
Uniqueness then implies that $(e_i'(\tilde{B},\tilde{A}),V)=(e_i'(B,A),I_m)$. Thus, the $i$-th equation of $(B,A)$ is identified in $\Omega^{R_i}\cap\Omega_{n,m,\kappa,\lambda}^{EUI}$.

If there exists a solution $(e_i'(\tilde{B},\tilde{A}),V)\neq(e_i'(B,A),I_m)$ to \eqref{eq:identcriterionEUIi}, then $(\tilde{B},\tilde{A})\in\Omega_{n,m,\kappa,\lambda}^{EU}$ and Theorem \ref{thm:oeEUI} (i) implies that $e_i'(\tilde{B},\tilde{A})\in e_i'\left(\big((B,A)/\sim\big)\cap\;\Omega_{n,m,\kappa,\lambda}^{EUI}\right)$. The restriction equations of system \eqref{eq:identcriterionEUIi} imply that $(\tilde{B},\tilde{A})\in\Omega^R$. It follows that $e_i'\left(\tilde{B},\tilde{A}\right)\in e_i'\left(\big((B,A)/\sim\big)\cap\;\Omega^R\cap\Omega_{n,m,\kappa,\lambda}^{EU}\right)$ so that $(B,A)$ is not identified in $\Omega^R\cap\Omega_{n,m,\kappa,\lambda}^{EUI}$. 
\end{proof}

\begin{proof}[Proof of Theorem \ref{thm:localEUIi}]
The result follows, again, from Theorem 5.A.1 of \cite{fisher} because \eqref{eq:nonlinearcriterionEUIi} can be expressed as
\begin{align*}
J_i(\mathrm{evec}_{\Omega_{n,m,\kappa,\lambda}}(\mathsf{dx},\mathsf{dy}),\mathrm{vec}_{\mathbb{R}^{m\times m}}(\mathsf{dV}))=0,
\end{align*}
where $J_i$ is the Jacobian of the mapping $e_i'\Omega_{n,m,\kappa,\lambda}\times\mathbb{R}^{m\times m}\rightarrow \mathbb{R}^{r+m((n+1)\kappa+1)}\times \mathbb{R}^{m\times m}$,
\begin{align*}
((\mathsf{x},\mathsf{y}),\mathsf{V})\mapsto\left(\left[\begin{array}{c}
P(\mathsf{V})'\mathrm{evec}_{\Omega_{n,m,\kappa,\lambda}}(\mathsf{x},\mathsf{y})\\
R_i\left(\mathrm{evec}_{\Omega_{n,m,\kappa,\lambda}}(\mathsf{x},\mathsf{y})\right)
\end{array}\right]
-0_{r+m((n+1)\kappa+1)\times 1}, \mathsf{V}'\mathsf{V}-I_m\right),
\end{align*}
evaluated at $((\mathsf{x},\mathsf{y}),\mathsf{V})=(e_i'(B,A),I_m)$. The zeros in $e_i'\Omega_{n,m,\kappa,\lambda}^{EU}\times \mathbb{R}^{m\times m}$ of the mapping above are precisely $e_i'\left(\big((B,A)/\sim\big)\cap\;\Omega^R\cap\Omega_{n,m,\kappa,\lambda}^{EUI}\right)$.
\end{proof}

\begin{proof}[Proof of Theorem \ref{thm:identEUI0i}]
Suppose that $(\mathsf{x},\mathsf{y})=e_i'(B,A)$ is the unique solution to \eqref{eq:identcriterionEUI0i} in $e_i'\Omega_{n,m,\kappa,\lambda}^{EU}$. If $(\tilde{B},\tilde{A})\in\big((B,A)/\sim\big)\cap\;\Omega^{R_i}\cap\Omega_{n,m,\kappa,\lambda}^{EUI0}$ then Theorem \ref{thm:oeEUI0} (i) implies that $(\tilde{B},\tilde{A})$ satisfies \eqref{eq:identcriterionEUI0i}. In particular,
\begin{align*}
P(I_m)'\mathrm{evec}_{\Omega_{n,m,\kappa,\lambda}}(e_i'(\tilde{B},\tilde{A}))&=P(I_m)'(I_{n(\kappa+\lambda+1)+m(\kappa+1)}\otimes e_i')\mathrm{vec}_{\Omega_{n,m,\kappa,\lambda}}(\tilde{B},\tilde{A})\\
&=(P(I_m)'\otimes e_i')\mathrm{vec}_{\Omega_{n,m,\kappa,\lambda}}(\tilde{B},\tilde{A})\\
&=(I_{m((n+1)\kappa+1)}\otimes e_i')(P(I_m)'\otimes I_n)\mathrm{vec}_{\Omega_{n,m,\kappa,\lambda}}(\tilde{B},\tilde{A})\\
&=0_{m((n+1)\kappa+1)\times 1}.
\end{align*}
Uniqueness then implies that $e_i'(\tilde{B},\tilde{A})=e_i'(B,A)$. Thus, the $i$-th equation of $(B,A)$ is identified in $\Omega^{R_i}\cap\Omega_{n,m,\kappa,\lambda}^{EUI0}$.

If there exists a solution $e_i'(\tilde{B},\tilde{A})\neq e_i'(B,A)$ to \eqref{eq:identcriterionEUI0i}, then $(\tilde{B},\tilde{A})\in\Omega_{n,m,\kappa,\lambda}^{EU}$ and Theorem \ref{thm:oeEUI0} (i) implies that $e_i'(\tilde{B},\tilde{A})\in e_i'\left(\big((B,A)/\sim\big)\cap\;\Omega_{n,m,\kappa,\lambda}^{EUI}\right)$. The restriction equations of system \eqref{eq:identcriterionEUI0i} imply that $(\tilde{B},\tilde{A})\in\Omega^R$. It follows that $e_i'\left(\tilde{B},\tilde{A}\right)\in e_i'\left(\big((B,A)/\sim\big)\cap\;\Omega^R\cap\Omega_{n,m,\kappa,\lambda}^{EU}\right)$ so that $(B,A)$ is not identified in $\Omega^R\cap\Omega_{n,m,\kappa,\lambda}^{EUI}$. 
\end{proof}

\begin{proof}[Proof of Theorem \ref{thm:localEUI0i}]
The result follows, again, from Theorem 5.A.1 of \cite{fisher} because \eqref{eq:nonlinearcriterionEUI0i} can be expressed as
\begin{align*}
J_i(\mathrm{evec}_{\Omega_{n,m,\kappa,\lambda}}(\mathsf{dx},\mathsf{dy}))=0,
\end{align*}
where $J_i$ is the Jacobian of the mapping $e_i'\Omega_{n,m,\kappa,\lambda}\rightarrow \mathbb{R}^{r+m((n+1)\kappa+1)}$,
\begin{align*}
(\mathsf{x},\mathsf{y})\mapsto\left(\left[\begin{array}{c}
P(\mathsf{V})'\mathrm{evec}_{\Omega_{n,m,\kappa,\lambda}}(\mathsf{x},\mathsf{y})\\
R_i\left(\mathrm{evec}_{\Omega_{n,m,\kappa,\lambda}}(\mathsf{x},\mathsf{y})\right)
\end{array}\right]
-0_{r+m((n+1)\kappa+1)\times 1}\right),
\end{align*}
evaluated at $(\mathsf{x},\mathsf{y})=e_i'(B,A)$. The zeros in $e_i'\Omega_{n,m,\kappa,\lambda}^{EU}\times \mathbb{R}^{m\times m}$ of the mapping above are precisely $e_i'\left(\big((B,A)/\sim\big)\cap\;\Omega^R\cap\Omega_{n,m,\kappa,\lambda}^{EUI0}\right)$.
\end{proof}

\bibliographystyle{apacite}
\bibliography{ident}

\begin{thebibliography}{}

\bibitem [\protect \citeauthoryear {%
Al-Sadoon%
}{%
Al-Sadoon%
}{%
{\protect \APACyear {2018}}%
}]{%
linsys}
\APACinsertmetastar {%
linsys}%
\begin{APACrefauthors}%
Al-Sadoon, M\BPBI M.%
\end{APACrefauthors}%
\unskip\
\newblock
\APACrefYearMonthDay{2018}{June}{}.
\newblock
{\BBOQ}\APACrefatitle {{The Linear Systems Approach to Linear Rational Expectations Models}} {{The Linear Systems Approach to Linear Rational Expectations Models}}.{\BBCQ}
\newblock
\APACjournalVolNumPages{Econometric Theory}{34}{03}{628-658}.
\PrintBackRefs{\CurrentBib}

\bibitem [\protect \citeauthoryear {%
Al-Sadoon%
}{%
Al-Sadoon%
}{%
{\protect \APACyear {2024}}%
}]{%
spectral}
\APACinsertmetastar {%
spectral}%
\begin{APACrefauthors}%
Al-Sadoon, M\BPBI M.%
\end{APACrefauthors}%
\unskip\
\newblock
\APACrefYearMonthDay{2024}{}{}.
\newblock
{\BBOQ}\APACrefatitle {The Spectral Approach to Linear Rational Expectations Models} {The spectral approach to linear rational expectations models}.{\BBCQ}
\newblock
\APACjournalVolNumPages{Econometric Theory}{}{}{}.
\newblock
\APACrefnote{Forthcoming}
\PrintBackRefs{\CurrentBib}

\bibitem [\protect \citeauthoryear {%
Anderson%
, Deistler%
, Chen%
\BCBL {}\ \BBA {} Filler%
}{%
Anderson%
\ \protect \BOthers {.}}{%
{\protect \APACyear {2012}}%
}]{%
andersonetal2012}
\APACinsertmetastar {%
andersonetal2012}%
\begin{APACrefauthors}%
Anderson, B\BPBI D.%
, Deistler, M.%
, Chen, W.%
\BCBL {}\ \BBA {} Filler, A.%
\end{APACrefauthors}%
\unskip\
\newblock
\APACrefYearMonthDay{2012}{}{}.
\newblock
{\BBOQ}\APACrefatitle {{Autoregressive Models of Singular Spectral Matrices}} {{Autoregressive Models of Singular Spectral Matrices}}.{\BBCQ}
\newblock
\APACjournalVolNumPages{Automatica}{48}{11}{2843 - 2849}.
\PrintBackRefs{\CurrentBib}

\bibitem [\protect \citeauthoryear {%
Anderson%
, Deistler%
, Felsenstein%
\BCBL {}\ \BBA {} Koelbl%
}{%
Anderson%
\ \protect \BOthers {.}}{%
{\protect \APACyear {2016}}%
}]{%
andersonetal}
\APACinsertmetastar {%
andersonetal}%
\begin{APACrefauthors}%
Anderson, B\BPBI D.%
, Deistler, M.%
, Felsenstein, E.%
\BCBL {}\ \BBA {} Koelbl, L.%
\end{APACrefauthors}%
\unskip\
\newblock
\APACrefYearMonthDay{2016}{}{}.
\newblock
{\BBOQ}\APACrefatitle {{The Structure of Multivariate AR and ARMA Systems: Regular and Singular Systems; the Single and the Mixed Frequency Case}} {{The Structure of Multivariate AR and ARMA Systems: Regular and Singular Systems; the Single and the Mixed Frequency Case}}.{\BBCQ}
\newblock
\APACjournalVolNumPages{Journal of Econometrics}{192}{2}{366 - 373}.
\PrintBackRefs{\CurrentBib}

\bibitem [\protect \citeauthoryear {%
Bacchiocchi%
\ \BBA {} Kitagawa%
}{%
Bacchiocchi%
\ \BBA {} Kitagawa%
}{%
{\protect \APACyear {2019}}%
}]{%
darkside}
\APACinsertmetastar {%
darkside}%
\begin{APACrefauthors}%
Bacchiocchi, E.%
\BCBT {}\ \BBA {} Kitagawa, T.%
\end{APACrefauthors}%
\unskip\
\newblock
\APACrefYearMonthDay{2019}{}{}.
\newblock
\APACrefbtitle {{The Dark Side of the SVAR: a Trip Into the Local Identification World}} {{The Dark Side of the SVAR: a Trip Into the Local Identification World}}\ \APACbVolEdTR {}{mimeo}.
\PrintBackRefs{\CurrentBib}

\bibitem [\protect \citeauthoryear {%
Baggio%
\ \BBA {} Ferrante%
}{%
Baggio%
\ \BBA {} Ferrante%
}{%
{\protect \APACyear {2016}}%
}]{%
baggioferrante}
\APACinsertmetastar {%
baggioferrante}%
\begin{APACrefauthors}%
Baggio, G.%
\BCBT {}\ \BBA {} Ferrante, A.%
\end{APACrefauthors}%
\unskip\
\newblock
\APACrefYearMonthDay{2016}{April}{}.
\newblock
{\BBOQ}\APACrefatitle {On the Factorization of Rational Discrete-Time Spectral Densities} {On the factorization of rational discrete-time spectral densities}.{\BBCQ}
\newblock
\APACjournalVolNumPages{IEEE Transactions on Automatic Control}{61}{4}{969-981}.
\PrintBackRefs{\CurrentBib}

\bibitem [\protect \citeauthoryear {%
Blanchard%
}{%
Blanchard%
}{%
{\protect \APACyear {2018}}%
}]{%
blanchard}
\APACinsertmetastar {%
blanchard}%
\begin{APACrefauthors}%
Blanchard, O.%
\end{APACrefauthors}%
\unskip\
\newblock
\APACrefYearMonthDay{2018}{01}{}.
\newblock
{\BBOQ}\APACrefatitle {{On the Future of Macroeconomic Models}} {{On the Future of Macroeconomic Models}}.{\BBCQ}
\newblock
\APACjournalVolNumPages{Oxford Review of Economic Policy}{34}{1-2}{43-54}.
\PrintBackRefs{\CurrentBib}

\bibitem [\protect \citeauthoryear {%
Bollen%
}{%
Bollen%
}{%
{\protect \APACyear {1989}}%
}]{%
bollen}
\APACinsertmetastar {%
bollen}%
\begin{APACrefauthors}%
Bollen, K\BPBI A.%
\end{APACrefauthors}%
\unskip\
\newblock
\APACrefYear{1989}.
\newblock
\APACrefbtitle {Structural Equations with Latent Variables} {Structural equations with latent variables}.
\newblock
\APACaddressPublisher{}{John Wiley \& Sons, Inc., New York}.
\newblock
\begin{APACrefURL} \url{https://doi.org/10.1002/9781118619179} \end{APACrefURL}
\newblock
\APACrefnote{A Wiley-Interscience Publication}
\newblock
\begin{APACrefDOI} \doi{10.1002/9781118619179} \end{APACrefDOI}
\PrintBackRefs{\CurrentBib}

\bibitem [\protect \citeauthoryear {%
Brockwell%
\ \BBA {} Davis%
}{%
Brockwell%
\ \BBA {} Davis%
}{%
{\protect \APACyear {1991}}%
}]{%
bd}
\APACinsertmetastar {%
bd}%
\begin{APACrefauthors}%
Brockwell, P\BPBI J.%
\BCBT {}\ \BBA {} Davis, R\BPBI A.%
\end{APACrefauthors}%
\unskip\
\newblock
\APACrefYear{1991}.
\newblock
\APACrefbtitle {Time Series: Theory and Methods, 2nd Edition} {Time series: Theory and methods, 2nd edition}.
\newblock
\APACaddressPublisher{New York, NY. USA}{Springer}.
\PrintBackRefs{\CurrentBib}

\bibitem [\protect \citeauthoryear {%
Canova%
}{%
Canova%
}{%
{\protect \APACyear {2011}}%
}]{%
canova}
\APACinsertmetastar {%
canova}%
\begin{APACrefauthors}%
Canova, F.%
\end{APACrefauthors}%
\unskip\
\newblock
\APACrefYear{2011}.
\newblock
\APACrefbtitle {Methods for Applied Macroeconomic Research} {Methods for applied macroeconomic research}.
\newblock
\APACaddressPublisher{Princeton, USA}{Princeton University Press}.
\PrintBackRefs{\CurrentBib}

\bibitem [\protect \citeauthoryear {%
Canova%
\ \BBA {} Sala%
}{%
Canova%
\ \BBA {} Sala%
}{%
{\protect \APACyear {2009}}%
}]{%
canovasala}
\APACinsertmetastar {%
canovasala}%
\begin{APACrefauthors}%
Canova, F.%
\BCBT {}\ \BBA {} Sala, L.%
\end{APACrefauthors}%
\unskip\
\newblock
\APACrefYearMonthDay{2009}{May}{}.
\newblock
{\BBOQ}\APACrefatitle {{Back to Square One: Identification Issues in DSGE Models}} {{Back to Square One: Identification Issues in DSGE Models}}.{\BBCQ}
\newblock
\APACjournalVolNumPages{Journal of Monetary Economics}{56}{4}{431-449}.
\PrintBackRefs{\CurrentBib}

\bibitem [\protect \citeauthoryear {%
Clancey%
\ \BBA {} Gohberg%
}{%
Clancey%
\ \BBA {} Gohberg%
}{%
{\protect \APACyear {1981}}%
}]{%
cg}
\APACinsertmetastar {%
cg}%
\begin{APACrefauthors}%
Clancey, K\BPBI F.%
\BCBT {}\ \BBA {} Gohberg, I.%
\end{APACrefauthors}%
\unskip\
\newblock
\APACrefYear{1981}.
\newblock
\APACrefbtitle {Factorization of Matrix Functions and Singular Integral Operators} {Factorization of matrix functions and singular integral operators}.
\newblock
\APACaddressPublisher{Boston, USA}{Birkh\"{a}user Verlag Basel}.
\PrintBackRefs{\CurrentBib}

\bibitem [\protect \citeauthoryear {%
Cochrane%
}{%
Cochrane%
}{%
{\protect \APACyear {2011}}%
}]{%
cochrane}
\APACinsertmetastar {%
cochrane}%
\begin{APACrefauthors}%
Cochrane, J\BPBI H.%
\end{APACrefauthors}%
\unskip\
\newblock
\APACrefYearMonthDay{2011}{}{}.
\newblock
{\BBOQ}\APACrefatitle {Determinacy and Identification with Taylor Rules} {Determinacy and identification with taylor rules}.{\BBCQ}
\newblock
\APACjournalVolNumPages{Journal of Political Economy}{119}{3}{565--615}.
\PrintBackRefs{\CurrentBib}

\bibitem [\protect \citeauthoryear {%
Deistler%
}{%
Deistler%
}{%
{\protect \APACyear {1975}}%
}]{%
deistler75}
\APACinsertmetastar {%
deistler75}%
\begin{APACrefauthors}%
Deistler, M.%
\end{APACrefauthors}%
\unskip\
\newblock
\APACrefYearMonthDay{1975}{}{}.
\newblock
{\BBOQ}\APACrefatitle {z-Transform and identification of linear econometric models with autocorrelated errors} {z-transform and identification of linear econometric models with autocorrelated errors}.{\BBCQ}
\newblock
\APACjournalVolNumPages{Metrika}{22}{}{13--25}.
\PrintBackRefs{\CurrentBib}

\bibitem [\protect \citeauthoryear {%
Deistler%
}{%
Deistler%
}{%
{\protect \APACyear {1976}}%
}]{%
deistler76}
\APACinsertmetastar {%
deistler76}%
\begin{APACrefauthors}%
Deistler, M.%
\end{APACrefauthors}%
\unskip\
\newblock
\APACrefYearMonthDay{1976}{}{}.
\newblock
{\BBOQ}\APACrefatitle {The Identifiability of Linear Econometric Models with Autocorrelated Errors} {The identifiability of linear econometric models with autocorrelated errors}.{\BBCQ}
\newblock
\APACjournalVolNumPages{International Economic Review}{17}{1}{26-46}.
\PrintBackRefs{\CurrentBib}

\bibitem [\protect \citeauthoryear {%
Deistler%
}{%
Deistler%
}{%
{\protect \APACyear {1978}}%
}]{%
deistler78}
\APACinsertmetastar {%
deistler78}%
\begin{APACrefauthors}%
Deistler, M.%
\end{APACrefauthors}%
\unskip\
\newblock
\APACrefYearMonthDay{1978}{}{}.
\newblock
{\BBOQ}\APACrefatitle {The structural identifiability of linear models with autocorrelated errors in the case of cross-equation restrictions} {The structural identifiability of linear models with autocorrelated errors in the case of cross-equation restrictions}.{\BBCQ}
\newblock
\APACjournalVolNumPages{Journal of Econometrics}{8}{1}{23-31}.
\PrintBackRefs{\CurrentBib}

\bibitem [\protect \citeauthoryear {%
Deistler%
}{%
Deistler%
}{%
{\protect \APACyear {1983}}%
}]{%
deistler83}
\APACinsertmetastar {%
deistler83}%
\begin{APACrefauthors}%
Deistler, M.%
\end{APACrefauthors}%
\unskip\
\newblock
\APACrefYearMonthDay{1983}{}{}.
\newblock
{\BBOQ}\APACrefatitle {The Properties of the Parameterization of Armax Systems and Their Relevance for Structural Estimation and Dynamic Specification} {The properties of the parameterization of armax systems and their relevance for structural estimation and dynamic specification}.{\BBCQ}
\newblock
\APACjournalVolNumPages{Econometrica}{51}{4}{1187--1207}.
\PrintBackRefs{\CurrentBib}

\bibitem [\protect \citeauthoryear {%
Deistler%
\ \BBA {} Scherrer%
}{%
Deistler%
\ \BBA {} Scherrer%
}{%
{\protect \APACyear {2022}}%
}]{%
ds}
\APACinsertmetastar {%
ds}%
\begin{APACrefauthors}%
Deistler, M.%
\BCBT {}\ \BBA {} Scherrer, W.%
\end{APACrefauthors}%
\unskip\
\newblock
\APACrefYear{2022}.
\newblock
\APACrefbtitle {Time Series Models} {Time series models}\ (\BNUM~224).
\newblock
\APACaddressPublisher{Cham, Switzerland}{Springer Nature}.
\PrintBackRefs{\CurrentBib}

\bibitem [\protect \citeauthoryear {%
Deistler%
\ \BBA {} Schrader%
}{%
Deistler%
\ \BBA {} Schrader%
}{%
{\protect \APACyear {1979}}%
}]{%
ds79}
\APACinsertmetastar {%
ds79}%
\begin{APACrefauthors}%
Deistler, M.%
\BCBT {}\ \BBA {} Schrader, J.%
\end{APACrefauthors}%
\unskip\
\newblock
\APACrefYearMonthDay{1979}{}{}.
\newblock
{\BBOQ}\APACrefatitle {{Linear Models with Autocorrelated Errors: Structural Identifiability in the Absence of Minimality Assumptions}} {{Linear Models with Autocorrelated Errors: Structural Identifiability in the Absence of Minimality Assumptions}}.{\BBCQ}
\newblock
\APACjournalVolNumPages{Econometrica}{47}{2}{495-504}.
\PrintBackRefs{\CurrentBib}

\bibitem [\protect \citeauthoryear {%
Deistler%
\ \BBA {} Wang%
}{%
Deistler%
\ \BBA {} Wang%
}{%
{\protect \APACyear {1989}}%
}]{%
deistlerwang}
\APACinsertmetastar {%
deistlerwang}%
\begin{APACrefauthors}%
Deistler, M.%
\BCBT {}\ \BBA {} Wang, L.%
\end{APACrefauthors}%
\unskip\
\newblock
\APACrefYearMonthDay{1989}{}{}.
\newblock
{\BBOQ}\APACrefatitle {The common structure of parametrizations for linear systems} {The common structure of parametrizations for linear systems}.{\BBCQ}
\newblock
\APACjournalVolNumPages{Linear Algebra and its Applications}{122}{}{921--941}.
\PrintBackRefs{\CurrentBib}

\bibitem [\protect \citeauthoryear {%
DeJong%
\ \BBA {} Dave%
}{%
DeJong%
\ \BBA {} Dave%
}{%
{\protect \APACyear {2011}}%
}]{%
dd}
\APACinsertmetastar {%
dd}%
\begin{APACrefauthors}%
DeJong, D.%
\BCBT {}\ \BBA {} Dave, C.%
\end{APACrefauthors}%
\unskip\
\newblock
\APACrefYear{2011}.
\newblock
\APACrefbtitle {Structural Macroeconometrics: (Second Edition)} {Structural macroeconometrics: (second edition)}.
\newblock
\APACaddressPublisher{Princeton, USA}{Princeton University Press}.
\PrintBackRefs{\CurrentBib}

\bibitem [\protect \citeauthoryear {%
Dhrymes%
}{%
Dhrymes%
}{%
{\protect \APACyear {1994}}%
}]{%
dhrymesbook}
\APACinsertmetastar {%
dhrymesbook}%
\begin{APACrefauthors}%
Dhrymes, P\BPBI J.%
\end{APACrefauthors}%
\unskip\
\newblock
\APACrefYear{1994}.
\newblock
\APACrefbtitle {Topics In Advanced Econometrics: Volume II Linear and Nonlinear Simultaneous Equations} {Topics in advanced econometrics: Volume ii linear and nonlinear simultaneous equations}.
\newblock
\APACaddressPublisher{New York, USA}{Springer-Verlag}.
\PrintBackRefs{\CurrentBib}

\bibitem [\protect \citeauthoryear {%
Dufour%
\ \BBA {} Renault%
}{%
Dufour%
\ \BBA {} Renault%
}{%
{\protect \APACyear {1998}}%
}]{%
dr}
\APACinsertmetastar {%
dr}%
\begin{APACrefauthors}%
Dufour, J\BHBI M.%
\BCBT {}\ \BBA {} Renault, E.%
\end{APACrefauthors}%
\unskip\
\newblock
\APACrefYearMonthDay{1998}{September}{}.
\newblock
{\BBOQ}\APACrefatitle {{Short Run and Long Run Causality in Time Series: Theory}} {{Short Run and Long Run Causality in Time Series: Theory}}.{\BBCQ}
\newblock
\APACjournalVolNumPages{Econometrica}{66}{5}{1099-1125}.
\PrintBackRefs{\CurrentBib}

\bibitem [\protect \citeauthoryear {%
Evard%
\ \BBA {} Jafari%
}{%
Evard%
\ \BBA {} Jafari%
}{%
{\protect \APACyear {1994}}%
}]{%
components}
\APACinsertmetastar {%
components}%
\begin{APACrefauthors}%
Evard, J\BHBI C.%
\BCBT {}\ \BBA {} Jafari, F.%
\end{APACrefauthors}%
\unskip\
\newblock
\APACrefYearMonthDay{1994}{}{}.
\newblock
{\BBOQ}\APACrefatitle {The Set of All $m \times n$ Rectangular Real Matrices of Rank $r$ is Connected by Analytic Regular Arcs} {The set of all $m \times n$ rectangular real matrices of rank $r$ is connected by analytic regular arcs}.{\BBCQ}
\newblock
\APACjournalVolNumPages{Proceedings of the American Mathematical Society}{120}{2}{413--419}.
\PrintBackRefs{\CurrentBib}

\bibitem [\protect \citeauthoryear {%
Fisher%
}{%
Fisher%
}{%
{\protect \APACyear {1966}}%
}]{%
fisher}
\APACinsertmetastar {%
fisher}%
\begin{APACrefauthors}%
Fisher, F.%
\end{APACrefauthors}%
\unskip\
\newblock
\APACrefYear{1966}.
\newblock
\APACrefbtitle {The Identification Problem in Econometrics} {The identification problem in econometrics}.
\newblock
\APACaddressPublisher{New York, USA}{McGraw-Hill}.
\PrintBackRefs{\CurrentBib}

\bibitem [\protect \citeauthoryear {%
Fuller%
}{%
Fuller%
}{%
{\protect \APACyear {1987}}%
}]{%
fuller}
\APACinsertmetastar {%
fuller}%
\begin{APACrefauthors}%
Fuller, W\BPBI A.%
\end{APACrefauthors}%
\unskip\
\newblock
\APACrefYear{1987}.
\newblock
\APACrefbtitle {Measurement Error Models} {Measurement error models}.
\newblock
\APACaddressPublisher{New York, USA}{John Wiley \& Sons, Inc.}
\PrintBackRefs{\CurrentBib}

\bibitem [\protect \citeauthoryear {%
I.~Gohberg%
, Lancaster%
\BCBL {}\ \BBA {} Rodman%
}{%
I.~Gohberg%
\ \protect \BOthers {.}}{%
{\protect \APACyear {2006}}%
}]{%
glr}
\APACinsertmetastar {%
glr}%
\begin{APACrefauthors}%
Gohberg, I.%
, Lancaster, P.%
\BCBL {}\ \BBA {} Rodman, L.%
\end{APACrefauthors}%
\unskip\
\newblock
\APACrefYear{2006}.
\newblock
\APACrefbtitle {Invariant Subspaces of Matrices with Applications} {Invariant subspaces of matrices with applications}.
\newblock
\APACaddressPublisher{Philadelphia, USA}{SIAM}.
\PrintBackRefs{\CurrentBib}

\bibitem [\protect \citeauthoryear {%
I\BPBI C.~Gohberg%
\ \BBA {} Fel'dman%
}{%
I\BPBI C.~Gohberg%
\ \BBA {} Fel'dman%
}{%
{\protect \APACyear {1974}}%
}]{%
gf}
\APACinsertmetastar {%
gf}%
\begin{APACrefauthors}%
Gohberg, I\BPBI C.%
\BCBT {}\ \BBA {} Fel'dman, I\BPBI A.%
\end{APACrefauthors}%
\unskip\
\newblock
\APACrefYear{1974}.
\newblock
\APACrefbtitle {Convolution Equations and Projection Methods for Their Solution} {Convolution equations and projection methods for their solution}\ (\BVOL~41).
\newblock
\APACaddressPublisher{Providence, USA}{American Mathematical Society}.
\PrintBackRefs{\CurrentBib}

\bibitem [\protect \citeauthoryear {%
Guillemin%
\ \BBA {} Pollack%
}{%
Guillemin%
\ \BBA {} Pollack%
}{%
{\protect \APACyear {1974}}%
}]{%
gp}
\APACinsertmetastar {%
gp}%
\begin{APACrefauthors}%
Guillemin, V.%
\BCBT {}\ \BBA {} Pollack, A.%
\end{APACrefauthors}%
\unskip\
\newblock
\APACrefYear{1974}.
\newblock
\APACrefbtitle {Differential Topology} {Differential topology}.
\newblock
\APACaddressPublisher{Englewood Cliffs, NJ, USA}{Prentice Hall, Inc.}
\PrintBackRefs{\CurrentBib}

\bibitem [\protect \citeauthoryear {%
Hannan%
}{%
Hannan%
}{%
{\protect \APACyear {1969}}%
}]{%
hannan69}
\APACinsertmetastar {%
hannan69}%
\begin{APACrefauthors}%
Hannan, E\BPBI J.%
\end{APACrefauthors}%
\unskip\
\newblock
\APACrefYearMonthDay{1969}{}{}.
\newblock
{\BBOQ}\APACrefatitle {The Identification of Vector Mixed Autoregressive-Moving Average Systems} {The identification of vector mixed autoregressive-moving average systems}.{\BBCQ}
\newblock
\APACjournalVolNumPages{Biometrika}{56}{1}{223--225}.
\PrintBackRefs{\CurrentBib}

\bibitem [\protect \citeauthoryear {%
Hannan%
}{%
Hannan%
}{%
{\protect \APACyear {1971}}%
}]{%
hannan71}
\APACinsertmetastar {%
hannan71}%
\begin{APACrefauthors}%
Hannan, E\BPBI J.%
\end{APACrefauthors}%
\unskip\
\newblock
\APACrefYearMonthDay{1971}{}{}.
\newblock
{\BBOQ}\APACrefatitle {{The Identification Problem for Multiple Equation Systems with Moving Average Errors}} {{The Identification Problem for Multiple Equation Systems with Moving Average Errors}}.{\BBCQ}
\newblock
\APACjournalVolNumPages{Econometrica}{39}{5}{751-765}.
\PrintBackRefs{\CurrentBib}

\bibitem [\protect \citeauthoryear {%
Hannan%
\ \BBA {} Deistler%
}{%
Hannan%
\ \BBA {} Deistler%
}{%
{\protect \APACyear {1988}}%
}]{%
hd}
\APACinsertmetastar {%
hd}%
\begin{APACrefauthors}%
Hannan, E\BPBI J.%
\BCBT {}\ \BBA {} Deistler, M.%
\end{APACrefauthors}%
\unskip\
\newblock
\APACrefYear{1988}.
\newblock
\APACrefbtitle {The Statistical Theory of Linear Systems} {The statistical theory of linear systems}.
\newblock
\APACaddressPublisher{Philadelphia, PA, USA}{Society for Industrial and Applied Mathematics}.
\newblock
\APACrefnote{(reprinted 2012)}
\PrintBackRefs{\CurrentBib}

\bibitem [\protect \citeauthoryear {%
Hansen%
\ \BBA {} Sargent%
}{%
Hansen%
\ \BBA {} Sargent%
}{%
{\protect \APACyear {1981}}%
}]{%
hansensargent}
\APACinsertmetastar {%
hansensargent}%
\begin{APACrefauthors}%
Hansen, L\BPBI P.%
\BCBT {}\ \BBA {} Sargent, T\BPBI J.%
\end{APACrefauthors}%
\unskip\
\newblock
\APACrefYearMonthDay{1981}{}{}.
\newblock
{\BBOQ}\APACrefatitle {{Formulating and Estimating Dynamic Linear Rational Expectations Models}} {{Formulating and Estimating Dynamic Linear Rational Expectations Models}}.{\BBCQ}
\newblock
\BIn{} R\BPBI E.~Lucas\ \BBA {} T\BPBI J.~Sargent\ (\BEDS), \APACrefbtitle {{Rational Expectations and Econometric Practice}} {{Rational Expectations and Econometric Practice}}\ (\BPGS\ 91--125).
\newblock
\APACaddressPublisher{Minneapolis}{The University of Minnesota Press}.
\PrintBackRefs{\CurrentBib}

\bibitem [\protect \citeauthoryear {%
Herbst%
\ \BBA {} Schorfheide%
}{%
Herbst%
\ \BBA {} Schorfheide%
}{%
{\protect \APACyear {2016}}%
}]{%
hs}
\APACinsertmetastar {%
hs}%
\begin{APACrefauthors}%
Herbst, E\BPBI P.%
\BCBT {}\ \BBA {} Schorfheide, F.%
\end{APACrefauthors}%
\unskip\
\newblock
\APACrefYear{2016}.
\newblock
\APACrefbtitle {Bayesian Estimation of DSGE Models} {Bayesian estimation of dsge models}.
\newblock
\APACaddressPublisher{Princeton, USA}{Princeton University Press}.
\PrintBackRefs{\CurrentBib}

\bibitem [\protect \citeauthoryear {%
Horn%
\ \BBA {} Johnson%
}{%
Horn%
\ \BBA {} Johnson%
}{%
{\protect \APACyear {1985}}%
}]{%
hj1}
\APACinsertmetastar {%
hj1}%
\begin{APACrefauthors}%
Horn, R\BPBI A.%
\BCBT {}\ \BBA {} Johnson, C\BPBI R.%
\end{APACrefauthors}%
\unskip\
\newblock
\APACrefYear{1985}.
\newblock
\APACrefbtitle {{Matrix Analysis}} {{Matrix Analysis}}.
\newblock
\APACaddressPublisher{Cambridge, United Kingdom}{Cambridge University Press}.
\PrintBackRefs{\CurrentBib}

\bibitem [\protect \citeauthoryear {%
Horn%
\ \BBA {} Johnson%
}{%
Horn%
\ \BBA {} Johnson%
}{%
{\protect \APACyear {1991}}%
}]{%
hj2}
\APACinsertmetastar {%
hj2}%
\begin{APACrefauthors}%
Horn, R\BPBI A.%
\BCBT {}\ \BBA {} Johnson, C\BPBI R.%
\end{APACrefauthors}%
\unskip\
\newblock
\APACrefYear{1991}.
\newblock
\APACrefbtitle {Topics in Matrix Analysis} {Topics in matrix analysis}.
\newblock
\APACaddressPublisher{Cambridge, United Kingdom}{Cambridge University Press}.
\PrintBackRefs{\CurrentBib}

\bibitem [\protect \citeauthoryear {%
Hsiao%
}{%
Hsiao%
}{%
{\protect \APACyear {1983}}%
}]{%
hsiao}
\APACinsertmetastar {%
hsiao}%
\begin{APACrefauthors}%
Hsiao, C.%
\end{APACrefauthors}%
\unskip\
\newblock
\APACrefYearMonthDay{1983}{September}{}.
\newblock
{\BBOQ}\APACrefatitle {{Identification}} {{Identification}}.{\BBCQ}
\newblock
\BIn{} Z.~Griliches\ \BBA {} M\BPBI D.~Intriligator\ (\BEDS), \APACrefbtitle {{Handbook of Econometrics}} {{Handbook of Econometrics}}\ (\BVOL~1, \BPG~223-283).
\newblock
\APACaddressPublisher{}{Elsevier}.
\PrintBackRefs{\CurrentBib}

\bibitem [\protect \citeauthoryear {%
Iskrev%
}{%
Iskrev%
}{%
{\protect \APACyear {2010}}%
}]{%
iskrev}
\APACinsertmetastar {%
iskrev}%
\begin{APACrefauthors}%
Iskrev, N.%
\end{APACrefauthors}%
\unskip\
\newblock
\APACrefYearMonthDay{2010}{March}{}.
\newblock
{\BBOQ}\APACrefatitle {{Local Identification in DSGE Models}} {{Local Identification in DSGE Models}}.{\BBCQ}
\newblock
\APACjournalVolNumPages{Journal of Monetary Economics}{57}{2}{189-202}.
\PrintBackRefs{\CurrentBib}

\bibitem [\protect \citeauthoryear {%
Kailath%
}{%
Kailath%
}{%
{\protect \APACyear {1980}}%
}]{%
kailath}
\APACinsertmetastar {%
kailath}%
\begin{APACrefauthors}%
Kailath, T.%
\end{APACrefauthors}%
\unskip\
\newblock
\APACrefYear{1980}.
\newblock
\APACrefbtitle {Linear Systems} {Linear systems}.
\newblock
\APACaddressPublisher{Englewood Cliffs, NJ}{Prentice-Hall}.
\PrintBackRefs{\CurrentBib}

\bibitem [\protect \citeauthoryear {%
Koci\c{e}cki%
\ \BBA {} Kolasa%
}{%
Koci\c{e}cki%
\ \BBA {} Kolasa%
}{%
{\protect \APACyear {2018}}%
}]{%
kk}
\APACinsertmetastar {%
kk}%
\begin{APACrefauthors}%
Koci\c{e}cki, A.%
\BCBT {}\ \BBA {} Kolasa, M.%
\end{APACrefauthors}%
\unskip\
\newblock
\APACrefYearMonthDay{2018}{}{}.
\newblock
{\BBOQ}\APACrefatitle {{Global Identification of Linearized DSGE Models}} {{Global Identification of Linearized DSGE Models}}.{\BBCQ}
\newblock
\APACjournalVolNumPages{Quantitative Economics}{9}{3}{1243-1263}.
\PrintBackRefs{\CurrentBib}

\bibitem [\protect \citeauthoryear {%
Koci\c{e}cki%
\ \BBA {} Kolasa%
}{%
Koci\c{e}cki%
\ \BBA {} Kolasa%
}{%
{\protect \APACyear {2023}}%
}]{%
kk23}
\APACinsertmetastar {%
kk23}%
\begin{APACrefauthors}%
Koci\c{e}cki, A.%
\BCBT {}\ \BBA {} Kolasa, M.%
\end{APACrefauthors}%
\unskip\
\newblock
\APACrefYearMonthDay{2023}{}{}.
\newblock
{\BBOQ}\APACrefatitle {A solution to the global identification problem in DSGE models} {A solution to the global identification problem in dsge models}.{\BBCQ}
\newblock
\APACjournalVolNumPages{Journal of Econometrics}{236}{2}{105477}.
\PrintBackRefs{\CurrentBib}

\bibitem [\protect \citeauthoryear {%
Komunjer%
\ \BBA {} Ng%
}{%
Komunjer%
\ \BBA {} Ng%
}{%
{\protect \APACyear {2011}}%
}]{%
kn}
\APACinsertmetastar {%
kn}%
\begin{APACrefauthors}%
Komunjer, I.%
\BCBT {}\ \BBA {} Ng, S.%
\end{APACrefauthors}%
\unskip\
\newblock
\APACrefYearMonthDay{2011}{November}{}.
\newblock
{\BBOQ}\APACrefatitle {{Dynamic Identification of Dynamic Stochastic General Equilibrium Models}} {{Dynamic Identification of Dynamic Stochastic General Equilibrium Models}}.{\BBCQ}
\newblock
\APACjournalVolNumPages{Econometrica}{79}{6}{1995-2032}.
\newblock
\begin{APACrefURL} \url{https://ideas.repec.org/a/ecm/emetrp/v79y2011i6p1995-2032.html} \end{APACrefURL}
\PrintBackRefs{\CurrentBib}

\bibitem [\protect \citeauthoryear {%
Lubik%
\ \BBA {} Schorfheide%
}{%
Lubik%
\ \BBA {} Schorfheide%
}{%
{\protect \APACyear {2003}}%
}]{%
ls}
\APACinsertmetastar {%
ls}%
\begin{APACrefauthors}%
Lubik, T\BPBI A.%
\BCBT {}\ \BBA {} Schorfheide, F.%
\end{APACrefauthors}%
\unskip\
\newblock
\APACrefYearMonthDay{2003}{November}{}.
\newblock
{\BBOQ}\APACrefatitle {{Computing sunspot equilibria in linear rational expectations models}} {{Computing sunspot equilibria in linear rational expectations models}}.{\BBCQ}
\newblock
\APACjournalVolNumPages{Journal of Economic Dynamics and Control}{28}{2}{273-285}.
\newblock
\begin{APACrefURL} \url{http://ideas.repec.org/a/eee/dyncon/v28y2003i2p273-285.html} \end{APACrefURL}
\PrintBackRefs{\CurrentBib}

\bibitem [\protect \citeauthoryear {%
L\"{u}tkepohl%
}{%
L\"{u}tkepohl%
}{%
{\protect \APACyear {2005}}%
}]{%
lutkepohl}
\APACinsertmetastar {%
lutkepohl}%
\begin{APACrefauthors}%
L\"{u}tkepohl, H.%
\end{APACrefauthors}%
\unskip\
\newblock
\APACrefYear{2005}.
\newblock
\APACrefbtitle {New Introduction to Multiple Time Series Analysis} {New introduction to multiple time series analysis}.
\newblock
\APACaddressPublisher{Berlin, Germany}{Springer}.
\PrintBackRefs{\CurrentBib}

\bibitem [\protect \citeauthoryear {%
Mavroeidis%
}{%
Mavroeidis%
}{%
{\protect \APACyear {2004}}%
}]{%
mavroeidis04}
\APACinsertmetastar {%
mavroeidis04}%
\begin{APACrefauthors}%
Mavroeidis, S.%
\end{APACrefauthors}%
\unskip\
\newblock
\APACrefYearMonthDay{2004}{}{}.
\newblock
{\BBOQ}\APACrefatitle {Weak Identification of Forward-looking Models in Monetary Economics} {Weak identification of forward-looking models in monetary economics}.{\BBCQ}
\newblock
\APACjournalVolNumPages{Oxford Bulletin of Economics and Statistics}{66}{s1}{609-635}.
\PrintBackRefs{\CurrentBib}

\bibitem [\protect \citeauthoryear {%
Mavroeidis%
}{%
Mavroeidis%
}{%
{\protect \APACyear {2005}}%
}]{%
mavroeidis05}
\APACinsertmetastar {%
mavroeidis05}%
\begin{APACrefauthors}%
Mavroeidis, S.%
\end{APACrefauthors}%
\unskip\
\newblock
\APACrefYearMonthDay{2005}{}{}.
\newblock
{\BBOQ}\APACrefatitle {Identification Issues in Forward-Looking Models Estimated by GMM, with an Application to the Phillips Curve} {Identification issues in forward-looking models estimated by gmm, with an application to the phillips curve}.{\BBCQ}
\newblock
\APACjournalVolNumPages{Journal of Money, Credit and Banking}{37}{3}{421--448}.
\PrintBackRefs{\CurrentBib}

\bibitem [\protect \citeauthoryear {%
McDonald%
\ \BBA {} Shalizi%
}{%
McDonald%
\ \BBA {} Shalizi%
}{%
{\protect \APACyear {2022}}%
}]{%
mcdonaldshalizi}
\APACinsertmetastar {%
mcdonaldshalizi}%
\begin{APACrefauthors}%
McDonald, D\BPBI J.%
\BCBT {}\ \BBA {} Shalizi, C\BPBI R.%
\end{APACrefauthors}%
\unskip\
\newblock
\APACrefYearMonthDay{2022}{}{}.
\newblock
\APACrefbtitle {Empirical Macroeconomics and DSGE Modeling in Statistical Perspective.} {Empirical macroeconomics and dsge modeling in statistical perspective.}
\newblock
\begin{APACrefURL} \url{https://arxiv.org/abs/2210.16224} \end{APACrefURL}
\PrintBackRefs{\CurrentBib}

\bibitem [\protect \citeauthoryear {%
Muth%
}{%
Muth%
}{%
{\protect \APACyear {1981}}%
}]{%
muthestimation}
\APACinsertmetastar {%
muthestimation}%
\begin{APACrefauthors}%
Muth, J\BPBI F.%
\end{APACrefauthors}%
\unskip\
\newblock
\APACrefYearMonthDay{1981}{}{}.
\newblock
{\BBOQ}\APACrefatitle {{Estimation of Economic Relationships Containing Latent Expectations Variables}} {{Estimation of Economic Relationships Containing Latent Expectations Variables}}.{\BBCQ}
\newblock
\BIn{} R\BPBI E.~Lucas\ \BBA {} T\BPBI J.~Sargent\ (\BEDS), \APACrefbtitle {{Rational Expectations and Econometric Practice}} {{Rational Expectations and Econometric Practice}}\ (\BPGS\ 321--328).
\newblock
\APACaddressPublisher{Minneapolis}{The University of Minnesota Press}.
\PrintBackRefs{\CurrentBib}

\bibitem [\protect \citeauthoryear {%
Nikolski%
}{%
Nikolski%
}{%
{\protect \APACyear {2002}}%
}]{%
nikolski}
\APACinsertmetastar {%
nikolski}%
\begin{APACrefauthors}%
Nikolski, N\BPBI K.%
\end{APACrefauthors}%
\unskip\
\newblock
\APACrefYear{2002}.
\newblock
\APACrefbtitle {Operators, Functions, and Systems: An Easy Reading: Volume 1: Hardy, Hankel, and Toeplitz} {Operators, functions, and systems: An easy reading: Volume 1: Hardy, hankel, and toeplitz}\ (\BVOL~92).
\newblock
\APACaddressPublisher{Providence, RI, USA}{American Mathematical Society}.
\PrintBackRefs{\CurrentBib}

\bibitem [\protect \citeauthoryear {%
Onatski%
}{%
Onatski%
}{%
{\protect \APACyear {2006}}%
}]{%
onatski}
\APACinsertmetastar {%
onatski}%
\begin{APACrefauthors}%
Onatski, A.%
\end{APACrefauthors}%
\unskip\
\newblock
\APACrefYearMonthDay{2006}{February}{}.
\newblock
{\BBOQ}\APACrefatitle {{Winding Number Criterion for Existence and Uniqueness of Equilibrium in Linear Rational Expectations Models}} {{Winding Number Criterion for Existence and Uniqueness of Equilibrium in Linear Rational Expectations Models}}.{\BBCQ}
\newblock
\APACjournalVolNumPages{Journal of Economic Dynamics and Control}{30}{2}{323-345}.
\newblock
\begin{APACrefURL} \url{http://ideas.repec.org/a/eee/dyncon/v30y2006i2p323-345.html} \end{APACrefURL}
\PrintBackRefs{\CurrentBib}

\bibitem [\protect \citeauthoryear {%
M.~Pesaran%
}{%
M.~Pesaran%
}{%
{\protect \APACyear {1981}}%
}]{%
pesaran81}
\APACinsertmetastar {%
pesaran81}%
\begin{APACrefauthors}%
Pesaran, M.%
\end{APACrefauthors}%
\unskip\
\newblock
\APACrefYearMonthDay{1981}{}{}.
\newblock
{\BBOQ}\APACrefatitle {{Identification of Rational Expectations Models}} {{Identification of Rational Expectations Models}}.{\BBCQ}
\newblock
\APACjournalVolNumPages{Journal of Econometrics}{16}{3}{375 - 398}.
\PrintBackRefs{\CurrentBib}

\bibitem [\protect \citeauthoryear {%
M\BPBI H.~Pesaran%
}{%
M\BPBI H.~Pesaran%
}{%
{\protect \APACyear {1987}}%
}]{%
pesaran1987}
\APACinsertmetastar {%
pesaran1987}%
\begin{APACrefauthors}%
Pesaran, M\BPBI H.%
\end{APACrefauthors}%
\unskip\
\newblock
\APACrefYear{1987}.
\newblock
\APACrefbtitle {The Limits to Rational Expectations} {The limits to rational expectations}.
\newblock
\APACaddressPublisher{New York, NY. USA}{Basil Blackwell Inc.}
\PrintBackRefs{\CurrentBib}

\bibitem [\protect \citeauthoryear {%
M\BPBI H.~Pesaran%
\ \BBA {} Smith%
}{%
M\BPBI H.~Pesaran%
\ \BBA {} Smith%
}{%
{\protect \APACyear {2011}}%
}]{%
pesaransmith}
\APACinsertmetastar {%
pesaransmith}%
\begin{APACrefauthors}%
Pesaran, M\BPBI H.%
\BCBT {}\ \BBA {} Smith, R\BPBI P.%
\end{APACrefauthors}%
\unskip\
\newblock
\APACrefYearMonthDay{2011}{}{}.
\newblock
{\BBOQ}\APACrefatitle {{Beyond the DSGE Straitjacket}} {{Beyond the DSGE Straitjacket}}.{\BBCQ}
\newblock
\APACjournalVolNumPages{The Manchester School}{79}{s2}{5-16}.
\PrintBackRefs{\CurrentBib}

\bibitem [\protect \citeauthoryear {%
Qu%
\ \BBA {} Tkachenko%
}{%
Qu%
\ \BBA {} Tkachenko%
}{%
{\protect \APACyear {2017}}%
}]{%
qt}
\APACinsertmetastar {%
qt}%
\begin{APACrefauthors}%
Qu, Z.%
\BCBT {}\ \BBA {} Tkachenko, D.%
\end{APACrefauthors}%
\unskip\
\newblock
\APACrefYearMonthDay{2017}{}{}.
\newblock
{\BBOQ}\APACrefatitle {{Global Identification in DSGE Models Allowing for Indeterminacy}} {{Global Identification in DSGE Models Allowing for Indeterminacy}}.{\BBCQ}
\newblock
\APACjournalVolNumPages{The Review of Economic Studies}{84}{3}{1306-1345}.
\PrintBackRefs{\CurrentBib}

\bibitem [\protect \citeauthoryear {%
Rayner%
}{%
Rayner%
}{%
{\protect \APACyear {1980}}%
}]{%
rayner}
\APACinsertmetastar {%
rayner}%
\begin{APACrefauthors}%
Rayner, J\BPBI R.%
\end{APACrefauthors}%
\unskip\
\newblock
\APACrefYear{1980}.
\unskip\
\newblock
\APACrefbtitle {{Expectations: Some Econometric Aspects}} {{Expectations: Some Econometric Aspects}}\ \APACtypeAddressSchool {\BUPhD}{}{}.
\unskip\
\newblock
\APACaddressSchool {}{Australian National University}.
\PrintBackRefs{\CurrentBib}

\bibitem [\protect \citeauthoryear {%
Romer%
}{%
Romer%
}{%
{\protect \APACyear {2016}}%
}]{%
romer}
\APACinsertmetastar {%
romer}%
\begin{APACrefauthors}%
Romer, P.%
\end{APACrefauthors}%
\unskip\
\newblock
\APACrefYearMonthDay{2016}{}{}.
\newblock
{\BBOQ}\APACrefatitle {{The Trouble with Macroeconomics}} {{The Trouble with Macroeconomics}}.{\BBCQ}
\newblock
\APACjournalVolNumPages{The American Economist}{20}{}{1--20}.
\PrintBackRefs{\CurrentBib}

\bibitem [\protect \citeauthoryear {%
Rosenblatt%
}{%
Rosenblatt%
}{%
{\protect \APACyear {2000}}%
}]{%
rosenblatt}
\APACinsertmetastar {%
rosenblatt}%
\begin{APACrefauthors}%
Rosenblatt, M.%
\end{APACrefauthors}%
\unskip\
\newblock
\APACrefYear{2000}.
\newblock
\APACrefbtitle {Gaussian and Non-Gaussian Linear Time Series and Random Fields} {Gaussian and non-gaussian linear time series and random fields}.
\newblock
\APACaddressPublisher{New York, USA}{Springer Verlag}.
\PrintBackRefs{\CurrentBib}

\bibitem [\protect \citeauthoryear {%
Rozanov%
}{%
Rozanov%
}{%
{\protect \APACyear {1967}}%
}]{%
rozanov}
\APACinsertmetastar {%
rozanov}%
\begin{APACrefauthors}%
Rozanov, Y\BPBI A.%
\end{APACrefauthors}%
\unskip\
\newblock
\APACrefYear{1967}.
\newblock
\APACrefbtitle {Stationary Random Processes} {Stationary random processes}.
\newblock
\APACaddressPublisher{San Francisco}{Holden-Day}.
\PrintBackRefs{\CurrentBib}

\bibitem [\protect \citeauthoryear {%
Sargan%
}{%
Sargan%
}{%
{\protect \APACyear {1988}}%
}]{%
sargan}
\APACinsertmetastar {%
sargan}%
\begin{APACrefauthors}%
Sargan, J\BPBI D.%
\end{APACrefauthors}%
\unskip\
\newblock
\APACrefYearMonthDay{1988}{}{}.
\newblock
{\BBOQ}\APACrefatitle {{The Identification and Estimation of Sets of Simultaneous Stochastic Equations}} {{The Identification and Estimation of Sets of Simultaneous Stochastic Equations}}.{\BBCQ}
\newblock
\BIn{} E.~Maasoumi\ (\BED), \APACrefbtitle {{Contributions to Econometrics}} {{Contributions to Econometrics}}\ (\BVOL~1, \BPGS\ 236--249).
\newblock
\APACaddressPublisher{Cambridge, UK}{Cambridge University Press}.
\PrintBackRefs{\CurrentBib}

\bibitem [\protect \citeauthoryear {%
Shiller%
}{%
Shiller%
}{%
{\protect \APACyear {1978}}%
}]{%
shiller}
\APACinsertmetastar {%
shiller}%
\begin{APACrefauthors}%
Shiller, R.%
\end{APACrefauthors}%
\unskip\
\newblock
\APACrefYearMonthDay{1978}{January}{}.
\newblock
{\BBOQ}\APACrefatitle {Rational Expectations and the Dynamic Structure of Macroeconomic Models: A Critical Review} {Rational expectations and the dynamic structure of macroeconomic models: A critical review}.{\BBCQ}
\newblock
\APACjournalVolNumPages{Journal of Monetary Economics}{}{4}{1-44}.
\PrintBackRefs{\CurrentBib}

\bibitem [\protect \citeauthoryear {%
Sims%
}{%
Sims%
}{%
{\protect \APACyear {1980}}%
}]{%
macroreality}
\APACinsertmetastar {%
macroreality}%
\begin{APACrefauthors}%
Sims, C\BPBI A.%
\end{APACrefauthors}%
\unskip\
\newblock
\APACrefYearMonthDay{1980}{}{}.
\newblock
{\BBOQ}\APACrefatitle {{Macroeconomics and Reality}} {{Macroeconomics and Reality}}.{\BBCQ}
\newblock
\APACjournalVolNumPages{Econometrica}{48}{1}{1-48}.
\PrintBackRefs{\CurrentBib}

\bibitem [\protect \citeauthoryear {%
Wallis%
}{%
Wallis%
}{%
{\protect \APACyear {1980}}%
}]{%
wallis80}
\APACinsertmetastar {%
wallis80}%
\begin{APACrefauthors}%
Wallis, K\BPBI F.%
\end{APACrefauthors}%
\unskip\
\newblock
\APACrefYearMonthDay{1980}{}{}.
\newblock
{\BBOQ}\APACrefatitle {Econometric Implications of the Rational Expectations Hypothesis} {Econometric implications of the rational expectations hypothesis}.{\BBCQ}
\newblock
\APACjournalVolNumPages{Econometrica}{48}{1}{49--73}.
\PrintBackRefs{\CurrentBib}

\bibitem [\protect \citeauthoryear {%
Wickens%
}{%
Wickens%
}{%
{\protect \APACyear {1982}}%
}]{%
wickens}
\APACinsertmetastar {%
wickens}%
\begin{APACrefauthors}%
Wickens, M\BPBI R.%
\end{APACrefauthors}%
\unskip\
\newblock
\APACrefYearMonthDay{1982}{}{}.
\newblock
{\BBOQ}\APACrefatitle {The Efficient Estimation of Econometric Models with Rational Expectations} {The efficient estimation of econometric models with rational expectations}.{\BBCQ}
\newblock
\APACjournalVolNumPages{The Review of Economic Studies}{49}{1}{55--67}.
\PrintBackRefs{\CurrentBib}

\bibitem [\protect \citeauthoryear {%
Zadrozny%
}{%
Zadrozny%
}{%
{\protect \APACyear {2022}}%
}]{%
zadrozny}
\APACinsertmetastar {%
zadrozny}%
\begin{APACrefauthors}%
Zadrozny, P\BPBI A.%
\end{APACrefauthors}%
\unskip\
\newblock
\APACrefYearMonthDay{2022}{October}{}.
\newblock
\APACrefbtitle {Linear Identification Of Linear Rational-Expectations Models By Exogenous Variables Reconciles Lucas And Sims} {Linear identification of linear rational-expectations models by exogenous variables reconciles lucas and sims}\ \APACbVolEdTR {}{Working Paper\ \BNUM~682}.
\newblock
\APACaddressInstitution{}{Center for Financial Studies}.
\PrintBackRefs{\CurrentBib}

\end{thebibliography}

\end{document}